\documentclass[global,envcountsect,envcountsame]{svjour}
\journalname{AAECC}
\usepackage{times,url}
\usepackage{amsmath,amssymb,bbm}
\usepackage{pstricks}
    \newgray{gray1}{0.8}
    \newgray{gray2}{0.5}
\usepackage{ifthen,float,algorithm}
\usepackage{algorithmic}
\newcommand\algorithmicfont[1]{\textbf{#1}}
\newcommand{\algorithmicreturn}{\algorithmicfont{return}\ }
\smartqed

\def\B{\mathcal{B}}
\def\bull{\hfill$\lhd$}
\def\CC{\mathbbm{C}}
\def\C{\mathcal{C}}
\def\cone#1#2{\C_{#1}(#2)}
\def\D{\mathcal{D}}
\def\ev{\vec{e}}
\def\F{\mathcal{F}}
\def\FF{\mathbbm{F}}
\def\gf{\mathfrak{g}}
\def\gr#1#2{\mathrm{gr}_{#1}{#2}}
\def\G{\mathcal{G}}
\def\H{\mathcal{H}}
\def\I{\mathcal{I}}
\def\id#1{\mathrm{id}_{#1}}
\def\idiv#1{\,|_{#1}\,}
\def\ispan#1#2{\lspan{#1}_{#2}}
\def\ispanl#1#2{\lspanl{#1}_{#2}}
\def\ispanr#1#2{\lspanr{#1}_{#2}}

\def\kk{{\mathbbm k}}
\def\L{\mathcal{L}}
\def\lc#1#2{\lco_{#1}{#2}}
\def\le#1#2{\leo_{#1}{#2}}
\def\lm#1#2{\lmo_{#1}{#2}}
\def\locp{\P_\prec}

\def\lspan#1{\langle{#1}\rangle}
\def\lspanl#1{\langle{#1}\rangle^{(l)}}
\def\lspanr#1{\langle{#1}\rangle^{(r)}}
\def\lspant#1{\langle\!\langle{#1}\rangle\!\rangle}
\def\lt#1#2{\lto_{#1}{#2}}
\def\M{\mathcal{M}}
\def\MuPAD{{\sffamily\itshape MuPAD\/}}
\def\mult#1#2#3{{#1}_{#2}(#3)}
\def\NF#1#2{\mathrm{NF}_{#1}(#2)}
\def\NN{{\mathbbm N}}
\def\N{\mathcal{N}}
\def\Nn{\NN_0^n}
\def\Nno{\NN_0^{n+1}}
\def\nmult#1#2#3{\bar{#1}_{#2}(#3)}

\def\P{\mathcal{P}}
\def\QQ{\mathbbm{Q}}
\def\R{\mathcal{R}}
\def\Rf{\mathfrak{R}}
\def\RR{\mathbbm{R}}
\def\S{\mathcal{S}}
\def\Sp{\S_\prec}

\def\Sv{\vec{S}}
\def\syz#1#2{\mathrm{Syz}^{#1}(#2)}
\def\T{\mathcal{T}}
\def\TT{\mathbbm{T}}
\def\Uf{\mathfrak{U}}

\def\WW{\mathbbm{W}}

\DeclareMathOperator{\cls}{cls}

\DeclareMathOperator{\ec}{\mathaccent"7013ecart}
\DeclareMathOperator{\lcm}{lcm}
\DeclareMathOperator{\lco}{lc}
\DeclareMathOperator{\leo}{le}
\DeclareMathOperator{\lmo}{lm}
\DeclareMathOperator{\lto}{lt}
\DeclareMathOperator{\supp}{supp}

\begin{document}
\title{A Combinatorial Approach to Involution and $\delta$-Regularity I:
  Involutive Bases in Polynomial Algebras of Solvable Type}
\titlerunning{Involution and $\delta$-Regularity I}
\author{Werner M.  Seiler}
\institute{AG ``Computational Mathematics'',
  Universit\"at Kassel, 34132 Kassel, Germany\\   
  \url{www.mathematik.uni-kassel.de/~seiler}\\ 
  \email{seiler@mathematik.uni-kassel.de}} 
\date{Received: date / Revised version: date}
\maketitle
\begin{abstract}
  Involutive bases are a special form of non-reduced Gr\"obner bases with
  additional combinatorial properties.  Their origin lies in the Janet-Riquier
  theory of linear systems of partial differential equations.  We study them
  for a rather general class of polynomial algebras including also
  non-commutative algebras like those generated by linear differential and
  difference operators or universal enveloping algebras of
  (finite-dimensional) Lie algebras.  We review their basic properties using
  the novel concept of a weak involutive basis and concrete algorithms for
  their construction.  As new original results, we develop a theory for
  involutive bases with respect to semigroup orders (as they appear in local
  computations) and over coefficient rings, respectively.  In both cases it
  turns out that generally only weak involutive bases exist.
\end{abstract}

\section{Introduction}

In the late $19$th and early $20$th century a number of French mathematicians
developed what is nowadays called the Janet-Riquier theory of differential
equations \cite{ja:edp,ja:mfa,ja:lec,mr:conv,riq:edp,th:ds,at:inv}.  It is a
theory for general systems of differential equations, i.\,e.\ also for under-
and overdetermined systems, and provides in particular a concrete algorithm
for the completion to a so-called passive\footnote{Sometimes the equivalent
  term ``involutive'' is used which seems to go back to Lie.} system.  In
recent times, interest in the theory has been rekindled mainly in the context
of Lie symmetry analysis, so that a number of references to modern works and
implementations are contained in the review \cite{her:sym}.

The defining property of passive systems is that they do not generate any
non-trivial integrability conditions.  As the precise definition of passivity
requires the introduction of a ranking on the set of all derivatives and as
every linear system of partial differential equations with constant
coefficients bijectively corresponds to a polynomial module, it appears
natural to relate this theory to the algebraic theory of Gr\"obner bases
\cite{al:gb,bw:groe}.

Essentially, the Janet-Riquier theory in its original form lacks only the
concept of reduction to a normal form; otherwise it contains all the
ingredients of Gr\"obner bases.  Somewhat surprisingly, a rigorous links has
been established only fairly recently first by Wu \cite{wu:invbas} and then by
Gerdt and collaborators who introduced a special form of non-reduced Gr\"obner
bases for polynomial ideals~\cite{gb:invbas,gb:minbas,zb:inv}, the
\emph{involutive bases} (Wu's ``well-behaved bases'' correspond to Thomas
bases in the language of \cite{gb:invbas}).  A slightly different approach to
involutive bases has been proposed by Apel \cite{apel:hilbert}; it will not be
used here.

The fundamental idea behind involutive bases (originating in the pioneering
work of Janet \cite{ja:edp,ja:mfa}) is to assign to each generator in a basis
a subset of all variables: its multiplicative variables.  This assignment is
called an involutive division, as it corresponds to a restriction of the usual
divisibility relation of terms.  We only permit to multiply each generator by
polynomials in its multiplicative variables.  As we will see later in Part II,
for appropriately prepared bases, this approach automatically leads to
combinatorial decompositions of polynomial modules.

Like Gr\"obner bases, involutive bases can be defined in many non-commutative
algebras.  We will work with a generalisation of the polynomial algebras of
solvable type introduced by Kandry-Rodi and Weispfenning \cite{krw:ncgb}.  It
is essentially equivalent to the generalisation discussed by Kredel
\cite{hk:solvpoly} or to the $G$-algebras considered by Apel \cite{apel:diss}
and Levandovskyy \cite{vl:ade,vl:diss}.  In contrast to some of these works,
we explicitly permit that the variables act on the coefficients, so that, say,
linear differential operators with \emph{variable} coefficients form a
polynomial algebra of solvable type in our sense.  Thus our framework
automatically includes the work of Gerdt \cite{vpg:ldo} on involutive bases
for linear differential equations.

This article is the first of two parts.  It reviews the basic theory of
involutive bases; this is immediately done in the framework of polynomial
algebras of solvable type, as it appears to be the most natural setting.
Indeed, we would like to stress that in our opinion the core of the involutive
bases theory is the monomial theory (in fact, we will formulate it in the
language of multi indices or exponent vectors, i.\,e.\ in the Abelian monoid
$(\Nn,+)$, in order to avoid problems with non-commuting variables) and the
subsequent extension to polynomials requires only straightforward normal form
considerations.

While much of the presented material may already be found scattered in the
literature (though not always in the generality presented here and sometimes
with incorrect proofs), the article also contains some original material.
Compared to Gerdt and Blinkov \cite{gb:invbas}, we give an alternative
definition of involutive bases which naturally leads to the new notion of a
\emph{weak} involutive basis.  While these weak bases are insufficient for the
applications studied in Part II, they extend the applicability of the
involutive completion algorithm to situations not covered before.

The main emphasis in the literature is on optimising the simple completion
algorithm of Section~\ref{sec:polycompl} and on providing fast
implementations; as the experiments reported in \cite{gby:janbas2}
demonstrate, the results have been striking.  We will, however, ignore this
rather technical topic and instead study in Part II a number of applications
of involutive bases (mainly for the special case of Pommaret bases) in the
structure analysis of polynomial modules.  This will include in particular the
relation between involutive bases and the above mentioned combinatorial
decompositions.  Note, however, that in these applications we will restrict to
the ordinary commutative polynomial ring.

This first part is organised as follows.  The next section defines involutive
divisions and bases within the Abelian monoid $(\Nn,+)$ of multi indices.  It
also introduces the two most important divisions named after Janet and
Pommaret, respectively.  Section~\ref{sec:solvalg} introduces the here used
concept of polynomial algebras of solvable type.  As the question whether
Hilbert's Basis Theorem remains valid is non-trivial if the coefficients form
only a ring and not a field, Section~\ref{sec:hilbert} is devoted to this
problem.  The following three sections define (weak) involutive bases and give
concrete algorithms for their construction.

The next four sections study some extensions of the basic theory.
Section~\ref{sec:two} analyses the relation between left and right ideals in
polynomial algebras of solvable type and the computation of bases for
two-sided ideals; this requires only a straightforward adaption of classical
Gr\"obner basis theory.  The following three sections contain original
results.  The first two ones generalise to semigroup orders and study the use
of the Mora normal form.  Finally, Section~\ref{sec:ring} considers involutive
bases over rings.  It turns out that in these more general situations usually
only weak bases exist.

In a short appendix we fix our conventions for term orders which are inverse
to the ones found in most textbooks on Gr\"obner bases.  We also prove an
elementary property of the degree reverse lexicographic term order that makes
it particularly natural for Pommaret bases.

\section{Involutive Divisions}\label{sec:invdiv}

We study the Abelian monoid $(\Nn,+)$ with the addition defined componentwise
and call its elements \emph{multi indices}.  They may be identified in a
natural way with the vertices of an $n$-dimensional integer lattice, so that
we can easily visualise subsets of $\Nn$.  For a multi index $\nu\in\Nn$ we
introduce its \emph{cone} $\C(\nu)=\nu+\Nn$, i.\,e.\ the set of all multi
indices that can be reached from $\nu$ by adding another multi index.  We say
that $\nu$ \emph{divides} $\mu$, written $\nu\mid\mu$, if $\mu\in\C(\nu)$.
Given a finite subset $\N\subset\Nn$, we define its \emph{span} as the monoid
ideal generated by $\N$:
\begin{equation}\label{eq:multispan}
  \lspan{\N}=\bigcup_{\nu\in\N}\C(\nu)\;.
\end{equation}

The basic idea of an involutive division is to introduce a restriction of the
cone of a multi index, the involutive cone: it is only allowed to add multi
indices certain entries of which vanish.  This is equivalent to a restriction
of the above defined divisibility relation.  The final goal will be having a
\emph{disjoint} union in (\ref{eq:multispan}) by using only these involutive
cones on the right hand side.  This will naturally lead to the combinatorial
decompositions discussed in Part II.

In order to finally give the definition of an involutive division, we need one
more notation: let $N\subseteq\{1,\dots,n\}$ be an arbitrary subset of the set
of the first $n$ integers; then we write $\NN^n_N=\bigl\{\nu\in\Nn\mid\forall
j\notin N:\nu_j=0\bigr\}$ for the set of all multi indices where the only
entries who may be non-zero are those whose indices are contained in $N$.

\begin{definition}\label{def:invdiv}
  An\/ \emph{involutive division} $L$ is defined on the Abelian monoid\/
  $(\Nn,+)$, if for any finite set\/ $\N\subset\Nn$ a subset\/
  $\mult{N}{L,\N}{\nu}\subseteq\{1,\dots,n\}$ of\/ \emph{multiplicative
    indices} is associated to every multi index\/ $\nu\in\N$ such that the
  following two conditions on the\/ \emph{involutive cones}
  $\cone{L,\N}{\nu}=\nu+\NN^n_{\mult{N}{L,\N}{\nu}}$ are satisfied.
  \begin{enumerate}
  \item If there exist two elements\/ $\mu,\nu\in\N$ with\/
    $\cone{L,\N}{\mu}\cap\cone{L,\N}{\nu}\neq\emptyset$, either\/
    $\cone{L,\N}{\mu}\subseteq\cone{L,\N}{\nu}$ or\/
    $\cone{L,\N}{\nu}\subseteq\cone{L,\N}{\mu}$ holds.
  \item If\/ $\N'\subset\N$, then\/
    $\mult{N}{L,\N}{\nu}\subseteq\mult{N}{L,\N'}{\nu}$ for all\/ $\nu\in\N'$.
  \end{enumerate}
  An arbitrary multi index\/ $\mu\in\Nn$ is\/ \emph{involutively divisible} by
  $\nu\in\N$, written $\nu\idiv{L,\N}\mu$, if\/ $\mu\in\cone{L,\N}{\nu}$.
\end{definition}

Before we discuss the precise meaning of this definition and in particular of
the two conditions contained in it, we should stress the following important
point: as indicated by the notation, involutive divisibility is always defined
with respect to both an involutive division $L$ and a fixed finite set
$\N\subset\Nn$: only an element of $\N$ can be an involutive divisor.
Obviously, involutive divisibility $\nu\idiv{L,\N}\mu$ implies ordinary
divisibility $\nu\mid\mu$.

The involutive cone $\cone{L,\N}{\nu}$ of any multi index $\nu\in\N$ is a
subset of the full cone $\C(\nu)$.  We are not allowed to add arbitrary multi
indices to $\nu$ but may increase only certain entries of $\nu$ determined by
the multiplicative indices.  The first condition in the above definition says
that involutive cones can intersect only trivially: if two intersect, one must
be a subset of the other.

The \emph{non-multiplicative indices} form the complement of
$\mult{N}{L,\N}{\nu}$ in $\{1,\dots,n\}$ and are denoted by
$\nmult{N}{L,\N}{\nu}$.  If we remove some elements from the set $\N$ and
determine the multiplicative indices of the remaining elements with respect to
the subset $\N'$, we obtain in general a different result than before.  The
second condition for an involutive division says that while it may happen that
a non-multiplicative index becomes multiplicative for some $\nu\in\N'$, the
converse cannot happen.

\begin{example}\label{ex:jandiv}
  A classical involutive division is the \emph{Janet division}~$J$.  In order
  to define it, we must introduce certain subsets of the given set
  $\N\subset\Nn$:
  \begin{equation}\label{eq:jancls}
    (d_k,\dots,d_n)=
        \bigl\{\,\nu\in\N\mid\nu_i=d_i\,,\ k\leq i\leq n\,\bigr\}\;.
  \end{equation}
  The index $n$ is multiplicative for $\nu\in\N$, if
  $\nu_n=\max_{\mu\in\N}\,\{\mu_n\}$, and $k<n$ is multiplicative for
  $\nu\in(d_{k+1},\dots,d_n)$, if
  $\nu_k=\max_{\mu\in(d_{k+1},\dots,d_n)}\,\{\mu_k\}$.  
  
  Obviously, this definition depends on the ordering of the variables
  $x_1,\dots,x_n$ and we obtain variants by applying an arbitrary but fixed
  permutation $\pi\in S_n$ to the variables.  In fact, Gerdt and Blinkov
  \cite{gb:invbas} use an ``inverse'' definition, i.\,e.\ they first apply the
  permutation {\small$\begin{pmatrix} 1&2&\cdots&n\\n&n-1&\cdots&1
  \end{pmatrix}$}.  Our convention is the original
  one of Janet \cite[pp.~16--17]{ja:lec}.
  
  Gerdt et al.~\cite{gby:janbas1} designed a special data structure, the Janet
  tree, for the fast determination of Janet multiplicative indices and for a
  number of other operations useful in the construction of Janet bases
  (Blinkov \cite{yab:tree} discusses similar tree structures also for other
  divisions).  As shown in \cite{wms:geocompl}, this data structure is based
  on a special relation between the Janet division and the lexicographic term
  order (see the appendix for our non-standard conventions).  This relation
  allows us to compute very quickly the multiplicative variables of any set
  $\N$ with Algorithm~\ref{alg:janet}.  The algorithm simply runs two pointers
  over the lexicographically ordered set $\N$ and changes accordingly the set
  $\M$ of potential multiplicative indices.\bull
\end{example}

\begin{algorithm}
  \caption{Multiplicative variables for the Janet division\label{alg:janet}}
  \begin{algorithmic}[1]
    \REQUIRE finite list $\N=\{\nu^{(1)},\dots,\nu^{(k)}\}$ of pairwise
             different multi indices from $\Nn$ 
    \ENSURE list $N=\bigl\{\mult{N}{J,\N}{\nu^{(1)}},\dots,
                  \mult{N}{J,\N}{\nu^{(k)}}\bigr\}$
            of lists with multiplicative variables
    \STATE $\N\leftarrow
            \mathtt{sort}(\N,\prec_{\mbox{\scriptsize\upshape lex}})$;
           \quad $\nu\leftarrow\N[1]$
    \STATE $p_1\leftarrow n$;\quad $\M\leftarrow\{1,\dots,n\}$;\quad 
           $N[1]\leftarrow\M$
    \FOR{$j$ \algorithmicfont{from} $2$ \algorithmicfont{to} $|\N|$}
        \STATE $p_2\leftarrow\mathtt{max}\,
                             \bigl\{i\mid (\nu-\N[j])_i\neq0\bigr\}$;\quad
               $\M\leftarrow\M\setminus\{p_2\}$
        \IF{$p_1<p_2$}
            \STATE $\M\leftarrow\M\cup\{p_1,\dots,p_2-1\}$
        \ENDIF
        \STATE $N[j]\leftarrow\M$;\quad$\nu\leftarrow\N[j]$;\quad 
               $p_1\leftarrow p_2$
    \ENDFOR
    \STATE \algorithmicreturn $N$
  \end{algorithmic}
\end{algorithm}

\begin{definition}\label{def:globdiv}
  The division\/ $L$ is\/ \emph{globally defined}, if the assignment of the
  multiplicative indices is independent of the set\/~$\N$; in this case we
  write simply\/ $\mult{N}{L}{\nu}$.
\end{definition}

\begin{example}\label{ex:pomdiv}
  Another very important division is the
  \emph{Pommaret\footnote{Historically seen, the terminology ``Pommaret
      division'' is a misnomer, as this division was already introduced by
      Janet \cite[p.~30]{ja:edp}, too.  However, the name has been generally
      accepted by now, so we stick to it.} division}~$P$.  It assigns the
  multiplicative indices according to a simple rule: if $1\leq
  k\leq n$ is the smallest index such that $\nu_k>0$ for some multi index
  $\nu\in\Nn\setminus\{[0,\dots,0]\}$, then we call $k$ the \emph{class} of
  $\nu$, written $\cls{\nu}$, and set $\mult{N}{P}{\nu}=\{1,\dots,k\}$.
  Finally, we define $\mult{N}{P}{[0,\dots,0]}=\{1,\dots,n\}$.  Hence $P$ is
  globally defined.  Like the Janet division it depends on the ordering of the
  variables $x_1,\dots,x_n$ and thus one may again introduce simple variants
  by applying a permutation.
  
  Above we have seen that the Janet division is in a certain sense related to
  the inverse lexicographic order.  The Pommaret division has a special
  relation to class respecting orders (recall that according to
  Lemma~\ref{lem:revlex} any class respecting term order coincides on terms of
  the same degree with the reverse lexicographic order).  Obviously, for
  homogeneous polynomials such orders always lead to maximal sets of
  multiplicative indices and thus to smaller bases.  But we will also see in
  Part II that from a theoretical point of view Pommaret bases with respect to
  such an order are particularly useful.\bull
\end{example}

Above we introduced the span of a set $\N\subset\Nn$ as the union of the cones
of its elements.  Given an involutive division it appears natural to consider
also the union of the involutive cones.  Obviously, this yields in general
only a subset (without any algebraic structure) of the monoid ideal
$\lspan{\N}$.

\begin{definition}\label{def:invspan}
  The\/ \emph{involutive span} of a finite set\/ $\N\subset\Nn$ is
  \begin{equation}\label{eq:ispanmono}
    \ispan{\N}{L}=\bigcup_{\nu\in\N}\cone{\N,L}{\nu}\;.
  \end{equation}
  The set\/ $\N$ is \emph{weakly involutive} for the division\/ $L$ or a\/
  \emph{weak involutive basis} of the monoid ideal\/ $\lspan{\N}$, if\/
  $\ispan{\N}{L}=\lspan{\N}$.  A weak involutive basis is a\/ \emph{strong
    involutive basis} or for short an\/ \emph{involutive basis}, if the union
  on the right hand side of (\ref{eq:ispanmono}) is disjoint, i.\,e.\ the
  intersections of the involutive cones are empty.  We call any finite set\/
  $\N\subseteq\bar\N\subset\Nn$ such that\/ $\ispan{\bar\N}{L}=\lspan{\N}$ a\/
  \emph{(weak) involutive completion} of\/ $\N$.  An\/ \emph{obstruction to
    involution} for the set\/ $\N$ is a multi index\/
  $\nu\in\lspan{\N}\setminus\ispan{\N}{L}$.
\end{definition}

\begin{remark}\label{rem:autored}
  An obvious necessary condition for a strong involutive basis is that no
  distinct multi indices $\mu,\nu\in\N$ exist such that $\mu\idiv{L,\N}\nu$.
  Sets with this property are called \emph{involutively autoreduced}.  One
  easily checks that the definition of the Janet division implies that
  $\cone{\N,J}{\mu}\cap\cone{\N,L}{\nu}=\emptyset$ whenever $\mu\neq\nu$.
  Hence for this particular division any set is involutively autoreduced.\bull
\end{remark}

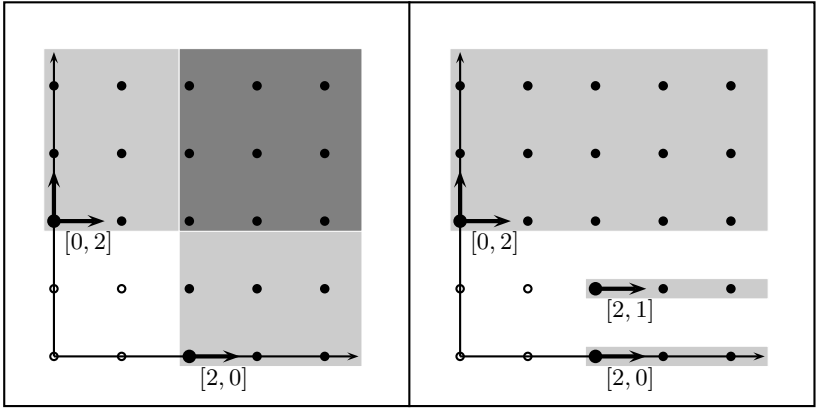
\begin{figure}[ht]
  \begin{center}
    \psset{unit=0.9cm}
    \begin{pspicture}(12,6)
      \psframe(0,0)(12,6)
      \psline(6,0)(6,6)
      \psframe[fillstyle=solid,fillcolor=gray1,linecolor=white,linewidth=0pt]%
          (0.6,2.6)(2.6,5.3)
      \psframe[fillstyle=solid,fillcolor=gray1,linecolor=white,linewidth=0pt]%
          (2.6,0.6)(5.3,2.6)
      \psframe[fillstyle=solid,fillcolor=gray2,linecolor=white,linewidth=0pt]%
          (2.6,2.6)(5.3,5.3)
      \psframe[fillstyle=solid,fillcolor=gray1,linecolor=white,linewidth=0pt]%
          (6.6,2.6)(11.3,5.3)
      \psframe[fillstyle=solid,fillcolor=gray1,linecolor=white,linewidth=0pt]%
          (8.6,0.6)(11.3,0.9)
      \psframe[fillstyle=solid,fillcolor=gray1,linecolor=white,linewidth=0pt]%
          (8.6,1.6)(11.3,1.9)
      \psline{->}(0.75,0.75)(5.25,0.75)
      \psline{->}(0.75,0.75)(0.75,5.25)
      \psline{->}(6.75,0.75)(11.25,0.75)
      \psline{->}(6.75,0.75)(6.75,5.25)
      \multips(0.75,0.75)(1,0){2}{\pscircle{0.07}}
      \multips(0.75,1.75)(1,0){2}{\pscircle{0.07}}
      \multips(2.75,0.75)(1,0){3}%
          {\pscircle[fillstyle=solid,fillcolor=black]{0.07}}
      \multips(2.75,1.75)(1,0){3}%
          {\pscircle[fillstyle=solid,fillcolor=black]{0.07}}
      \multips(0.75,2.75)(1,0){5}%
          {\pscircle[fillstyle=solid,fillcolor=black]{0.07}}
      \multips(0.75,3.75)(1,0){5}%
          {\pscircle[fillstyle=solid,fillcolor=black]{0.07}}
      \multips(0.75,4.75)(1,0){5}%
          {\pscircle[fillstyle=solid,fillcolor=black]{0.07}}
      \multips(6.75,0.75)(1,0){2}{\pscircle{0.07}}
      \multips(6.75,1.75)(1,0){2}{\pscircle{0.07}}
      \multips(8.75,0.75)(1,0){3}%
          {\pscircle[fillstyle=solid,fillcolor=black]{0.07}}
      \multips(8.75,1.75)(1,0){3}%
          {\pscircle[fillstyle=solid,fillcolor=black]{0.07}}
      \multips(6.75,2.75)(1,0){5}%
          {\pscircle[fillstyle=solid,fillcolor=black]{0.07}}
      \multips(6.75,3.75)(1,0){5}%
          {\pscircle[fillstyle=solid,fillcolor=black]{0.07}}
      \multips(6.75,4.75)(1,0){5}%
          {\pscircle[fillstyle=solid,fillcolor=black]{0.07}}
      \qdisk(0.75,2.75){0.1}
      \qdisk(2.75,0.75){0.1}
      \qdisk(6.75,2.75){0.1}
      \qdisk(8.75,0.75){0.1}
      \qdisk(8.75,1.75){0.1}
      \rput[tl]{0}(0.9,2.6){$[0,2]$}
      \rput[tl]{0}(2.9,0.6){$[2,0]$}
      \rput[tl]{0}(6.9,2.6){$[0,2]$}
      \rput[tl]{0}(8.9,0.6){$[2,0]$}
      \rput[tl]{0}(8.9,1.6){$[2,1]$}
      \psline[linewidth=1.6pt]{->}(0.75,2.75)(0.75,3.5)
      \psline[linewidth=1.6pt]{->}(0.75,2.75)(1.5,2.75)
      \psline[linewidth=1.6pt]{->}(2.75,0.75)(3.5,0.75)
      \psline[linewidth=1.6pt]{->}(6.75,2.75)(7.5,2.75)
      \psline[linewidth=1.6pt]{->}(6.75,2.75)(6.75,3.5)
      \psline[linewidth=1.6pt]{->}(8.75,0.75)(9.5,0.75)
      \psline[linewidth=1.6pt]{->}(8.75,1.75)(9.5,1.75)
    \end{pspicture}
    \caption{\emph{Left:} intersecting cones.  \emph{Right:} involutive cones.}
    \label{fig:cone}
  \end{center}
\end{figure}

\begin{example}\label{ex:inv}
  Figure~\ref{fig:cone} demonstrates the geometric interpretation of
  involutive divisions for $n=2$.  In both diagrams one can see the monoid
  ideal generated by the set $\N=\bigl\{[0,2],[2,0]\bigr\}$; the vertices
  belonging to it are marked by dark points.  The arrows represent the
  multiplicative indices, i.\,e.\ the ``allowed directions'', for both the
  Janet and the Pommaret division, as they coincide for this example.  The
  left diagram shows that the full cones of the two elements of $\N$ intersect
  in the darkly shaded area and that $\N$ is not (weakly) involutive, as the
  multi indices $[k,1]$ with $k\geq2$ are obstructions to involution.  The
  right diagram shows a strong involutive basis of $\lspan{\N}$ for both the
  Janet and the Pommaret division.  We must add to $\N$ the multi index
  $[2,1]$ and both for it and for $[2,0]$ only the index~$1$ is
  multiplicative.  One clearly sees how the span $\lspan{\N}$ is decomposed
  into three disjoint involutive cones: one of dimension~$2$, two of
  dimension~$1$.\bull
\end{example}

We are particularly interested in \emph{strong} involutive bases.  The
following result shows that in the ``monomial'' case any weak involutive basis
can be reduced to a strong one by simply eliminating some elements.

\begin{proposition}\label{prop:strongbas}
  If\/ $\N$ is a weakly involutive set, then a subset\/ $\N'\subseteq\N$
  exists such that\/ $\N'$ is a strong involutive basis of\/ $\lspan{\N}$.
\end{proposition}

\begin{proof}
  This proposition represents a nice motivation for the two conditions in
  Definition~\ref{def:invdiv} of an involutive division.  If $\N$ is not yet a
  strong involutive basis, the union in (\ref{eq:ispanmono}) is not disjoint
  and intersecting involutive cones exist.  By the first condition, this
  implies that some cones are contained in other ones; no other form of
  intersection is possible.  If we eliminate the tips of these cones
  from~$\N$, we get a subset~$\N'\subset\N$ which, by the second condition,
  has the same involutive span, as the remaining elements may only gain
  additional multiplicative indices.  Thus after a finite number of such
  eliminations we arrive at a strong involutive basis.\qed
\end{proof}

\begin{remark}\label{rem:sumbas}
  Let $\I_1$, $\I_2$ be two monoid ideals in $\Nn$ and $\N_1$, $\N_2$ (weak)
  involutive bases of them for some division $L$.  In general, we cannot
  expect that $\N_1\cup\N_2$ is again a weak involutive basis of the ideal
  $\I_1+\I_2$, as the involutive cones of the generators may shrink when taken
  with respect to the larger set $\N_1\cup\N_2$. Only for a global division we
  always obtain at least a weak involutive basis (which may then be reduced to
  a strong basis according to Proposition \ref{prop:strongbas}).\bull
\end{remark}

Recall that for arbitrary monoid ideals a basis $\N$ is called \emph{minimal},
if it is not possible to remove an element of $\N$ without losing the property
that we have a basis.  A similar notion can be naturally introduced for
involutive bases.

\begin{definition}\label{def:minbas}
  Let\/ $\I\subseteq\Nn$ be a monoid ideal and\/ $L$ an involutive division.
  An involutive basis\/ $\N$ of\/ $\I$ with respect to\/ $L$ is called\/
  \emph{minimal}, if any other involutive basis\/ $\N'$ of\/ $\I$ with respect
  to\/ $L$ satisfies\/ $\N\subseteq\N'$.
\end{definition}

Obviously, the minimal involutive basis of a monoid ideal is unique, if it
exists.  For globally defined divisions, any involutive basis is unique.

\begin{proposition}\label{prop:monomin}
  Let\/ $L$ be a globally defined division and\/ $\I\subseteq\Nn$ a monoid
  ideal.  If\/ $\I$ has an involutive basis for\/ $L$, then it is unique and
  thus minimal.
\end{proposition}

\begin{proof}
  Let $\N$ be the minimal basis of $\I$ and $\N_1$, $\N_2$ two distinct
  involutive bases of $\I$.  Both $\N_1\setminus\N_2$ and $\N_2\setminus\N_1$
  must be non-empty, as otherwise one basis was contained in the other one and
  thus the larger basis could not be involutively autoreduced with respect to
  the global division $L$.  Take an arbitrary multi index
  $\nu\in\N_1\setminus\N_2$.  The basis $\N_2$ contains a unique multi index
  $\mu$ such that $\mu\idiv{L}\nu$.  It cannot be an element of $\N_1$, as
  $\N_1$ is involutively autoreduced.  Thus $\N_1$ must contain a unique multi
  index $\lambda\neq\mu$ such that $\lambda\idiv{L}\mu$.  As $L$ is globally
  defined, this implies that $\lambda\idiv{L}\nu$ and $\lambda\neq\nu$, a
  contradiction.  \qed
\end{proof}

The algorithmic construction of (weak) involutive completions for a given set
$\N\subset\Nn$ will be discussed in detail in Section~\ref{sec:monocompl}.
For the moment we only note that we cannot expect that for an arbitrary set
$\N$ and an arbitrary involutive division $L$ an involutive basis $\N'$ of
$\lspan{\N}$ exists.

\begin{example}\label{ex:infbas}
  We consider the set $\N=\bigl\{[1,1]\bigr\}$ for the Pommaret division.  As
  $\cls{[1,1]}=1$, we get $\mult{N}{P}{[1,1]}=\{1\}$.  So
  $\cone{P}{[1,1]}\subsetneq\C([1,1])$.  But any multi index contained in
  $\lspan{\N}$ also has class $1$.  Hence no \emph{finite} involutive basis of
  $\lspan{\N}$ exists for the Pommaret division.  We can generate it
  involutively only with the infinite set $\bigl\{[1,k]\mid
  k\in\NN\bigr\}$.
  
  More generally, we may consider an \emph{irreducible} monoid ideal $\I$ in
  $\Nn$.  It is well-known that any such $\I$ has a minimal basis of the form
  $\{(\ell_1)_{i_1},\dots,(\ell_k)_{i_k}\}$ with $1\leq k\leq n$, $\ell_j>0$
  and $1\leq i_1<\cdots<i_k\leq n$.  Here $(\ell_j)_{i_j}$ is the multi index
  where all entries are zero except of the $i_j$th one which has the value
  $\ell_j$.  Such an ideal possesses a Pommaret basis, if and only if there
  are no ``gaps'' in the sequence $i_1<\cdots<i_k\leq n$, i.\,e.\ $i_k=n$ and
  $i_1=n-k+1$.  Indeed, if a gap exists, say between $i_j$ and $i_{j+1}$, then
  any Pommaret basis must contain the infinitely many multi indices of the
  form $(\ell_j)_{i_j}+(\ell)_{i_j+1}$ with $\ell>0$ and thus cannot be finite
  (obviously, in this case a simple renumbering of the variables suffices to
  remedy this problem).  Conversely, if no gaps appear, then it is easy to see
  that the set of all multi indices
  $[0,\dots,0,\ell_{i_j},\mu_{i_j+1},\dots,\mu_n]$ with $1\leq j\leq k$ and
  $0\leq\mu_i<\ell_{n-k+i}$ is a strong Pommaret basis of $\I$.  \bull
\end{example}

\begin{definition}\label{def:noether}
  An involutive division\/ $L$ is called\/ \emph{Noetherian}, if any finite
  subset\/ $\N\subset\Nn$ possesses a finite involutive completion with
  respect to $L$.
\end{definition}

\begin{lemma}\label{lem:jannoeth}
  The Janet division is Noetherian.
\end{lemma}

\begin{proof}
  Let $\N\subset\Nn$ be an arbitrary finite set.  We explicitly construct a
  Janet basis for $\lspan{\N}$.  Define the multi index $\mu=\lcm{\N}$ by
  $\mu_i=\max_{\nu\in\N}\nu_i$.  Then we claim that the (obviously finite) set
  \begin{equation}
    \bar\N=\bigl\{\bar\nu\in\lspan{\N}\mid\mu\in\C(\bar\nu)\bigr\}
  \end{equation}
  is an involutive completion of $\N$ with respect to the Janet division.
  Indeed, $\N\subseteq\bar\N$ and $\bar\N\subset\lspan{\N}$.  Let
  $\rho\in\lspan{\N}$ be an arbitrary element.  If $\rho\in\bar\N$, then
  trivially $\rho\in\ispan{\bar\N}{J}$.  Otherwise set
  $I=\{i\mid\rho_i>\mu_i\}$ and define the multi index $\bar\rho$ by
  $\bar\rho_i=\rho_i$ for $i\notin I$ and $\bar\rho_i=\mu_i$ for $i\in I$,
  i.\,e.{}\ $\bar\rho_i=\min\,\{\rho_i,\mu_i\}$.  By construction of the set
  $\bar\N$ and the definition of $\mu$, we have that $\bar\rho\in\bar\N$ and
  $I\subseteq\mult{N}{J,\bar\N}{\bar\rho}$.  But this implies that
  $\rho\in\cone{J,\bar\N}{\bar\rho}$ and thus $\bar\N$ is a finite Janet basis
  for $\lspan{\N}$.\qed
\end{proof}

\section{Polynomial Algebras of Solvable Type}\label{sec:solvalg}

We could identify multi indices and monomials and proceed to define involutive
bases for polynomial ideals.  But as the basic ideas remain unchanged in many
different situations, e.\,g.\ rings of linear differential or difference
operators, we generalise a concept originally introduced by Kandry-Rody and
Weispfenning \cite{krw:ncgb} and use \emph{polynomial algebras of solvable
  type} (the same kind of generalisation has already been intensively studied
by Kredel \cite{hk:solvpoly}).

Let $\P=\R[x_1,\dots,x_n]$ be a polynomial ring over a unitary
ring\footnote{For us a ring is always associative.} $\R$.  If $\R$ is
commutative, then $\P$ is a unitary commutative ring with respect to the usual
multiplication.  We equip the $\R$-module $\P$ with alternative
multiplications, in particular with non-commutative ones.  We allow that both
the variables~$x_i$ do not commute any more and that they operate on the
coefficients.  The usual multiplication is denoted either by a dot $\cdot$ or
by no symbol at all.  Alternative multiplications $\P\times\P\rightarrow\P$
are always written as $f\star g$.

Like Gr\"obner bases, involutive bases are always defined with respect to a
\emph{term order} $\prec$.  It selects in each polynomial $f\in\P$ a
\emph{leading term} $\lt{\prec}{f}=x^\mu$ with \emph{leading exponent}
$\le{\prec}{f}=\mu$.  The coefficient $r\in\R$ of $x^\mu$ in $f$ is the
\emph{leading coefficient} $\lc{\prec}{f}$ and the product $rx^\mu$ is the
\emph{leading monomial} $\lm{\prec}{f}$.  Based on the leading exponents we
associate to each finite set $\F\subset\P$ a set $\le{\prec}{\F}\subset\Nn$ to
which we may apply the theory developed in the previous section.  But this
requires a kind of compatibility between the multiplication $\star$ and the
chosen term order.

\begin{definition}\label{def:solvalg}
  $(\P,\star,\prec)$ is a\/ \emph{polynomial algebra of solvable type} over
  the unitary coefficient ring\/ $\R$ for the term order\/ $\prec$, if the
  multiplication\/ $\star:\P\times\P\rightarrow\P$ satisfies three axioms.
  \begin{description}
  \item[\phantom{ii}{\upshape (i)}] $(\P,\star)$ is a ring with
    unit\/ $1$.
  \item[\phantom{i}{\upshape (ii)}] $\forall r\in\R,\, f\in\P:\ r\star f=rf$.
  \item[{\upshape (iii)}] $\forall\mu,\nu\in\Nn,\, r\in\R\setminus\{0\}:\ 
    \le{\prec}{(x^\mu\star x^\nu)}=\mu+\nu\ \wedge\ \le{\prec}{(x^\mu\star
      r)}=\mu$.
  \end{description}
\end{definition}

Condition (i) ensures that arithmetics in $(\P,\star,\prec)$ obeys the usual
associative and distributive laws.  Because of Condition (ii),
$(\P,\star,\prec)$ is a left $\R$-module.  We do not require that it is a
right $\R$-module, as this would exclude the possibility that the variables
$x_i$ operate non-linearly on $\R$.  Condition (iii) ensures the compatibility
of the new multiplication $\star$ and the term order $\prec$; we say that the
multiplication $\star$ \emph{respects the term order} $\prec$.  It implies the
existence of injective maps $\rho_\mu:\R\rightarrow\R$, maps
$h_\mu:\R\rightarrow\P$ with $\le{\prec}{\bigl(h_\mu(r)\bigr)}\prec\mu$ for
all $r\in\R$, coefficients $r_{\mu\nu}\in\R\setminus\{0\}$ and polynomials
$h_{\mu\nu}\in\P$ with $\le{\prec}{h_{\mu\nu}}\prec\mu+\nu$ such that
\begin{subequations}\label{eq:crmono}
  \begin{gather}
    x^\mu\star r=\rho_\mu(r)x^\mu+h_\mu(r)\;,\label{eq:crmono1}\\
    x^\mu\star x^\nu=r_{\mu\nu}x^{\mu+\nu}+h_{\mu\nu}\;.\label{eq:crmono2}
  \end{gather}
\end{subequations}

\begin{lemma}\label{lem:rrho}
  The maps\/ $\rho_\mu$ and the coefficients\/ $r_{\mu\nu}$ satisfy for
  arbitrary multi indices\/ $\mu,\nu,\lambda\in\Nn$ and for arbitrary ring
  elements\/ $r\in\R$
  \begin{subequations}\label{eq:rrho}
    \begin{gather}
      \rho_\mu\bigl(\rho_\nu(r)\bigr)r_{\mu\nu}=
          r_{\mu\nu}\rho_{\mu+\nu}(r)\;,\label{eq:xxr}\\
      \rho_\mu(r_{\nu\lambda})r_{\mu,\nu+\lambda}=
          r_{\mu\nu}r_{\mu+\nu,\lambda}\;.\label{eq:xxx}
    \end{gather}
  \end{subequations}
  Furthermore, all maps\/ $\rho_\mu$ are ring endomorphisms.
\end{lemma}

\begin{proof}
  The first assertion is a trivial consequence of the associativity of the
  multiplication~$\star$.  The equations correspond to the leading
  coefficients of the equalities $x^\mu\star(x^\nu\star r)=(x^\mu\star
  x^\nu)\star r$ and $x^\mu\star(x^\nu\star x^\lambda)=(x^\mu\star x^\nu)\star
  x^\lambda$, respectively.  The second assertion follows mainly from
  Condition (i).\qed
\end{proof}

If $\R$ is a (skew) field, then for arbitrary polynomials $f,g\in\P$ an
element $r\in\R\setminus\{0\}$ and a polynomial $h\in\P$ satisfying
$\le{\prec}{h}\prec\le{\prec}{(f\cdot g)}$ exist such that
\begin{equation}\label{eq:crsolvalg}
  f\star g=r\,(f\cdot g)+h\;.
\end{equation}
Indeed, if $\lm{\prec}{f}=ax^\mu$ and $\lm{\prec}{g}=bx^\nu$, then a simple
computation yields that $r$ is the (unique) solution of the equation
$a\rho_\mu(b)r_{\mu\nu}=rab$ and $h$ is the difference $f\star g-r(f\cdot g)$.
Under this assumption we may reformulate (iii) as
\begin{description}
\item[{\upshape (iii)'}] $\forall f,g\in\P:\le{\prec}{(f\star
    g)}=\le{\prec}{f}+\le{\prec}{g}$.
\end{description}
For this special case the same class of non-commutative algebras was
introduced in \cite{bglc:quantum,bgv:algnoncomm} under the name \emph{PBW
  algebras} (see Example~\ref{ex:pbw} below for an explanation of this name).

\begin{proposition}\label{prop:fixmult}
  The product\/ $\star$ is fixed, as soon as the following data are given:
  constants\/ $r_{ij}\in\R\setminus\{0\}$, polynomials\/ $h_{ij}\in\P$ and
  maps\/ $\rho_i:\R\rightarrow\R$, $h_i:\R\rightarrow\P$ such that for\/
  $1\leq i\leq n$
  \begin{subequations}\label{eq:basiccr}
    \begin{gather}
      x_i\star r=\rho_i(r)x_i+h_i(r)\;,\quad\forall r\in\R\;,\label{eq:crop}\\
      x_i\star x_j=r_{ij}x_j\star x_i+h_{ij}\;,\quad\forall 1\leq j<i\;.
      \label{eq:crij}
    \end{gather}
  \end{subequations}  
\end{proposition}

\begin{proof}
  The set of all ``monomials'' $x_{i_1}\star x_{i_2}\star\cdots\star x_{i_q}$
  with $i_1\leq i_2\leq\cdots\leq i_q$ forms an $\R$-linear basis of $\P$, as
  because of (iii) the $\R$-linear map defined by $x_{i_1}\star
  x_{i_2}\star\cdots\star x_{i_q} \mapsto x_{i_1}\cdot x_{i_2}\cdots x_{i_q}$
  is an $\R$-module automorphism mapping the new basis into the standard
  basis.  Obviously, it is possible to evaluate any product $f\star g$ by
  repeated applications of the rewrite rules (\ref{eq:basiccr}) provided $f$
  and $g$ are expressed in the new basis.\qed
\end{proof}

Note that this proof is non-constructive in the sense that we are not able to
determine the multiplication in terms of the standard basis, as we do not know
explicitly the transformation between the new and the standard basis.  The
advantage of this proof is that it is valid for arbitrary coefficient rings
$\R$.  Making some assumptions on $\R$ (the simplest possibility is to require
that it is a field), one could use Lemma~\ref{lem:rrho} to express the
coefficients $r_{\mu\nu}$ and $\rho_\mu$ in (\ref{eq:crmono}) by the data in
(\ref{eq:basiccr}).  This would yield a constructive proof.

These considerations also show the main difference between our definition of
solvable algebras and related definitions that one can find at various places
in the literature.  Our Condition (iii) represents in most approaches a lemma.
Instead one usually imposes conditions like $x_i\star x_j=x_ix_j$ for $i<j$.
Then the rewrite rule (\ref{eq:crij}) suffices to obtain an explicit
transformation between the two bases used in the proof above.  This is for
instance the approach taken by Kandry-Rody and Weispfenning \cite{krw:ncgb}
and subsequently by Kredel \cite{hk:solvpoly}.

Of course, the data in Proposition \ref{prop:fixmult} cannot be chosen
arbitrarily.  Besides the obvious conditions on the leading exponents of the
polynomials $h_{ij}$ and $h_i(r)$ imposed by Condition (iii), each map
$\rho_i$ must be an injective $\R$-endomorphism and each map $h_i$ must
satisfy $h_i(r+s)=h_i(r)+h_i(s)$ and a kind of pseudo-Leibniz rule
$h_i(rs)=\rho_i(r)h_i(s)+h_i(r)\star s$.  The associativity of $\star$ imposes
further rather complicated conditions on the data.  For the case of a
$G$-algebra with the multiplication defined by rewrite rules they have been
explicitly determined by Levandovskyy \cite{vl:ade,vl:diss} who called them
\emph{non-degeneracy conditions} (see also the extensive discussion in
\cite[Sect.~3.3]{hk:solvpoly}).

\begin{example}\label{ex:ore}
  An important class of non-commutative polynomials was originally introduced
  by Noether and Schmeidler \cite{ns:ore} and later systematically studied by
  Ore \cite{oo:poly}; our exposition follows \cite{bp:pseudo}.  It includes in
  particular linear differential and difference operators (with variable
  coefficients).
  
  Let $\FF$ be an arbitrary commutative ring and $\sigma:\FF\rightarrow\FF$ an
  injective endomorphism.  A \emph{pseudo-derivation} with respect to $\sigma$
  is a map $\delta:\FF\rightarrow\FF$ such that
  (i)~$\delta(f+g)=\delta(f)+\delta(g)$ and (ii) $\delta(f\cdot
  g)=\sigma(f)\cdot\delta(g)+\delta(f)\cdot g$ for all $f,g\in\FF$.  If
  $\sigma=\id{\FF}$, the identity map, (ii) is the standard Leibniz rule for
  derivations.  If $\sigma\neq\id{\FF}$, one can show that there exists an
  $h\in\FF$ such that $\delta=h(\sigma-\id{\FF})$.  And conversely, if
  $\delta\neq0$, there exists an $h\in\FF$ such that
  $\sigma=h\delta+\id{\FF}$.  Ore called $\sigma(f)$ the \emph{conjugate} and
  $\delta(f)$ the \emph{derivative} of $f$.
  
  Given $\sigma$ and $\delta$, the ring $\FF[\partial;\sigma,\delta]$ of
  univariate \emph{Ore polynomials} consists of all formal polynomials in
  $\partial$ with coefficients in $\FF$, i.\,e.\ of expressions of the form
  $\theta=\sum_{i=0}^qf_i\partial^i$ with $f_i\in\FF$ and $q\in\NN_0$.  The
  addition is defined as usual.  The variable $\partial$ operates on an
  element $f\in\FF$ according to the rule
  \begin{equation}\label{eq:oreprod}
    \partial\star f=\sigma(f)\partial+\delta(f)
  \end{equation}
  which is extended associatively and distributively to define the
  multiplication in $\FF[\partial;\sigma,\delta]$: given two elements
  $\theta_1,\theta_2\in\FF[\partial;\sigma,\delta]$, we can transform the
  product $\theta_1\star\theta_2$ to the above normal form by repeatedly
  applying (\ref{eq:oreprod}).  The injectivity of the endomorphism $\sigma$
  ensures that $\deg{(\theta_1\star\theta_2)}=\deg{\theta_1}+\deg{\theta_2}$.
  We call $\FF[\partial;\sigma,\delta]$ the \emph{Ore extension} of $\FF$
  generated by $\sigma$ and $\delta$.
  
  A simple concrete example is given by choosing for $\FF$ some ring of
  differentiable functions in the real variable $x$, say $\FF=\QQ[x]$,
  $\delta=\frac{d}{dx}$ and $\sigma=\id{\FF}$ yielding \emph{linear ordinary
    differential operators} with polynomial functions as coefficients (i.\,e.\ 
  the Weyl algebra over $\QQ$).  Similarly, we obtain \emph{linear recurrence}
  and \emph{difference operators}.  We set $\FF=\CC(n)$, the space of
  sequences with complex elements, and take for $\sigma$ the shift operator,
  i.\,e.\ the automorphism mapping $s_n$ to $s_{n+1}$.  Then
  $\Delta=\sigma-\id{\FF}$ is a pseudo-derivation.  $\FF[E;\sigma,0]$ consists
  of linear ordinary recurrence operators, $\FF[E;\sigma,\Delta]$ of linear
  ordinary difference operators.
  
  For multivariate Ore polynomials we take a set
  $\Sigma=\{\sigma_1,\dots,\sigma_n\}$ of $\FF$-endomorphisms and a set
  $\Delta=\{\delta_1,\dots,\delta_n\}$ where each $\delta_i$ is a
  pseudo-derivation with respect to~$\sigma_i$.  For each pair
  $(\sigma_i,\delta_i)$ we introduce a variable~$\partial_i$ satisfying a
  commutation rule (\ref{eq:oreprod}).  If we require that all the maps
  $\sigma_i,\delta_j$ commute with each other, i.\,e.\ 
  $\sigma_i\circ\sigma_j=\sigma_j\circ\sigma_i$,
  $\delta_i\circ\delta_j=\delta_j\circ\delta_i$ and
  $\sigma_i\circ\delta_j=\delta_j\circ\sigma_i$ for all $i\neq j$, one easily
  checks that $\partial_i\star\partial_j=\partial_j\star\partial_i$, i.\,e.\ 
  the variables $\partial_i$ commute.  Setting
  $\D=\{\partial_1,\dots,\partial_n\}$, we denote by $\FF[\D;\Sigma,\Delta]$
  the ring of multivariate Ore polynomials.  Because of the commutativity of
  the $\partial_i$ we may write the terms as $\partial^\mu$ with multi indices
  $\mu\in\Nn$, so that it indeed makes sense to speak of a polynomial ring.
  Comparing with Proposition~\ref{prop:fixmult}, we see that we are in the
  special case where the maps $h_i$ always yield constant polynomials and the
  variables $x_i$ commute.
  
  Finally, we show that $\bigl(\FF[\D;\Sigma,\Delta],\star,\prec\bigr)$ is an
  algebra of solvable type for any term order $\prec$.  The product of two
  monomial operators $a\partial^\mu$ and $b\partial^\nu$ is given by
  \begin{equation}\label{eq:multdiffop}
    a\partial^\mu\star b\partial^\nu=
        \sum_{\lambda+\kappa=\mu}\binom{\mu}{\lambda}a
            \sigma^\lambda\bigl(\delta^\kappa(b)\bigr)\partial^{\lambda+\nu}
  \end{equation}
  where $\binom{\mu}{\lambda}$ is a shorthand for
  $\prod_{i=1}^n\binom{\mu_i}{\lambda_i}$,
  $\sigma^\lambda=\sigma_1^{\lambda_1}\circ\cdots\circ\sigma_n^{\lambda_n}$
  and similarly for $\delta$.  By the properties of a term order this implies
  \begin{equation}
    \le{\prec}{\bigl(a\partial^\mu\star b\partial^\nu\bigr)}=\mu+\nu=
        \le{\prec}{\bigl(a\partial^\mu\bigr)}+
        \le{\prec}{\bigl(b\partial^\nu\bigr)}\;,
  \end{equation}
  as any term $\partial^{\lambda+\nu}$ appearing on the right hand side of
  (\ref{eq:multdiffop}) divides $\partial^{\mu+\nu}$ and thus
  $\partial^{\lambda+\nu}\preceq\partial^{\mu+\nu}$ for any term order
  $\prec$.\bull
\end{example}

\begin{example}\label{ex:pbw}
  Bell and Goodearl \cite{bg:pbw} introduced the
  \emph{Poincar\'e-Birkhoff-Witt extension} (for short \emph{PBW extension})
  of a ring $\R$ as a ring $\P\supseteq\R$ containing a finite number of
  elements $x_1,\dots,x_n\in\P$ such that (i) $\P$ is freely generated as a
  left $\R$-module by the monomials $x^\mu$ with $\mu\in\Nn$, (ii) $x_i\star
  r-r\star x_i\in\R$ for all $r\in\R$ and (iii) $x_i\star x_j-x_j\star
  x_i\in\R+\R x_1+\cdots\R x_n$.  Obviously, any such extension is a
  polynomial algebra of solvable type in the sense of
  Definition~\ref{def:solvalg} for any degree compatible term order.  Other
  term orders generally do not respect the multiplication in $\P$.
  
  The classical example of such a {\small PBW} extension is the
  \emph{universal enveloping algebra} $\Uf(\gf)$ of a finite-dimensional Lie
  algebra $\gf$ which also explains the name: the Poincar\'e-Birkhoff-Witt
  theorem asserts that the monomials form a basis of these algebras
  \cite{var:lie}.  They still fit into the framework developed by Kandry-Rody
  and Weispfenning \cite{krw:ncgb}, as the $x_i$ do not act on the
  coefficients.  This is no longer the case for the more general \emph{skew
    enveloping algebras} $\Rf\#\Uf(\gf)$ where $\Rf$ is a $\kk$-algebra on
  which the elements of $\gf$ act as derivations
  \cite[Sect.~1.7.10]{mr:ncrings}.\bull
\end{example}

\begin{example}\label{ex:quantum}
  In all these examples, the coefficients $r_{\mu\nu}$ appearing in
  (\ref{eq:crmono}) are one; thus (\ref{eq:crij}) are classical commutation
  relations.  This is no longer true in the \emph{quantised enveloping
    algebras} $\Uf_h(\gf)$ introduced by Drinfeld \cite{vgd:qybe} and Jimbo
  \cite{mj:qybe}.  For these algebra it is non-trivial that a
  Poincar\'e-Birkhoff-Witt theorem holds; it was shown for general Lie
  algebras $\gf$ by Lusztig \cite{gl:quantum}.  Berger \cite{rb:qpbw}
  generalised this result later to a larger class of associative algebras, the
  so-called \emph{$q$-algebras}.  They are characterised by the fact that the
  polynomials $h_{ij}$ in (\ref{eq:crij}) are at most quadratic with the
  additional restriction that $h_{ij}$ may contain only those quadratic terms
  $x_kx_\ell$ that satisfy $i<k\leq\ell<j$ and $k-i=j-\ell$.  Thus any such
  algebra is a polynomial algebra of solvable type for any degree compatible
  term order.
  
  A simple concrete example is the \emph{$q$-Heisenberg algebra} for a real
  $q>0$ (and $q\neq1$).  Let $f$ be a function of a real variable $x$ lying in
  some appropriate function space.  Then we introduce the operators
  \begin{equation}\label{eq:qheisen}
    \delta_qf(x)=\frac{f(x)-f(qx)}{(1-q)x}\,,\quad
    \tau_qf(x)=f(qx)\,,\quad \hat xf(x)=xf(x)\,.
  \end{equation}
  It is straightforward to verify that these three operators satisfy the
  following $q$-deformed form of the classical Heisenberg commutation rules
  \begin{equation}\label{eq:qheiscr}
    \delta_q\star\hat x=\hat x\star\delta_q+\tau_q\,,\quad
    \delta_q\star\tau_q=q\tau_q\star\delta_q\,,\quad
    \tau_q\star\hat x=q\hat x\star\tau_q\,.
  \end{equation}
  Hence the algebra $\kk[\delta_q,\tau_q,\hat x]$ is a polynomial algebra of
  solvable type for any degree compatible term order (but also for any
  lexicographic order with $\tau_q\prec\delta_q$ and $\tau_q\prec\hat
  x$).\bull
\end{example}

\begin{example}\label{ex:grp}
  Let $(\P,\star,\prec)$ be a polynomial algebra of solvable type with a
  degree compatible term order $\prec$.  Then $\P$ is a filtered ring with
  respect to the standard filtration $\Sigma_q=\bigoplus_{i=0}^q\P_i$ and we
  may introduce the \emph{associated graded algebra} by setting
  $(\gr{\Sigma}{\P})_q=\Sigma_q/\Sigma_{q-1}$.  It is easy to see that
  $\gr{\Sigma}{\P}$ is again a polynomial algebra of solvable type for
  $\prec$.  If in (\ref{eq:basiccr}) $\deg{h_i(r)}=0$, $\deg{h_{ij}}\leq1$,
  $\rho_i=\id{\R}$ and $r_{ij}=1$ (which is for example the case for all
  Poincar\'e-Birkhoff-Witt extensions), then in fact
  $\gr{\Sigma}{\P}=(\P,\cdot)$, the commutative polynomial ring.  Such
  algebras are sometimes called \emph{almost commutative}
  \cite[Sect.~8.4.2]{mr:ncrings}.\bull
\end{example}

\begin{proposition}\label{prop:ore}
  If the ring\/ $\R$ is an integral domain, then any polynomial algebra\/
  $(\P,\star,\prec)$ of solvable type over it is an integral domain, too, and
  a left Ore domain.
\end{proposition}

\begin{proof}
  The first assertion is a trivial consequence of (\ref{eq:crsolvalg}): if
  $\R$ has no zero divisors, then $f\cdot g\neq0$ implies $f\star g\neq0$.
  Hence $\P$ does not contain any zero divisors.
  
  For the second one we must verify the \emph{left Ore conditions}
  \cite{pmc:alg2,oo:linear}: we must show that one can find for any two
  polynomials $f,g\in\P$ with $f\star g\neq0$ two further polynomials
  $\phi,\psi\in\P\setminus\{0\}$ such that $\phi\star f=\psi\star g$.  We
  describe now a concrete algorithm for this task.
  
  We set $\F_0=\{f,g\}$ and choose coefficients $r_0,s_0\in\R$ such that in
  the difference $r_0g\star f-s_0f\star g=\bar h_1$ the leading terms cancel.
  Then we perform a pseudoreduction of $\bar h_1$ with respect to $\F_0$.  It
  leads with an appropriately chosen coefficient $t_0\in\R$ to an equation of
  the form
  \begin{equation}\label{eq:ore1}
    t_0\bar h_1 = \phi_0\star f + \psi_0\star g + h_1
  \end{equation}
  where the remainder $h_1$ satisfies
  $\le{\prec}{h_1}\notin\lspan{\le{\prec}{\F_0}}$.  If $h_1=0$, we are done
  and the polynomials $\phi=t_0r_0g-\phi_0$ and $\psi=t_0s_0f+\psi_0$ form a
  solution of our problem.  By Part (iii) of Definition \ref{def:solvalg} we
  have $\le{\prec}{\bar h_1}\prec\le{\prec}{f}+\le{\prec}{g}$.  This implies
  by the monotonicity of term orders that
  $\le{\prec}{\phi_0}\prec\le{\prec}{g}$ and
  $\le{\prec}{\psi_0}\prec\le{\prec}{f}$.  Thus we have found a non-trivial
  solution.
  
  Otherwise we set $\F_1=\F_0\cup\{h_1\}$ and choose coefficients
  $r_1,s_1\in\R$ such that in the difference $r_1f\star h_1-s_1h_1\star f=\bar
  h_2$ the leading terms cancel.  Now we perform a pseudoreduction of $\bar
  h_2$ with respect to $\F_1$.  This computation yields a coefficient
  $t_1\in\R$ and polynomials $\phi_1,\psi_1,\rho_1\in\P$ such that
  \begin{equation}\label{eq:ore2}
    t_1\bar h_2 = \phi_1\star f + \psi_1\star g + \rho_1\star h_1 + h_2
  \end{equation} 
  where the remainder $h_2$ satisfies
  $\le{\prec}{h_2}\notin\lspan{\le{\prec}{\F_1}}$.  If $h_2=0$, then we are
  done, as we can substitute $h_1$ from (\ref{eq:ore1}) and obtain thus for
  our problem the solution
  $\phi=(t_1r_1f-\rho_1)\star(t_0r_0g-\phi_0)-t_1s_1h_1+\phi_1$ and
  $\psi=(t_1r_1f-\rho_1)\star(t_0s_0f-\psi_0)+\psi_1$.  By the same reasoning
  on the leading exponents as above, it is a non-trivial one.
  
  Otherwise we iterate: we set $\F_2=\F_1\cup\{h_2\}$, choose coefficients
  $r_2,s_2\in\R$ such that in the difference $r_2f\star h_2-s_2h_2\star f=\bar
  h_3$ the leading terms cancel, compute the remainder $h_3$ of a pseudo
  reduction of $\bar h_3$ with respect to $\F_2$ and so on.  If the iteration
  stops, i.\,e.\ if the remainder $h_N$ vanishes for some value $N\in\NN$,
  then we can construct non-zero polynomials $\phi$, $\psi$ with $\phi\star
  f=\psi\star g$ by substituting all remainders $h_i$ by their defining
  equations.  The iteration terminates by a simple Noetherian argument:
  $\lspan{\le{\prec}{\F_0}}\subset\lspan{\le{\prec}{\F_1}}\subset
  \lspan{\le{\prec}{\F_2}}\subset\cdots$ is a strictly ascending chain of
  monoid ideals in $\Nn$ and thus cannot be infinite.\qed
\end{proof}

Obviously, we can show by the same argument that $\P$ is a right Ore domain.
The Ore multipliers $\phi$, $\psi$ constructed in the proof above are not
unique.  Instead of always analysing differences of the form $r_{i}f\star
h_{i}-s_{i}h_{i}\star f$, we could have used differences of the form
$r_{i}g\star h_{i}-s_{i}h_{i}\star g$ or we could have alternated between
using $f$ and $g$ and so on.  In general, each ansatz will lead to different
multipliers.

We have given here a direct and in particular constructive proof that $\P$
satisfies the left and right Ore conditions.  Instead we could have tried to
invoke Theorem~2.1.15 of \cite{mr:ncrings} stating that any right Noetherian
integral domain is also a right Ore domain.  However, as we will see in the
next section, if the coefficient ring $\R$ of $\P$ is not a field, then the
question whether or not $\P$ is (left or right) Noetherian becomes nontrivial
in general.

\begin{example}\label{ex:so3}
  In the commutative polynomial ring one has always the trivial solution
  $\phi=g$ and $\psi=f$.  One might expect that in the non-commutative case
  one only has to add some lower terms to it.  However, this is not the case.
  Consider the universal enveloping algebra of the Lie algebra
  $\mathfrak{so}(3)$.  We may write it as
  $\Uf\bigl(\mathfrak{so}(3)\bigr)=\kk[x_1,x_2,x_3]$ with the multiplication
  $\star$ defined by the relations:
  \begin{equation}\label{eq:so3mult}
    \begin{aligned}
      x_1\star x_2&=x_1x_2\;,\qquad & x_2\star x_1&=x_1x_2-x_3\;,\\
      x_1\star x_3&=x_1x_3\;, & x_3\star x_1&=x_1x_3+x_2\;,\\
      x_2\star x_3&=x_2x_3\;, & x_3\star x_2&=x_2x_3-x_1\;.
    \end{aligned}
  \end{equation}
  This multiplication obviously respects any degree compatible term order but
  not the lexicographic order.  Choosing $f=x_1$ and $g=x_2$, possible
  solutions for $\phi\star f=\psi\star g$ are $\phi=x_2^2-1$ and
  $\psi=x_1x_2-2x_3$ or $\phi=x_1x_2+x_3$ and $\psi=x_1^2-1$.  They are easily
  constructed using the algorithm of the proof of Proposition \ref{prop:ore}
  once with $f$ and once with $g$.  Here we must use polynomials of degree
  $2$; it is not possible to find a solution of degree $1$.  \bull
\end{example}

\section{Hilbert's Basis Theorem for Solvable Algebras}\label{sec:hilbert}

A classical property of the ordinary polynomial ring $\P=\R[x_1,\dots,x_n]$,
which is crucial in the theory of Gr\"obner bases, is Hilbert's Basis Theorem.
For our more general class of polynomial algebras, it remains true only under
additional assumptions.  As $\P$ is generally non-commutative, we must
distinguish left, right and two-sided ideals and thus also study separately
whether $\P$ is left or right Noetherian.

With the exception of Section~\ref{sec:two}, we will exclusively work with
left ideals and thus do not introduce special notations.  This restriction to
left ideals is not only for convenience but stems from the fundamental
left-right asymmetry of Definition~\ref{def:solvalg} of a polynomial algebra
of solvable type where products $r\star x^\mu$ and $x^\mu\star r$ are treated
completely different.  For this reason we discuss only the question when $\P$
is left Noetherian (see also Remark~\ref{rem:rightnoether} below).

Most classical proofs of Hilbert's Basis Theorem consider only the univariate
case and then extend inductively to an arbitrary (but finite) number of
variables.  However, this inductive approach is not possible in arbitrary
polynomial algebras of solvable type, as the multiplication $\star$ does not
necessarily restrict to a subalgebra with fewer variables.  A simple
counterexample is provided by the universal enveloping algebra
$\Uf\bigl(\mathfrak{so}(3)\bigr)$ introduced in Example~\ref{ex:so3} where
$\star$ cannot be restricted to the subspace $\kk[x_1,x_2]$ since $x_2\star
x_1=x_1x_2-x_3$.  This observation motivates the following definition.

\begin{definition}\label{def:iteralg}
  $(\P,\star,\prec)$ is an\/ \emph{iterated polynomial algebra of solvable
    type}, if (i)\/ $\P=\R[x_1][x_2]\cdots[x_n]$ where each intermediate
  ring\/ $\P_k=\R[x_1][x_2]\cdots[x_k]$ is again solvable for the
  corresponding restrictions of the multiplication\/ $\star$ and the term
  order\/ $\prec$ and (ii) we have the equality\/
  $x_k\star\P_{k-1}+\P_{k-1}=\P_{k-1}\star x_k+\P_{k-1}$.
\end{definition}

For iterated polynomial algebras of solvable type we may apply the usual
inductive technique for proving a basis theorem, but we must still impose
further conditions on the multiplication $\star$.  The following result is
proven in \cite[Theorem~1.2.9]{mr:ncrings} for Ore algebras, but it is a
trivial exercise to verify that the proof remains valid for more general
algebras.

\begin{theorem}\label{thm:iternoether}
  If\/ $(\P,\star,\prec)$ is an iterated polynomial algebra of solvable type
  over a left Noetherian ring\/ $\R$, then $\P$ is a left Noetherian ring,
  too.
\end{theorem}

The second condition in Definition \ref{def:iteralg} cannot be omitted, if a
basis theorem is to hold.  McConnell and Robson
\cite[Example~1.2.11]{mr:ncrings} provide a concrete counterexample of a
univariate polynomial ring of solvable type which violates this condition and
which is neither left nor right Noetherian.

With some complications, the central (univariate) arguments in the proof of
Theorem \ref{thm:iternoether} can be directly generalised to multivariate
polynomial rings.  However, this requires again certain assumptions on the
commutation relations (\ref{eq:crmono}) in order to ensure that all needed
computations are possible.

\begin{definition}\label{def:ccr}
  The polynomial algebra of solvable type\/ $(\P,\star,\prec)$ has\/
  \emph{centred commutation relations}, if (i) there exists a field\/
  $\kk\subseteq\R$ lying in the centre of\/ $\R$, (ii) the functions\/
  $\rho_\mu$ in (\ref{eq:crmono1}) are of the form\/
  $\rho_\mu(r)=\bar\rho_\mu(r)r$ with functions\/
  $\bar\rho_\mu:\R\rightarrow\kk$ and (iii) we have\/ $r_{\mu\nu}\in\kk$ in
  (\ref{eq:crmono2}).
\end{definition}

Using K\"onig's Lemma, Kredel proved in his thesis \cite{hk:solvpoly} the
following version of Hilbert's Basis Theorem.

\begin{theorem}\label{thm:subfieldnoether}
  Let\/ $(\P,\star,\prec)$ be a polynomial algebra of solvable type with
  centred commutation relations over a left Noetherian coefficient ring\/
  $\R$.  Then\/ $\P$ is left Noetherian, too.
\end{theorem}

A third proof assumes that the ring $\P$ possesses a \emph{filtration}
$\Sigma$.  Using an approach detailed in \cite{jeb:rdo} for the special case
of the Weyl algebra (but which does not use any special properties of the Weyl
algebra), one obtains the following general result where it is not even
necessary to assume that $\P$ is a polynomial ring.  As discussed in
Example~\ref{ex:grp}, this result covers many of the polynomial algebras of
solvable type that have appeared in the literature so far.

\begin{theorem}
  Let\/ $\Sigma$ be a filtration on the ring\/ $\P$.  If the associated graded
  ring\/ $\gr{\Sigma}{\P}$ is left Noetherian, then\/ $\P$ is left Noetherian,
  too.
\end{theorem}

Because of Condition (iii) in Definition \ref{def:solvalg} we can define
Gr\"obner bases for ideals in algebras of solvable type.  In the case that
$\R$ is a (commutative) field $\kk$, this is straightforward and from now on
we will restrict to this case; the general case will be discussed only in
Section~\ref{sec:ring}.

\begin{definition}\label{def:gbsolv}
  Let\/ $(\P,\star,\prec)$ be a polynomial algebra of solvable type over a
  field\/~$\kk$ and\/ $\I\subseteq\P$ a left ideal.  A finite set\/
  $\G\subset\P$ is a\/ \emph{Gr\"obner basis} of\/ $\I$ (for the term order\/
  $\prec$), if\/ $\lspan{\le{\prec}{\G}}=\le{\prec}{\I}$.
\end{definition}

For the ordinary multiplication this definition reduces to the classical one.
The decisive point, explaining the conditions imposed in
Definition~\ref{def:solvalg}, is that normal forms with respect to a finite
set $\F\subset\P$ may be computed in algebras of solvable type in precisely
the same way as in the ordinary polynomial ring.  Assume we are given a
polynomial $f\in\P$ such that $\le{\prec}{g}\mid\le{\prec}{f}$ for some
$g\in\G$ and set $\mu=\le{\prec}{f}-\le{\prec}{g}$.  If we consider
$g_\mu=x^\mu\star g$, then by (iii) $\le{\prec}{g_\mu}=\le{\prec}{f}$.
Setting $d=\lc{\prec}{f}/\lc{\prec}{g_\mu}$, we find by (ii) that
$\le{\prec}{(f-dg_\mu)}\prec\le{\prec}{f}$.  Hence we may use the usual
algorithms for computing normal form; in particular, they always terminate by
the same argument as in the ordinary case.  Note that in general
$d\neq\lc{\prec}{f}/\lc{\prec}{g}$, if $r\neq1$ in (\ref{eq:crsolvalg}), and
that normal form computations are typically more expensive due to the
appearance of the additional polynomial $h$ in (\ref{eq:crsolvalg}).

The classical Gr\"obner basis theory can be extended straightforwardly to
polynomial algebras of solvable type
\cite{apel:diss,bglc:quantum,bgv:algnoncomm,krw:ncgb,hk:solvpoly,vl:ade,vl:diss},
as most proofs are based on the computation of normal forms.  The remaining
arguments mostly take place in the monoid $\Nn$ and thus can be applied
without changes.  In particular, a trivial adaption of the standard approach
leads to the following result crucial for the termination of Buchberger's
algorithm.

\begin{theorem}\label{thm:noetherfield}
  Let\/ $(\P,\star,\prec)$ be a polynomial algebra of solvable type over a
  field.  Then $\P$ is a left Noetherian ring and every left ideal\/
  $\I\subseteq\P$ possesses a Gr\"obner basis with respect to\/ $\prec$.
\end{theorem}

\begin{remark}\label{rem:rightnoether}
  Even in the case of a coefficient field we cannot generally expect $\P$ to
  be a \emph{right} Noetherian ring, too; a concrete counterexample is
  provided again by McConnell and Robson \cite[Example~1.2.11]{mr:ncrings}.
  In the proof of Theorem~\ref{thm:noetherfield} one essentially uses that in
  normal form computations one always multiplies with elements of $\P$ from
  the left.  Because of the already above mentioned left-right asymmetry of
  Definition~\ref{def:solvalg}, right ideals show in general a completely
  different behaviour.  In order to obtain right Noetherian rings we must
  either adapt correspondingly our definition of a solvable algebra or impose
  additional conditions on the commutation relations (\ref{eq:crmono}).
  
  The simplest possibility is to require that all the maps $\rho_{\mu}$ in
  (\ref{eq:crmono}) are automorphisms (by Proposition~\ref{prop:fixmult} it
  suffices, if the maps $\rho_{i}$ in (\ref{eq:crop}) satisfy this condition).
  In this case we have $\kk\star x_{i}+\kk=x_{i}\star\kk+\kk$ for all
  variables $x_{i}$ implying that we can rewrite any polynomial $f=\sum_{\mu}
  c_{\mu}x^{\mu}$ in the ``reverse'' form $f=\sum_{\mu} x^{\mu}\star \tilde
  c_{\mu}$.  Now a straightforward adaption of the classical proof of Theorem
  \ref{thm:noetherfield} shows that the ring $\P$ is also right Noetherian.
  \bull
\end{remark}

We do not give more details on Gr\"obner bases, as they can be found in the
above cited references.  Instead we will present in the next section a
completely different approach leading to involutive bases.

\section{Involutive Bases}\label{sec:invbas}

We proceed to define involutive bases for left ideals in polynomial algebras
of solvable type.  In principle, we could at once consider submodules of free
modules over such an algebra.  As this only complicates the notation, we
restrict to the ideal case and the extension to submodules goes as for
Gr\"obner bases.

\begin{definition}\label{def:invbas}
  Let\/ $(\P,\star,\prec)$ be a polynomial algebra of solvable type over a
  field\/~$\kk$ and\/ $\I\subseteq\P$ a non-zero left ideal.  A finite
  subset\/ $\H\subset\I$ is a\/ \emph{weak involutive basis} of\/ $\I$ for an
  involutive division\/ $L$ on\/ $\Nn$, if its leading exponents\/
  $\le{\prec}{\H}$ form a weak involutive basis of the monoid ideal\/
  $\le{\prec}{\I}$.  The subset\/ $\H$ is a\/ \emph{(strong) involutive basis}
  of\/ $\I$, if\/ $\le{\prec}{\H}$ is a strong involutive basis of\/
  $\le{\prec}{\I}$ and no two distinct elements of\/ $\H$ have the same
  leading exponents.
\end{definition}

\begin{remark}
  This definition of an involutive basis is different from the original one
  given by Gerdt and Blinkov \cite{gb:invbas}.  Firstly, the notion of a weak
  basis is new.  Secondly, we do not require that an involutive basis is
  completely autoreduced.  Finally, our approach is a natural extension of
  Definition~\ref{def:gbsolv} of a Gr\"obner basis in the ring $\P$ whereas
  the approach of Gerdt and Blinkov \cite{gb:invbas} is closer to the
  constructive characterisation of Gr\"obner bases via $S$-polynomials.
  However, we will see below that essentially both approaches are
  equivalent.\bull
\end{remark}

Definition \ref{def:invbas} implies immediately that any weak involutive basis
is a Gr\"obner basis.  As in Section~\ref{sec:invdiv}, we call any finite set
$\F\subset\P$ \emph{(weakly) involutive}, if it is a (weak) involutive basis
of the ideal~$\lspan{\F}$ generated by it.

\begin{definition}\label{def:invspanpoly}
  Let\/ $\F\subset\P\setminus\{0\}$ be a finite set and\/ $L$ an involutive
  division on\/ $\Nn$.  We assign to each element\/ $f\in\F$ a set of\/
  \emph{multiplicative variables}
  \begin{equation}\label{eq:multvar}
    \mult{X}{L,\F,\prec}{f}=\bigl\{x_i\mid
    i\in\mult{N}{L,\le{\prec}{\F}}{\le{\prec}{f}}\bigr\}\;.   
  \end{equation} 
  The\/ \emph{involutive span} of\/ $\F$ is then the set
  \begin{equation}\label{eq:ispanpoly}
    \ispan{\F}{L,\prec}=
        \sum_{f\in\F}\kk\bigl[\mult{X}{L,\F,\prec}{f}\bigr]\star f
        \subseteq\lspan{\F}\;.
  \end{equation}
\end{definition}

An important aspect of Gr\"obner bases is the existence of standard
representations for ideal elements.  For (weak) involutive bases a similar
characterisation exists and in the case of strong bases we even obtain unique
representations.

\begin{theorem}\label{thm:invnormalrep}
  Let\/ $\I\subseteq\P$ be a non-zero ideal,\/ $\H\subset\I\setminus\{0\}$ a
  finite set and\/ $L$ an involutive division on\/ $\Nn$.  Then the following
  two statements are equivalent.
  \begin{description}
  \item[\phantom{i}{\upshape (i)}] The set\/ $\H$ is a weak involutive basis
    of\/ $\I$ with respect to\/ $L$ and\/ $\prec$.
  \item[{\upshape (ii)}] Every polynomial\/ $f\in\I$ can be written in the
    form
    \begin{equation}\label{eq:invnormalrep}
      f=\sum_{h\in\H}P_h\star h
    \end{equation}
    where the coefficients\/ $P_h\in\kk[\mult{X}{L,\H,\prec}{h}]$ satisfy\/
    $\le{\prec}{(P_h\star h)}\preceq\le{\prec}{f}$ for all polynomials\/
    $h\in\H$.
  \end{description}
  $\H$ is a strong involutive basis, if and only if the representation
  (\ref{eq:invnormalrep}) is unique.
\end{theorem}

\begin{proof}
  Let us first assume that the set $\H$ is a weak involutive basis.  Take an
  arbitrary polynomial $f\in\I$.  According to Definition~\ref{def:invbas},
  its leading exponent $\le{\prec}{f}$ lies in the involutive cone
  $\cone{L,\le{\prec}{\H}}{h}$ of at least one element $h\in\H$.  Let
  $\mu=\le{\prec}{f}-\le{\prec}{h}$ and set $f_1=f-cx^\mu\star h$ where the
  coefficient $c\in\kk$ is chosen such that the leading terms cancel.
  Obviously, $f_1\in\I$ and $\le{\prec}{f_1}\prec\le{\prec}{f}$.  Iteration
  yields a sequence of polynomials $f_i\in\I$.  After a finite number of steps
  we must reach $f_N=0$, as the leading exponents are always decreasing and by
  assumption the leading exponent of \emph{any} polynomial in $\I$ possesses
  an involutive divisor in $\le{\prec}{\H}$.  But this implies the existence
  of a representation of the form (\ref{eq:invnormalrep}).
  
  Now assume that $\H$ is even a strong involutive basis and take an
  involutive standard representation (\ref{eq:invnormalrep}).  By definition
  of a strong basis, there exists one and only one generator $h\in\H$ such
  that $\le{\prec}{(P_h\star h)}=\le{\prec}{f}$.  This fact determines
  uniquely $\le{\prec}{P_h}$.  Applying the same argument to
  $f-(\lt{\prec}{P_h})\star h$ shows by recursion that the representation
  (\ref{eq:invnormalrep}) is indeed unique.
  
  For the converse note that (ii) trivially implies that
  $\le{\prec}{f}\in\ispan{\le{\prec}{\H}}{L,\prec}$ for any polynomial
  $f\in\I$.  Thus $\le{\prec}{\I}\subseteq\ispan{\le{\prec}{\H}}{L,\prec}$.
  As the converse inclusion is obvious, we have in fact an equality and $\H$
  is a weak involutive basis.
  
  Now let us assume that the set $\H$ is only a weak but not a strong
  involutive basis of $\I$.  This implies the existence of two generators
  $h_1,h_2\in\H$ such that $\cone{L,\le{\prec}{\H}}{\le{\prec}{h_2}}\subset
  \cone{L,\le{\prec}{\H}}{\le{\prec}{h_1}}$.  Hence we have
  $\lm{\prec}{h_2}=\lm{\prec}{(cx^\mu\star h_1)}$ for suitably chosen
  $c\in\kk$ and $\mu\in\Nn$.  Consider the polynomial $h_2-cx^\mu\star
  h_1\in\I$.  If it vanishes, we have found a non-trivial involutive standard
  representation of~$0$.  Otherwise an involutive standard representation
  $h_2-cx^\mu\star h_1=\sum_{h\in\H}P_h\star h$ with
  $P_h\in\kk[\mult{X}{L,\H,\prec}{h}]$ exists.  Setting $P'_h=P_h$ for all
  generators $h\neq h_1,h_2$ and $P'_{h_1}=P_{h_1}+cx^\mu$,
  $P'_{h_2}=P_{h_2}-1$ yields again a non-trivial involutive standard
  representation $0=\sum_{h\in\H}P'_h\star h$.  The existence of such a
  non-trivial representation of $0$ immediately implies that
  (\ref{eq:invnormalrep}) cannot be unique.  Thus only for a strong involutive
  basis the involutive standard representation is always unique.\qed
\end{proof}

\begin{corollary}\label{cor:charinvbas}
  Let the set\/ $\H$ be a weak involutive basis of the left ideal\/
  $\I\subseteq\P$.  Then\/ $\ispan{\H}{L,\prec}=\I$.
\end{corollary}

\begin{example}\label{ex:weakbasis}
  It is \emph{not} true that any set $\F$ with $\ispan{\F}{L,\prec}=\I$ is a
  weak involutive basis of the ideal $\I$.  Consider in the ordinary
  polynomial ring $\kk[x,y]$ the ideal $\I$ generated by the two polynomials
  $f_1=y^2$ and $f_2=y^2+x^2$.  If we order the variables as $x_1=x$ and
  $x_2=y$, then the set $\F=\{f_1,f_2\}$ trivially satisfies
  $\ispan{\F}{J,\prec}=\I$, as with respect to the Janet division all
  variables are multiplicative for each generator.  However,
  $\le{\prec}{\F}=\{[0,2]\}$ does \emph{not} generate $\le{\prec}{\I}$, as
  obviously $[2,0]\in\le{\prec}{\I}\setminus\lspan{\{[0,2]\}}$.  Thus $\F$ is
  not a weak Janet basis (neither is the autoreduced set $\F'=\{y^2,x^2\}$, as
  $x^2y\notin\ispan{\F'}{J,\prec}$).\bull
\end{example}

\begin{proposition}\label{prop:weakbasis}
  Let\/ $\I\subseteq\P$ be an ideal and $\H\subset\P$ a weak involutive basis
  of it for the involutive division\/ $L$.  Then there exists a subset\/
  $\H'\subseteq\H$ which is a strong involutive basis of\/ $\I$.
\end{proposition}

\begin{proof}
  If the set $\le{\prec}{\H}$ is already a strong involutive basis of
  $\le{\prec}{\I}$, we are done.  Otherwise $\H$ contains polynomials $h_1$,
  $h_2$ such that $\le{\prec}{h_1}\idiv{L,\le{\prec}{\H}}\le{\prec}{h_2}$.
  Consider the subset $\H'=\H\setminus\{h_2\}$.  As in the proof of
  Proposition~\ref{prop:strongbas} one easily shows that
  $\le{\prec}{\H'}=\le{\prec}{\H}\setminus\{\le{\prec}{h_2}\}$ is still a weak
  involutive basis of $\le{\prec}{\I}$ and thus $\H'$ is still a weak
  involutive basis of $\I$.  After a finite number of such eliminations we
  must reach a strong involutive basis.\qed
\end{proof}

Given this result, one may wonder why we have introduced the notion of a weak
basis.  The reason is that in more general situations like computations in
local rings or polynomial algebras over coefficient rings (treated in later
sections) strong bases rarely exist.

\begin{definition}\label{def:reduce}
  Let\/ $\F\subset\P$ be a finite set and\/ $L$ an involutive division.  A
  polynomial\/ $g\in\P$ is\/ \emph{involutively reducible} with respect to\/
  $\F$, if it contains a term\/ $x^\mu$ such that\/
  $\le{\prec}{f}\idiv{L,\le{\prec}{\F}}\mu$ for some\/ $f\in\F$.  It is in\/
  \emph{involutive normal form} with respect to $\F$, if it is not
  involutively reducible.  The set\/ $\F$ is\/ \emph{involutively
    autoreduced}, if no polynomial\/ $f\in\F$ contains a term\/ $x^\mu$ such
  that another polynomial\/ $f'\in\F\setminus\{f\}$ exists with\/
  $\le{\prec}{f'}\idiv{L,\le{\prec}{\F}}\mu$.
\end{definition}

\begin{remark}\label{rem:auto}
  The definition of an involutively autoreduced set \emph{cannot} be
  formulated more concisely by saying that each $f\in\F$ is in involutive
  normal form with respect to $\F\setminus\{f\}$.  If we are not dealing with
  a global division, the removal of $f$ from $\F$ will generally change the
  assignment of the multiplicative indices and thus affect the involutive
  divisibility.\bull
\end{remark}

An \emph{obstruction to involution} is a polynomial
$g\in\lspan{\F}\setminus\ispan{\F}{L,\prec}$ possessing a (necessarily
non-involutive) standard representation with respect to $\F$.  We will later
see that these elements make the difference between an involutive and an
arbitrary Gr\"obner basis.

\begin{example}\label{ex:obstr}
  Consider the set $\F=\{f_1,f_2,f_3\}\subset\kk[x,y,z]$ with the polynomials
  $f_1=z^2-xy$, $f_2=yz-x$ and $f_3=y^2-z$.  For any degree compatible term
  order, the leading terms of $f_2$ and $f_3$ are unique.  For $f_1$ we have
  two possibilities: if we use the degree lexicographic order (i.\,e.\ for
  $x\prec y\prec z$), it is $z^2$, for the degree inverse lexicographic order
  (i.\,e.\ for $x\succ y\succ z$) the leading term is $xy$.
  
  In the first case, $\ispan{\F}{J,\prec_{\mbox{\scriptsize
        deglex}}}=\lspan{\F}$, so that for this term order $\F$ is a Janet
  basis, i.\,e.\ an involutive basis with respect to the Janet division,
  although we have not yet the necessary tools to prove this fact.  In the
  second case, $f_4=z^3-x^2=zf_1+xf_2\in\lspan{\F}$ does not possess a
  standard representation and $\F$ is not even a Gr\"obner basis.  Adding
  $f_4$ to $\F$ yields a Gr\"obner basis $\G$ of $\lspan{\F}$, as one may
  easily check.  But this makes $z$ non-multiplicative for $f_2$ and
  $f_5=zf_2$ is now an obstruction to involution of $\G$, as it is not
  involutively reducible with respect to the Janet division.  In fact, the set
  $\F'=\{f_1,f_2,f_3,f_4,f_5\}$ is the smallest Janet basis of $\I$ for this
  term order, as it is not possible to remove an element.  Note that this
  second basis is not only larger but also contains polynomials of higher
  degree.\bull
\end{example}

\begin{remark}\label{rem:invgb}
  If $\G$ is a Gr\"obner basis of the ideal $\I$, then any element of $\I$ has
  a standard representation.  But this does not imply that for a given
  division $L$ the ideal $\I$ is free of obstructions to involution. In order
  to obtain at least a weak involutive basis, we must add further elements of
  $\I$ to $\G$ until $\ispan{\le{\prec}{\G}}{L}=\le{\prec}{\I}$.  Obviously,
  this observation allows us to reduce the construction of a polynomial
  involutive basis to a Gr\"obner basis computation plus a monomial
  completion.  But we will see later that better possibilities exist.
  
  It follows that in general involutive bases are not reduced Gr\"obner bases,
  as we already observed in Example~\ref{ex:obstr}.  For
  $\prec_{\mbox{\scriptsize deglex}}$ the set $\F$ was simultaneously a Janet
  basis and a reduced Gr\"obner basis.  But for $\prec_{\mbox{\scriptsize
      deginvlex}}$ the reduced Gr\"obner basis is $\F\cup\{f_4\}$, whereas a
  Janet basis requires in addition the polynomial $f_5$.  We will see in
  Part~II that this ``redundancy'' in involutive bases is the key for their
  use in the structure analysis of polynomial ideals and modules.\bull
\end{remark}

It often suffices, if one does not consider all terms in $g$ but only the
leading term $\lt{\prec}{g}$: the polynomial $g$ is \emph{involutively head
  reducible}, if $\le{\prec}{f}\idiv{L,\le{\prec}{\F}}\le{\prec}{g}$ for
some\/ $f\in\F$.  Similarly, the set $\F$ is \emph{involutively head
  autoreduced}, if no leading exponent of an element $f\in\F$ is involutively
divisible by the leading exponent of another element $f'\in\F\setminus\{f\}$.
Note that the definition of a strong involutive basis immediately implies that
it is involutively head autoreduced.

As involutive reducibility is a restriction of ordinary reducibility,
involutive normal forms can be determined with trivial adaptions of the
familiar algorithms.  The termination follows by the same argument as usual,
namely that any term order is a well-order.  If $g'$ is an involutive normal
form of $g\in\P$ with respect to the set $\F$ for the division $L$, then we
write $g'=\NF{\F,L,\prec}{g}$, although involutive normal forms are in general
not unique (like ordinary normal forms).  Depending on the order in which
reductions are applied different results are obtained.

The ordinary normal form is unique, if and only if it is computed with respect
to a Gr\"obner basis; this property is often used as an alternative definition
of Gr\"obner bases.  The situation is somewhat different for the involutive
normal form.

\begin{lemma}\label{lem:head} 
  The sum in\/ (\ref{eq:ispanpoly}) is direct, if and only if the finite set\/
  $\F\subset\P\setminus\{0\}$ is involutively head autoreduced with respect to
  the involutive division\/ $L$.
\end{lemma}

\begin{proof}
  One direction is obvious.  For the converse, let $f_1$, $f_2$ be two
  distinct elements of $\F$ and $X_i=\mult{X}{L,\F,\prec}{f_i}$ their
  respective sets of multiplicative variables for the division $L$.  Assume
  that two polynomials $P_i\in\kk[X_i]$ exist with $P_1\star f_1=P_2\star f_2$
  and hence $\le{\prec}{(P_1\star f_1)}=\le{\prec}{(P_2\star f_2)}$.  As the
  multiplication $\star$ respects the term order $\prec$, this implies that
  $\cone{L,\le{\prec}{\F}}{\le{\prec}{f_1}}\cap
  \cone{L,\le{\prec}{\F}}{\le{\prec}{f_2}}\neq\emptyset$.  Thus one of the
  involutive cones is completely contained in the other one and either
  $\le{\prec}{f_1}\idiv{L,\le{\prec}{\F}}\le{\prec}{f_2}$ or
  $\le{\prec}{f_2}\idiv{L,\le{\prec}{\F}}\le{\prec}{f_1}$ contradicting that
  $\F$ is involutively head autoreduced.\qed
\end{proof}

\begin{proposition}\label{prop:invnormalform}
  If the finite set\/ $\F\subset\P\setminus\{0\}$ is involutively head
  autoreduced, every polynomial\/ $g\in\P$ has a unique involutive normal
  form\/ $\NF{\F,L,\prec}{g}$.
\end{proposition}

\begin{proof}
  If $0$ is an involutive normal form of $g$, then obviously
  $g\in\ispan{\F}{L,\prec}$.  Conversely, assume that
  $g\in\ispan{\F}{L,\prec}$, i.\,e.\ the polynomial $g$ can be written in the
  form $g=\sum_{f\in\F}P_f\star f$ with $P_f\in\kk[\mult{X}{L,\F,\prec}{f}]$.
  As $\F$ is involutively head autoreduced, the leading terms of the summands
  never cancel (see the proof of Lemma~\ref{lem:head}).  Thus
  $\le{\prec}{g}=\le{\prec}{(P_f\star f)}$ for some $f\in\F$ and any
  polynomial $g\in\ispan{\F}{L,\prec}$ is involutively head reducible with
  respect to $\F$.  Each reduction step in an involutive normal form algorithm
  leads to a new polynomial $g'\in\ispan{\F}{L,\prec}$ with
  $\le{\prec}{g'}\preceq\le{\prec}{g}$.  If the leading term is reduced, we
  even get $\le{\prec}{g'}\prec\le{\prec}{g}$.  As each terminating normal
  form algorithm must sooner or later reduce the leading term, we eventually
  obtain $0$ as unique involutive normal form of any
  $g\in\ispan{\F}{L,\prec}$.
  
  Let $g_1$ and $g_2$ be two involutive normal forms of the polynomial $g$.
  Obviously, $g_1-g_2\in\ispan{\F}{L,\prec}$.  By definition of a normal form,
  neither $g_1$ nor $g_2$ contain any term involutively reducible with respect
  to $\F$ and the same holds for $g_1-g_2$.  Hence the difference $g_1-g_2$ is
  also in involutive normal form and by our considerations above we must have
  $g_1-g_2=0$.\qed
\end{proof}

\begin{proposition}\label{prop:normalform}
  The ordinary and the involutive normal form of any polynomial\/ $g\in\P$
  with respect to a finite weakly involutive set\/ $\F\subset\P\setminus\{0\}$
  are identical.
\end{proposition}

\begin{proof}
  Recalling the proof of the previous proposition, we see that we used the
  assumption that $\F$ was involutively head autoreduced only for proving the
  existence of a generator $f\in\F$ such that
  $\le{\prec}{f}\idiv{L,\le{\prec}{\F}}\le{\prec}{g}$ for every polynomial
  $g\in\ispan{\F}{L,\prec}$.  But obviously this property is also implied by
  the definition of a weak involutive basis.  Thus by the same argument as
  above, we conclude that the involutive normal form with respect to a weakly
  involutive set is unique.  For Gr\"obner bases the uniqueness of the
  ordinary normal form is a classical property and any weak involutive basis
  is also a Gr\"obner basis.  As a polynomial in ordinary normal form with
  respect to $\F$ is trivially in involutive normal form with respect to $\F$,
  too, the two normal forms must coincide.\qed
\end{proof}

Finally, we extend the notion of a minimal involutive basis from $\Nn$ to
$\P$.  This is done in the same manner as in the theory of Gr\"obner bases.

\begin{definition}\label{def:polyminbas}
  Let\/ $\I\subseteq\P$ be a non-zero ideal and\/ $L$ an involutive division.  An
  involutive basis\/ $\H$ of\/ $\I$ with respect to\/ $L$ is\/ \emph{minimal},
  if\/ $\le{\prec}{\H}$ is the minimal involutive basis of the monoid ideal\/
  $\le{\prec}{\I}$ for the division\/ $L$.
\end{definition}

By Proposition \ref{prop:monomin}, we find that for a globally defined
division like the Pommaret division any involutive basis is minimal.
Uniqueness requires two additional assumptions.  First of all, our definition
of an involutive basis requires only that it is involutively head autoreduced;
for uniqueness we obviously need a full involutive autoreduction.  Secondly,
we must normalise the leading coefficients to one, i.\,e.\ we must take a
\emph{monic} basis.

\begin{proposition}\label{prop:uniminbas}
  Let\/ $\I\subseteq\P$ be a non-zero ideal and\/ $L$ an involutive division.
  Then\/ $\I$ possesses at most one monic, involutively autoreduced, minimal
  involutive basis for the division\/ $L$.
\end{proposition}

\begin{proof}
  Assume that $\H_1$ and $\H_2$ are two different monic, involutively
  autoreduced, minimal involutive bases of $\I$ with respect to $L$ and
  $\prec$.  By definition of a minimal involutive bases, this implies that
  $\le{\prec}{\H_1}=\le{\prec}{\H_2}$.  As $\H_1$ and $\H_2$ are not
  identical, we must have two polynomials $h_1\in\H_1$ and $h_2\in\H_2$ such
  that $\le{\prec}{h_1}=\le{\prec}{h_2}$ but $h_1\neq h_2$.  Now consider the
  polynomial $h=h_1-h_2\in\I$.  Its leading exponent must lie in the
  involutive span of $\le{\prec}{\H_1}=\le{\prec}{\H_2}$.  On the other hand,
  the term $\lt{\prec}{h}$ must be contained in either $h_1$ or $h_2$.  But
  this implies that either $\H_1$ or $\H_2$ is not involutively
  autoreduced.\qed
\end{proof}

\section{Monomial Completion}\label{sec:monocompl}

We turn to the question of the actual construction of involutive bases.
Unfortunately, for arbitrary involutive division no satisfying solution is
known so far.  In the monomial case, one may follow a brute force approach,
namely performing a breadth first search through the tree of all possible
completions.  Obviously, it terminates only, if a finite basis exists.  But
for divisions satisfying some additional properties one can design a fairly
efficient completion algorithm.

The first problem in constructing an involutive completion of a finite subset
$\N\subset\Nn$ for a division $L$ is to check \emph{effectively} whether $\N$
is already involutive.  The trouble is that we do not know a priori where
obstructions to involution might lie.  If we denote by $1_j\in\Nn$ the multi
index where all entries are zero except the $j$th one which is one, then the
multi indices $\nu+1_j$ with $\nu\in\N$ and $j\in\nmult{N}{L,\N}{\nu}$ are a
natural first guess.

\begin{definition}\label{def:localinv}
  The finite set\/ $\N\subset\Nn$ is\/ \emph{locally involutive} for the
  involutive division\/ $L$, if\/ $\nu+1_j\in\ispan{\N}{L}$ for every
  non-multiplicative index\/ $j\in\nmult{N}{L,\N}{\nu}$ of every multi index\/
  $\nu\in\N$.
\end{definition}

Obviously, local involution is easy to check effectively.  However, while
(weak) involution obviously implies local involution, the converse does not
necessarily hold.  A concrete counter example was given by Gerdt and Blinkov
\cite{gb:invbas}.  But they also discovered that for many divisions the
converse is in fact true and thus for such divisions we can effectively decide
involution.

\begin{definition}\label{def:contdiv}
  Let\/ $L$ be an involutive division and\/ $\N\subset\Nn$ a finite set.  Let
  furthermore\/ $(\nu^{(1)},\dots,\nu^{(t)})$ be a finite sequence of elements
  of\/ $\N$ where every multi index\/ $\nu^{(k)}$ with\/ $k<t$ has a
  non-multiplicative index\/ $j_k\in\nmult{N}{L,\N}{\nu^{(k)}}$ such that\/
  $\nu^{(k+1)}\idiv{L,\N}\nu^{(k)}+1_{j_k}$.  The division\/ $L$ is\/
  \emph{continuous}, if any such sequence consists only of distinct elements,
  i.\,e.\ if\/ $\nu^{(k)}\neq\nu^{(\ell)}$ for all\/ $k\neq\ell$.
\end{definition}

\begin{proposition}\label{prop:localinv}
  For a continuous division\/ $L$, any locally involutive set\/ $\N\subset\Nn$
  is weakly involutive.
\end{proposition}

\begin{proof}
  Let the set $\Sigma$ contain those obstructions to involution that are of
  minimal length.\footnote{The length $|\nu|$ of a multi index $\nu\in\Nn$ is
    the sum of its entries, i.\,e.\ the degree of the monomial $x^\nu$.}  We
  claim that for a continuous division $L$ all multi indices $\sigma\in\Sigma$
  are of the form $\nu+1_j$ with $\nu\in\N$ and $j\in\nmult{N}{L,\N}{\nu}$.
  This immediately implies our proposition: since for a locally involutive set
  all such multi indices are contained in $\ispan{\N}{L}$, we must have
  $\Sigma=\emptyset$ and thus $\lspan{\N}=\ispan{\N}{L}$.
  
  In order to prove our claim, we choose a $\sigma\in\Sigma$ for which no
  $\nu\in\N$ exists with $\sigma=\nu+1_j$.  We collect in $\N_\sigma$ all
  divisors $\nu\in\N$ of $\sigma$ of maximal length.  Let $\nu^{(1)}$ be an
  element of $\N_\sigma$; by assumption the multi index
  $\mu^{(1)}=\sigma-\nu^{(1)}$ satisfies $|\mu^{(1)}|>1$ and at least one
  non-multiplicative index $j_1\in\nmult{N}{L,\N}{\nu^{(1)}}$ exists with
  $\mu^{(1)}_{j_1}>0$.  By the definition of $\Sigma$ we have
  $\nu^{(1)}+1_{j_1}\in\ispan{\N}{L}$.  Thus a multi index $\nu^{(2)}\in\N$
  exists with $\nu^{(2)}\idiv{L,\N}\nu^{(1)}+1_{j_1}$.  This implies
  $\nu^{(2)}\mid\sigma$ and we set $\mu^{(2)}=\sigma-\nu^{(2)}$.  By the
  definition of the set $\N_\sigma$ we have $|\nu^{(2)}|\leq|\nu^{(1)}|$.
  Hence $\nu^{(2)}+1_j\in\ispan{\N}{L}$ for all $j$.
  
  Choose a non-multiplicative index $j_2\in\nmult{N}{L,\N}{\nu^{(2)}}$ with
  $\mu^{(2)}_{j_2}>0$.  Such an index exists as otherwise
  $\sigma\in\ispan{\N}{L}$.  By the same arguments as above, a multi index
  $\nu^{(3)}\in\N$ exists with $\nu^{(3)}\idiv{L,\N}\nu^{(2)}+1_{j_2}$ and
  $|\nu^{(3)}|\leq|\nu^{(2)}|$.  We can iterate this process and produce an
  infinite sequence $(\nu^{(1)}, \nu^{(2)}, \dots)$ where each multi index
  satisfies $\nu^{(i)}\in\N$ and $\nu^{(i+1)}\idiv{L,\N}\nu^{(i)}+1_{j_i}$
  with $j_i\in\nmult{N}{L,\N}{\nu^{(i)}}$.  As $\N$ is a finite set, the
  elements of the sequence cannot be all different.  This contradicts our
  assumption that $L$ is a continuous division: by taking a sufficiently large
  part of this sequence we obtain a finite sequence with all properties
  mentioned in Definition~\ref{def:contdiv} but containing some identical
  elements.  Hence a multi index $\nu\in\N$ must exist such that
  $\sigma=\nu+1_j$.\qed
\end{proof}

\begin{lemma}\label{lem:cont}
  The Janet and the Pommaret division are continuous.
\end{lemma}

\begin{proof}
  Let $\N\subseteq\Nn$ be a finite set and $(\nu^{(i)},\dots,\nu^{(t)})$ a
  finite sequence where $\nu^{(i+1)}\idiv{L,\N}\nu^{(i)}+1_j$ with
  $j\in\nmult{N}{L,\N}{\nu^{(i)}}$ for $1\leq i<t$.
  
  We claim that for $L=J$, the Janet division,
  $\nu^{(i+1)}\succ_{\mbox{\scriptsize lex}}\nu^{(i)}$ implying that the
  sequence cannot contain any identical entries.  Set
  $k=\max\,\{i\mid\mu_i\neq\nu_i\}$.  Then $j\leq k$, as otherwise
  $j\in\mult{N}{J,\N}{\nu^{(i+1)}}$ entails $j\in\mult{N}{J,\N}{\nu^{(i)}}$
  contradicting our assumption that $j$ is non-multiplicative for the multi
  index $\nu^{(i)}$.  But $j<k$ is also not possible, as then
  $\nu^{(i+1)}_k<\nu^{(i)}_k$ and so $k$ cannot be multiplicative for
  $\nu^{(i+1)}$.  There remains as only possibility $j=k$.  In this case
  $\nu^{(i+1)}_j=\nu^{(i)}_j+1$, as otherwise $j$ could not be multiplicative
  for $\nu^{(i+1)}$.  Thus we conclude that
  $\nu^{(i+1)}\succ_{\mbox{\scriptsize lex}}\nu^{(i)}$ and the Janet division
  is continuous.
  
  The proof for the case $L=P$, the Pommaret division, is slightly more
  subtle.\footnote{It is tempting to tackle the Pommaret division in the same
    manner as the Janet division using $\prec_{\mbox{\scriptsize revlex}}$
    instead of $\prec_{\mbox{\scriptsize lex}}$; in fact, such a ``proof'' can
    be found in the literature.  Unfortunately, it is not correct: if
    $\nu^{(i+1)}=\nu^{(i)}+1_j$, then $\nu^{(i+1)}\prec_{\mbox{\scriptsize
        revlex}}\nu^{(i)}$ although the latter multi index is a divisor of the
    former one ($\prec_{\mbox{\scriptsize revlex}}$ is \emph{not} a term
    order!).  Thus the sequences considered in the application of
    Definition~\ref{def:contdiv} to the Pommaret division are in general not
    strictly ascending with respect to $\prec_{\mbox{\scriptsize revlex}}$.}
  The condition $j\in\nmult{N}{P}{\nu^{(i)}}$ implies that
  $\cls{(\nu^{(i)}+1_j)}=\cls{\nu^{(i)}}$ and if
  $\nu^{(i+1)}\idiv{P}\nu^{(i)}+1_j$, then
  $\cls{\nu^{(i+1)}}\geq\cls{\nu^{(i)}}$, i.\,e.\ the class of the elements of
  the sequence is monotonously increasing.  If
  $\cls{\nu^{(i+1)}}=\cls{\nu^{(i)}}=k$, then the involutive divisibility
  requires that $\nu^{(i+1)}_k\leq\nu^{(i)}_k$, i.\,e.\ among the elements of
  the sequence of the same class the corresponding entry is monotonously
  decreasing.  And if finally $\nu^{(i+1)}_k=\nu^{(i)}_k$, then we must have
  $\nu^{(i+1)}=\nu^{(i)}+1_j$, i.\,e.\ the length of the elements is strictly
  increasing.  Hence all elements of the sequence are different and the
  Pommaret division is continuous.\qed
\end{proof}

\begin{remark}\label{rem:prodbas}
  In Remark \ref{rem:sumbas} we discussed that for a global division a weak
  involutive basis of the sum $\I_1+\I_2$ of two monoid ideals is obtained by
  simply taking the union of (weak) involutive bases of $\I_1$ and $\I_2$.  As
  a more theoretical application of the concept of continuity, we prove now a
  similar statement for the product $\I_1\cdot\I_2$ and the intersection
  $\I_1\cap\I_2$ in the special case of the Pommaret division.  Let $\N_1$ be
  a (weak) Pommaret basis of $\I_1$ and $\N_2$ of $\I_2$.  We claim that the
  set $\N=\{\mu+\nu\mid\mu\in\N_1, \nu\in\N_2\}$ is a weak Pommaret basis of
  $\I_1\cdot\I_2$ and that the set $\hat{\N}=\{\lcm{(\mu,\nu)}\mid\mu\in\N_1,
  \nu\in\N_2\}$ is a weak Pommaret basis of $\I_1\cap\I_2$.
  
  By Proposition \ref{prop:localinv}, it suffices to show that the sets $\N$
  and $\hat{\N}$, respectively, are locally involutive for the Pommaret
  division.  Thus we take a generator $\mu+\nu\in\N$, where we assume for
  definiteness that $\cls{\mu}\leq\cls{\nu}$, and a non-multiplicative index
  $j_1>\cls{(\mu+\nu)}=\cls{\mu}$ of it.  Then $j_1$ is also
  non-multiplicative for $\mu\in\N_1$ alone and the Pommaret basis $\N_1$ must
  contain a multi index $\mu^{(1)}$ which involutively divides $\mu+1_{j_1}$.
  If we are lucky, then the generator $\mu^{(1)}+\nu\in\N$ is an involutive
  divisor of $\mu+\nu+1_{j_1}$, too, and we are done.
  
  Otherwise, there exists an index $k_1>\cls{\nu}$ such that
  $(\mu-\mu^{(1)})_{k_1}>0$.  In this case the Pommaret basis $\N_2$ must
  contain a multi index $\nu^{(1)}$ which involutively divides $\nu+1_{k_1}$.
  Again, if we are lucky, then $\mu^{(1)}+\nu^{(1)}\in\N$ is an involutive
  divisor of $\mu+\nu+1_{j_1}$ and we are done.  Otherwise, there are two
  possibilities.  There could be an index $j_2>\cls{\mu^{(1)}}$ such that
  $(\mu+\nu+1_{j_1}-\mu^{(1)}+\nu^{(1)})_{j_2}>0$ entailing the existence of a
  further generator $\mu^{(2)}\in\N_1$ which involutively divides
  $\mu^{(1)}+1_{j_2}$.  Or there could exist an index $k_2>\cls{\nu^{(1)}}$
  such that $(\mu+\nu+1_{j_1}-\mu^{(1)}+\nu^{(1)})_{k_2}>0$ implying that
  there is a multi index $\nu^{(2)}\in\N_2$ involutively dividing
  $\nu^{(1)}+1_{k_2}$.
  
  Continuing in this manner, one easily sees that we build up two sequences
  $\bigl(\mu,\mu^{(1)},\mu^{(2)},\dots\bigr)\subseteq\N_1$ and
  $\bigl(\nu,\nu^{(1)},\nu^{(2)},\dots\bigr)\subseteq\N_2$ as in the
  definition of a continuous division.  Since both Pommaret bases are finite
  by definition and the Pommaret division is continuous by Lemma
  \ref{lem:cont}, no sequence may become infinite and the above described
  process must stop with an involutive divisor of $\mu+\nu+1_{j_1}$.  Hence
  $\N$ is locally involutive and a weak Pommaret basis of $\I_1\cdot\I_2$.
  The proof for $\hat{\N}$ goes completely analogously replacing at
  appropriate places the sum of two multi indices by their least common
  multiple.  \bull
\end{remark}

\begin{definition}\label{def:constrdiv}
  Let\/ $L$ be a continuous involutive division and\/ $\N\subset\Nn$ a finite
  set of multi indices.  Choose a multi index\/ $\nu\in\N$ and a
  non-multiplicative index\/ $j\in\nmult{N}{L,\N}{\nu}$ such that:
  \begin{description}
  \item[{\upshape (i)}] $\nu+1_j\notin\ispan{\N}{L}$;
  \item[{\upshape (ii)}] if there exists\/ $\mu\in\N$ and\/
    $k\in\nmult{N}{L,\N}{\mu}$ such that $\mu+1_k\mid\nu+1_j$ but\/
    $\mu+1_k\neq\nu+1_j$, then\/ $\mu+1_k\in\ispan{\N}{L}$.
  \end{description}
  The division\/ $L$ is\/ \emph{constructive}, if for any such set\/ $\N$ and
  any such multi index\/ $\nu+1_j$ no multi index\/ $\rho\in\ispan{\N}{L}$
  with\/ $\nu+1_j\in\cone{L,\N\cup\{\rho\}}{\rho}$ exists.
\end{definition}

In words, constructivity may roughly be explained as follows.  The conditions
imposed on $\nu$ and $j$ ensure a kind of minimality: no proper divisor of
$\nu+1_j$ is of the form $\mu+1_k$ for a $\mu\in\N$ and not contained in the
involutive span $\ispan{\N}{L}$.  The conclusion implies that it is useless to
add multi indices to $\N$ that lie in some involutive cone, as none of them
can be an involutive divisor of $\nu+1_j$.  An efficient completion algorithm
for a constructive division should consider only non-multiplicative indices.

\begin{lemma}\label{lem:constr} 
  Any globally defined division (and thus the Pommaret division) is
  constructive.  The Janet division is constructive, too.
\end{lemma}

\begin{proof}
  For a globally defined division the proof is very simple.  For any multi
  index $\rho\in\ispan{\N}{L}$ there exists a multi index $\mu\in\N$ such that
  $\rho\in\cone{L}{\mu}$.  As for a globally defined division the
  multiplicative indices are independent of the reference set, we must have by
  the definition of an involutive division that
  $\cone{L}{\rho}\subseteq\cone{L}{\mu}$.  Hence adding such a multi index to
  $\N$ cannot change the involutive span and if $\nu+1_j\notin\ispan{\N}{L}$,
  then also $\nu+1_j\notin\ispan{\N\cup\{\rho\}}{L}$.  This implies
  constructivity.
  
  The proof of the constructivity of the Janet division is more involved.  The
  basic idea is to show that if it was not constructive, it could not be
  continuous either.  Let $\N$, $\nu$, $j$ be as described in
  Definition~\ref{def:constrdiv}.  Assume for a contradiction that a multi
  index $\rho\in\ispan{\N}{J}$ exists with
  $\nu+1_j\in\cone{J,\N\cup\{\rho\}}{\rho}$.  We write $\rho=\nu^{(1)}+\mu$
  for a multi index $\nu^{(1)}\in\N$ with $\rho\in\cone{J,\N}{\nu^{(1)}}$.  As
  $\nu+1_j\notin\ispan{\N}{J}$, we must have $|\mu|>0$.  Set
  $\lambda=\nu+1_j-\rho$ and let $m$, $l$ be the maximal indices such that
  $\mu_m>0$ and $\lambda_l>0$, respectively.
  
  We claim that $j>\max\,\{m,l\}$.  Indeed, if $j\leq m$, then
  $\nu^{(1)}_m<\nu_m$ and, by definition of the Janet division, this implies
  that $m\notin\mult{N}{J,\N}{\nu^{(1)}}$, a contradiction.  Similarly, we
  cannot have $j<l$, as then $l\notin\mult{N}{J,\N\cup\{\rho\}}{\rho}$.
  Finally, $j=l$ is not possible.  As we know already that $j>m$, we have in
  this case that $\rho_i=\nu^{(1)}_i=\nu_i$ for all $i>j$ and
  $\rho_j\leq\nu_j$.  Hence $j\in\nmult{N}{J,\N\cup\{\rho\}}{\nu}$ and this
  implies furthermore $j\in\nmult{N}{J,\N\cup\{\rho\}}{\rho}$, a
  contradiction.
  
  We construct a sequence as in Definition~\ref{def:contdiv} of a continuous
  division.  Choose an index $j_1$ with $\lambda_{j_1}>0$ and
  $j_1\in\nmult{N}{J,\N}{\nu^{(1)}}$.  Such an index exists, as otherwise
  $\nu+1_j\in\cone{J,\N}{\nu^{(1)}}\subseteq\ispan{\N}{J}$.  We write
  $\nu+1_j=(\nu^{(1)}+1_{j_1})+\mu+\lambda-1_{j_1}$.  Because of $|\mu|>0$,
  the multi index $\nu^{(1)}+1_{j_1}$ is a proper divisor of $\nu+1_j$ and
  according to our assumptions $\nu^{(2)}\in\N$ exists with
  $\nu^{(1)}+1_{j_1}\in\cone{J,\N}{\nu^{(2)}}$.
  
  By the same arguments as above an index $j_2\in\nmult{N}{J,\N}{\nu^{(2)}}$
  must exist with $(\mu+\lambda-1_{j_1})_{j_2}>0$ and a multi index
  $\nu^{(3)}\in\N$ with $\nu^{(2)}+1_{j_2}\in\cone{J,\N}{\nu^{(3)}}$.  Thus we
  can iterate and produce an infinite sequence $(\nu^{(1)},\nu^{(2)},\dots)$
  such that everywhere $\nu^{(i+1)}\idiv{J,\N}\nu^{(i)}+1_{j_i}$ with
  $j_i\in\nmult{N}{J,\N}{\nu^{(i)}}$.  By the continuity of the Janet division
  all members of the sequence must be different.  However, every multi index
  $\nu^{(i)}$ is a divisor of $\nu+1_j$, so only finitely many of them can be
  different.  Thus the sequence must terminate which only happens, if
  $\nu+1_j\in\cone{J,\N}{\nu^{(i)}}$ for some $i$ contradicting our
  assumptions.  \qed
\end{proof}

We present now an algorithm for determining weak involutive completions of
finite sets $\N\subset\Nn$.  As mentioned above, for arbitrary involutive
divisions, nobody has so far been able to find a reasonable approach.  But if
we assume that the division is constructive, then a very simple completion
algorithm exists, the basic ideas of which go back to Janet.

\begin{algorithm}
  \caption{Completion in $(\Nn,+)$\label{alg:monocompl}}
  \begin{algorithmic}[1]
    \REQUIRE a finite set $\N\subset\Nn$, an involutive division $L$
    \ENSURE a weak involutive completion $\bar\N$ of $\N$
    \STATE $\bar\N\leftarrow\N$
    \LOOP
        \STATE $\S\leftarrow\left\{\nu+1_j\mid\nu\in\bar\N,\,
                                   j\in\nmult{N}{L,\bar\N}{\nu},\,
                                   \nu+1_j\notin\ispan{\bar\N}{L}\right\}$
        \IF{$\S=\emptyset$}
            \STATE \algorithmicreturn $\bar\N$
        \ELSE
            \STATE choose $\mu\in\S$ such that $\S$ does not contain a proper 
                   divisor of it
            \STATE $\bar\N\leftarrow\bar\N\cup\{\mu\}$
        \ENDIF
    \ENDLOOP
  \end{algorithmic}
\end{algorithm}

The strategy behind Algorithm~\ref{alg:monocompl} is fairly natural given the
results above.  It collects in a set $\S$ all obstructions to local
involution.  For a continuous division $L$, the set $\N$ is weakly involutive,
if and only if $\S=\emptyset$.  Furthermore, for a constructive division $L$
it does not make sense to add elements of $\ispan{\N}{L}$ to $\N$ in order to
complete.  Thus we add in Line /8/ an element of $\S$ which is minimal in the
sense that the set $\S$ does not contain a proper divisor of it.

\begin{proposition}\label{prop:monocompl}
  Let the finite set\/ $\N\subset\Nn$ possess a finite (weak) involutive
  completion with respect to the constructive division\/ $L$.  Then
  Algorithm~\ref{alg:monocompl} terminates with a weak involutive completion\/
  $\bar\N$ of\/ $\N$.
\end{proposition}

\begin{proof}
  If Algorithm~\ref{alg:monocompl} terminates, its correctness is obvious
  under the made assumptions.  The criterion for its termination,
  $\S=\emptyset$, is equivalent to local involution of $\bar\N$.  By
  Proposition~\ref{prop:localinv}, local involution implies for a continuous
  division weak involution.  Thus the result $\bar\N$ is a weak involutive
  completion of $\N$, as by construction $\N\subseteq\bar\N\subset\lspan{\N}$.
  
  If the input set $\N$ is already involutive, Algorithm \ref{alg:monocompl}
  leaves it unchanged and thus obviously terminates.  Let us assume that $\N$
  is not yet involutive.  In the first iteration of the \texttt{loop} a multi
  index of the form $\mu=\nu+1_j$ is added to $\N$.  It is not contained in
  $\ispan{\N}{L}$ and $\S$ does not contain a proper divisor of it.  If $\N_L$
  is an arbitrary involutive completion of $\N$, it must contain a multi index
  $\lambda\notin\N$ such that $\lambda\idiv{L,\N_L}\mu$.  We claim that
  $\lambda=\mu$.
  
  Assume on the contrary that $\lambda\neq\mu$.  Since
  $\N_L\subset\lspan{\N}$, the multi index $\lambda$ must lie in the cone of a
  generator $\nu^{(1)}\in\N$.  We will show that, because of the continuity of
  $L$, $\lambda\in\ispan{\N}{L}$, contradicting the constructivity of $L$.  If
  $\nu^{(1)}\idiv{L,\N}\lambda$, we are done.  Otherwise we write
  $\lambda=\nu^{(1)}+\rho^{(1)}$ for some multi index $\rho^{(1)}\in\Nn$.  By
  construction, a non-multiplicative index $j_1\in\nmult{N}{L,\N}{\nu^{(1)}}$
  with $\rho^{(1)}_{j_1}>0$ must exist.  Consider the multi index
  $\nu^{(1)}+1_{j_1}$.  Because of $\nu^{(1)}+1_{j_1}\mid\lambda$, the multi
  index $\nu^{(1)}+1_{j_1}$ is a proper divisor of $\mu$.  Since the set $\S$
  does not contain any proper divisor of $\mu$, we must have
  $\nu^{(1)}+1_{j_1}\in\ispan{\N}{L}$.  Thus a multi index $\nu^{(2)}\in\N$
  exists such that $\nu^{(2)}\idiv{L,\N}\nu^{(1)}+1_{j_1}$.
  
  By iteration of this argument, we obtain a sequence $\bigl(\nu^{(1)},
  \nu^{(2)}, \dots\bigr)$ where each element $\nu^{(i)}\in\N$ is a divisor of
  $\lambda$ and where $\nu^{(i+1)}\idiv{L,\N}\nu^{(i)}+1_{j_i}$ with a
  non-multiplicative index $j_i\in\nmult{N}{L,\N}{\nu^{(i)}}$.  This sequence
  cannot become infinite for a continuous division, as $\lambda$ possesses
  only finitely many different divisors and all the multi indices $\nu^{(i)}$
  in arbitrary finite pieces of the sequence must be different.  But the
  sequence will only stop, if some $\nu^{(i)}\in\N$ exists such that
  $\nu^{(i)}\idiv{L,\N}\lambda$ and hence we must have that
  $\lambda\in\ispan{\N}{L}$.
    
  Thus \emph{every} weak involutive completion $\N_L$ of the given set $\N$
  must contain the multi index $\nu+1_j$.  In the next iteration of the
  \texttt{loop}, Algorithm~\ref{alg:monocompl} treats the enlarged set
  $\N_1=\N\cup\{\nu+1_j\}$.  It follows from our considerations above that any
  weak involutive completion $\N_L$ of $\N$ is also a weak involutive
  completion of $\N_1$ and hence we may apply the same argument again.  As a
  completion $\N_L$ is by definition a finite set, we must reach after a
  finite number $k$ of iterations a weak involutive basis $\N_k$ of
  $\lspan{\N}$.\qed
\end{proof}

Note the crucial difference between this result and the termination proof of
Buchberger's algorithm for the construction of Gr\"obner bases.  In the latter
case, we can show the termination for arbitrary input, i.\,e.\ the theorem
provides a constructive proof for the existence of such a basis.  Here we are
only able to prove the termination under the assumption that a finite (weak)
involutive basis exists; the existence has to be shown separately.  For
example, Lemma~\ref{lem:jannoeth} guarantees us that any monoid ideal
possesses a finite weak Janet basis.

Recall that by Proposition~\ref{prop:strongbas} any weak involutive basis can
be made strongly involutive by simply eliminating some redundant elements.
Thus we obtain an algorithm for the construction of a strong involutive basis
of $\lspan{\N}$ by adding an involutive autoreduction as last step to
Algorithm~\ref{alg:monocompl}.  Alternatively, we could perform the involutive
autoreduction as first step.  Indeed, if the input set $\N$ is involutively
autoreduced, then all intermediate sets $\bar\N$ constructed by
Algorithm~\ref{alg:monocompl} are also involutively autoreduced.  This is a
simple consequence of the second condition in Definition~\ref{def:invdiv} of
an involutive division that involutive cones may only shrink, if we add
elements to the set $\N$.

\begin{remark}\label{rem:complautored}\footnote{The following considerations
    are joint work with Vladimir Gerdt.}  While we just stated that it
  suffices to perform an involutive autoreduction as either first or last step
  in Algorithm~\ref{alg:monocompl}, we now analyse for later use what happens,
  if we involutively autoreduce $\bar\N$ every time a new element has been
  added to it.  The termination argument given in the proof of
  Proposition~\ref{prop:monocompl} does not remain valid after this
  modification and we must provide an alternative proof.
  
  Let again $\N_L=\bigl\{\mu^{(1)},\dots,\mu^{(r)}\bigr\}$ be a weak
  involutive completion of the input set $\N$.  If we denote by $\bar\N_i$ the
  value of $\bar\N$ after the $i$th iteration of the \texttt{loop}, then it
  was shown in the proof of Proposition~\ref{prop:monocompl} that $\N_L$ is
  also a weak involutive completion of any set $\bar\N_i$.  As by definition
  $\N_L$ is finite and each $\bar\N_i$ is a subset of it, the only possibility
  for non-termination is the appearance of a cycle, i.\,e.\ the existence of
  values $k_0$, $\ell$ such that $\bar\N_{k+\ell}=\bar\N_k$ for all $k\geq
  k_0$.
  
  Assume that in some iteration of the \texttt{loop} the multi index
  $\mu^{(k)}$ is added to $\bar\N$ and that in the subsequent involutive
  autoreduction some elements of $\bar\N$ are eliminated (in order to have a
  cycle this must indeed happen).  The first step in the autoreduction must be
  that some multi index $\mu^{(\ell)}$ is eliminated, because $\mu^{(k)}$ is
  an involutive divisor of it.  Indeed, by Condition (ii) in
  Definition~\ref{def:invdiv}, any other reduction would have been possible
  already before the insertion of $\mu^{(k)}$ and thus the previous involutive
  autoreduction would not have been finished.
  
  Since $\mu^{(k)}$ has been added to $\bar\N$, there must exist some multi
  index $\mu^{(a_1)}\in\N$ such that $\mu^{(k)}=\mu^{(a_1)}+\rho$.
  Furthermore, we know that $\mu^{(\ell)}=\mu^{(k)}+\tilde\sigma$ for some
  multi index $\tilde\sigma$ with $|\tilde\sigma|>0$ and thus
  $\mu^{(\ell)}=\mu^{(a_1)}+\sigma$ with $\sigma=\tilde\sigma+\rho$ and
  $|\sigma|>1$.  As we are in a cycle, the multi index $\mu^{(\ell)}$ must
  have been added to $\bar\N$ in a previous iteration of the \texttt{loop},
  say when analysing $\bar\N_i$.  Thus $\mu^{(\ell)}$ cannot be involutively
  divisible by $\mu^{(a_1)}$ and we must have $\sigma_{j_1}>0$ for a
  non-multiplicative index $j_1\in\nmult{N}{L,\bar\N_i}{\mu^{(a_1)}}$.  It
  cannot be that $\mu^{(a_1)}+1_{j_1}=\mu^{(\ell)}$, as $|\sigma|>1$, and
  therefore $\mu^{(a_1)}+1_{j_1}$ is a proper divisor of $\mu^{(\ell)}$.
  Hence $\bar\N_i$ must contain an involutive divisor $\mu^{(a_2)}$ of
  $\mu^{(a_1)}+1_{j_1}$, as otherwise this multi index would have been added
  to $\bar\N$ instead of $\mu^{(\ell)}$.
  
  Obviously, $\mu^{(a_2)}\mid\mu^{(k)}$ and, decomposing
  $\mu^{(k)}=\mu^{(a_2)}+\pi$, we conclude by the same reasoning as above that
  $\pi_{j_2}>0$ for some non-multiplicative index
  $j_2\in\nmult{N}{L,\bar\N_i}{\mu^{(a_2)}}$.  Iteration of this argument
  yields an infinite sequence $\bigl(\mu^{(a_1)},\mu^{(a_2)},\dots\bigr)$ as
  in Definition~\ref{def:contdiv} of a continuous division.  However, since
  $L$ is a continuous division and $\N_L$ a finite set, we arrive at a
  contradiction.  Thus even with involutive autoreductions after each step
  Algorithm~\ref{alg:monocompl} terminates.\bull
\end{remark}

In some sense our description of Algorithm~\ref{alg:monocompl} is not
complete, as we have not specified how one should choose the multi index $\mu$
in Line /7/, if several choices are possible.  One would expect that different
involutive completions are obtained for different choices.  However, an
interesting aspect of our proof of Proposition~\ref{prop:monocompl} is that it
shows that this is not the case.  The choice affects only the order in which
multi indices are added but not which or how many multi indices are added
during the completion.  A simple method for choosing $\mu$ consists of taking
an arbitrary term order $\prec$ (which also could be changed in each iteration
of the \texttt{loop}) and setting $\mu=\min_{\prec}{\S}$.

\begin{corollary}\label{cor:unicompl}
  If Algorithm~\ref{alg:monocompl} terminates, its output\/ $\bar\N$ is
  independent of the manner in which\/ $\mu$ is chosen.  Furthermore, if\/
  $\N_L$ is any weak involutive completion of\/ $\N$ with respect to the
  division\/ $L$, then $\bar\N\subseteq\N_L$.
\end{corollary}

\begin{proof}
  Consider the set $\L(\N)$ of all weak involutive completions of $\N$ with
  respect to the division $L$ and define
  \begin{equation}\label{eq:Nmin}
    \tilde\N=\bigcap_{\N_L\in\L(\N)}\N_L\;.
  \end{equation}
  We claim that Algorithm~\ref{alg:monocompl} determines this set $\tilde\N$
  independent of the used term order.  Obviously, this implies our corollary.
  
  In the proof of Proposition~\ref{prop:monocompl} we showed that the multi
  indices added in Algorithm~\ref{alg:monocompl} are contained in \emph{every}
  weak involutive completion of $\N$.  Thus all these multi indices are
  elements of $\tilde\N$.  As our algorithm terminates with a weak involutive
  completion, its output is also an element of $\L(\N)$ and hence must be
  $\tilde\N$.\qed
\end{proof}

Any monoid ideal in $\Nn$ has a \emph{unique} minimal basis: take an arbitrary
basis and eliminate all multi indices having a divisor in the basis.
Obviously, these eliminations do not change the span and the result is a
minimal basis.  Similarly we have seen in Section~\ref{sec:invdiv} that if a
monoid ideal $\I\subseteq\Nn$ has a finite involutive basis for a given
division $L$, then a unique minimal involutive basis exists.  By the same
argument as in the proof of Corollary~\ref{cor:unicompl}, it can easily be
constructed by taking the unique minimal basis of $\I$ as input for
Algorithm~\ref{alg:monocompl}.

\section{Polynomial Completion}\label{sec:polycompl}

An obvious way to compute an involutive basis for an ideal $\I$ in a
polynomial algebra $(\P,\star,\prec)$ of solvable type goes as follows: we
determine first a Gr\"obner basis $\G$ of $\I$ and then with
Algorithm~\ref{alg:monocompl} an involutive completion of $\le{\prec}{\G}$.
In fact, a similar method is proposed by Sturmfels and White \cite{sw:comb}
for the construction of Stanley decompositions (cf.~Part~II).  However, we
prefer to extend the ideas behind Algorithm~\ref{alg:monocompl} to a direct
completion algorithm for polynomial ideals, as we believe that this approach
is more efficient.

First, we need two subalgorithms: \emph{involutive normal forms} and
\emph{involutive head autoreductions}.  The design of an algorithm
$\mathtt{NormalForm}_{L,\prec}(g,\H)$ determining an involutive normal form of
the polynomial $g$ with respect to the finite set $\H\subset\P$ is trivial.
We may use the standard algorithm for normal forms in the Gr\"obner theory, if
we replace the ordinary divisibility by its involutive version.  Obviously,
this does not affect the termination.  Actually, for our purposes it is not
even necessary to compute a full normal form; we may stop as soon as we have
obtained a polynomial that is not involutively head reducible.

The design of an algorithm $\mathtt{InvHeadAutoReduce}_{L,\prec}(\F)$ for an
involutive head autoreduction of a finite set $\F$ is also obvious.  Again one
may use the standard algorithm with the ordinary divisibility replaced by
involutive divisibility.

Based on these two subalgorithms, we propose Algorithm~\ref{alg:polycompl} for
the computation of involutive bases in $\P$.  It follows the same strategy as
the monomial algorithm.  We multiply each generator by its non-multiplicative
variables.  Then we decide whether or not the result is already contained in
the involutive span of the basis; if not, it is added.  This decision is
effectively made via an involutive normal form computation: the involutive
normal form of a polynomial is zero, if and only if the polynomial lies in the
involutive span.  As our goal is a strong involutive basis, we take care that
our set is always involutively head autoreduced.

\begin{algorithm}
  \caption{Completion in $(\P,\star,\prec)$\label{alg:polycompl}}
  \begin{algorithmic}[1]
    \REQUIRE a finite set $\F\subset\P$, an involutive division $L$
    \ENSURE an involutive basis $\H$ of $\I=\lspan{\F}$ with respect to $L$
    and $\prec$
    \STATE $\H\leftarrow\mathtt{InvHeadAutoReduce}_{L,\prec}(\F)$
    \LOOP
        \STATE $\S\leftarrow\left\{x_j\star h\mid h\in\H,\,
                                   x_j\in\nmult{X}{L,\H,\prec}{h},\,
                                   x_j\star h\notin\ispan{\H}{L,\prec}\right\}$
        \IF{$\S=\emptyset$}
            \STATE \algorithmicreturn $\H$
        \ELSE
            \STATE choose $\bar g\in\S$ such that 
                   $\le{\prec}{\bar g}=\min_\prec\S$
            \STATE $g\leftarrow\mathtt{NormalForm}_{L,\prec}(\bar g,\H)$
            \STATE $\H\leftarrow
                    \mathtt{InvHeadAutoReduce}_{L,\prec}(\H\cup\{g\})$
        \ENDIF
    \ENDLOOP
  \end{algorithmic}
\end{algorithm}

The manner in which we choose in Line /7/ the next polynomial $\bar g$ to be
treated (we briefly write $\min_{\prec}\S$ for the minimal leading exponent of
an element of $\S$) corresponds to the normal selection strategy in the theory
of Gr\"obner bases.  There, this strategy is known to work very well for
degree compatible term orders but not so well for other orders like the purely
lexicographic one.  Whereas for Gr\"obner bases the selection strategy
concerns only the efficiency of the computation, we will see below that here
the use of this particular strategy is important for our termination proof.
With more refined and optimised versions of the basic completion Algorithm
\ref{alg:polycompl} one can circumvent this restriction
\cite{ah:unred,cg:invdir,vpg:opt}, but we will not discuss this highly
technical question here.

\begin{definition}\label{def:localinvpoly}
  A finite set\/ $\F\subset\P$ is\/ \emph{locally involutive} for the
  division\/~$L$, if for every polynomial\/ $f\in\F$ and for every
  non-multiplicative variable\/ $x_j\in\nmult{X}{L,\F,\prec}{f}$ the product\/
  $x_j\star f$ has an involutive standard representation with respect to\/
  $\F$.
\end{definition}

Note that for an involutively head autoreduced set $\F$, we may equivalently
demand that $x_j\star f\in\ispan{\F}{L,\prec}$; because of
Lemma~\ref{lem:head} this automatically implies the existence of an involutive
standard representation.  In fact, the criterion appears in this form in Line
/3/ of Algorithm~\ref{alg:polycompl}.  In any case, local involution may be
effectively verified by computing an involutive normal form of $x_j\star f$ in
the usual manner, i.\,e.\ always performing head reductions.

\begin{proposition}\label{prop:localinvpoly}
  If the finite set\/ $\F\subset\P$ is locally involutive for the continuous
  division\/~$L$, then\/ $\ispan{\F}{L,\prec}=\lspan{\F}$.
\end{proposition}

\begin{proof}
  We claim that if the set $\F$ is locally involutive (with respect to the
  continuous division $L$), then every product $x^\mu\star f_1$ of an
  arbitrary term $x^\mu$ with a polynomial $f_1\in\F$ possesses an involutive
  standard representation.  This claim trivially entails our proposition, as
  any polynomial in $\lspan{\F}$ consists of a finite linear combination of
  such products: adding the corresponding involutive standard representations
  shows that the polynomial is contained in $\ispan{\F}{L,\prec}$.
  
  In order to prove our claim, it suffices to show the existence of a
  representation
  \begin{equation}\label{eq:involhead}
    x^\mu\star f_1=\sum_{f\in\F}\Bigl(P_f\star f+
      \sum_{\nu\in\Nn}c_{\nu,f}x^\nu\star f\Bigr)
  \end{equation}
  where $P_f\in\kk[\mult{X}{L,\F,\prec}{f}]$ and $\le{\prec}{(P_f\star
    f)}=\le{\prec}{(x^\mu\star f_1)}$ (or $P_f=0$) and where the coefficients
  $c_{\nu,f}\in\kk$ vanish for all multi indices $\nu\in\Nn$ such that
  $\le{\prec}{(x^\nu\star f)}\succeq\le{\prec}{(x^\mu\star f_1)}$.  Our claim
  follows then by an obvious induction.
  
  If $x^\mu\in\kk[\mult{X}{L,\F,\prec}{f_1}]$, i.\,e.\ it contains only
  variables that are multiplicative for $\le{\prec}{f_1}$, nothing has to be
  shown.  Otherwise we choose a non-multiplicative index
  $j_1\in\nmult{N}{L,\le{\prec}{\F}}{\le{\prec}{f_1}}$ such that
  $\mu_{j_1}>0$.  As $\F$ is locally involutive, an involutive standard
  representation $x_{j_1}\star f_1=\sum_{f\in\F}P_f^{(1)}\star f$ exists.  Let
  $\F_2\subseteq\F$ contain all polynomials $f_2$ such that
  $\le{\prec}{(P_{f_2}^{(1)}\star f_2)}=\le{\prec}{(x_{j_1}\star f_1)}$.  If
  we have $x^{\mu-1_{j_1}}\in\kk[\mult{X}{L,\F,\prec}{f_2}]$ for all
  polynomials $f_2\in\F_2$, then we are done, as at least
  $\lm{\prec}{(x^{\mu-1_{j_1}}\star
    P_{f_2}^{(1)})}\in\kk[\mult{X}{L,\F,\prec}{f_2}]$.
  
  Otherwise we consider the subset $\F_2'\subseteq\F_2$ of polynomials $f_2$
  for which $x^{\mu-1_{j_1}}\notin\kk[\mult{X}{L,\F,\prec}{f_2}]$ and iterate
  over it.  For each polynomial $f_2\in\F_2'$ we choose a non-multiplicative
  index $j_2\in\nmult{N}{L,\le{\prec}{\F}}{\le{\prec}{f_2}}$ such that
  $(\mu-1_{j_1})_{j_2}>0$.  Again the local involution of the set $\F$ implies
  the existence of an involutive standard representation $x_{j_2}\star
  f_2=\sum_{f\in\F}P_f^{(2)}\star f$.  We collect in $\F_3\subseteq\F$ all
  polynomials $f_3$ such that $\le{\prec}{(P_{f_3}^{(2)}\star
    f_3)}=\le{\prec}{(x_{j_2}\star f_2)}$.  If we introduce the multi index
  $\nu=\le{\prec}{(x_{j_1}\star f_1)}-\le{\prec}{f_2}$, then
  $\le{\prec}{(x^\mu\star f_1)}=\le{\prec}{(x^{\mu+\nu-1_{j_1}-1_{j_2}}\star
    f_3)}$ for all $f_3\in\F_3$.  If
  $x^{\mu+\nu-1_{j_1}-1_{j_2}}\in\kk[\mult{X}{L,\F,\prec}{f_3}]$ for all
  $f_3\in\F_3$, we are done.
  
  Otherwise we continue in the same manner: we collect in a subset
  $\F_3'\subseteq\F_3$ all polynomials $f_3$ which are multiplied by
  non-multiplicative variables, for each of them we choose a
  non-multiplicative index $j_3\in\kk[\mult{X}{L,\F,\prec}{f_3}]$ such that
  $(\mu-1_{j_1}-1_{j_2})_{j_3}>0$, determine an involutive standard
  representation of $x_{j_3}\star f_3$ and analyse the leading terms.  If they
  are still multiplied with non-multiplicative variables, this leads to sets
  $\F_4'\subseteq\F_4$ and so on.  This process yields a whole tree of cases
  and each branch leads to a sequence
  $\bigl(\nu^{(1)}=\le{\prec}{f_1},\nu^{(2)}=\le{\prec}{f_2},\dots\bigr)$
  where all contained multi indices $\nu^{(k)}$ are elements of the finite set
  $\le{\prec}{\F}$ and where to each $\nu^{(k)}$ a non-multiplicative index
  $j_k\in\nmult{N}{L,\le{\prec}{\F}}{\nu^{(k)}}$ exists such that
  $\nu^{(k+1)}\idiv{L,\le{\prec}{\F}}\nu^{(k)}+1_{j_k}$.  By the definition of
  a continuous division, this sequence cannot become infinite and thus each
  branch must terminate.  But this implies that we may construct for each
  polynomial $f_1\in\F$ and each non-multiplicative variables
  $x_j\in\nmult{X}{L,\F,\prec}{f_1}$ a representation of the claimed form
  (\ref{eq:involhead}).  \qed
\end{proof}

Note that the proposition only asserts that the involutive span equals the
normal span.  It does \emph{not} say that $\F$ is weakly involutive (indeed,
the set $\F$ studied in Example~\ref{ex:weakbasis} would be a simple
counterexample).  If $g=\sum_{\mu\in\Nn}\sum_{f\in\F} c_{\mu,f}x^\mu\star f$
is an arbitrary polynomial in $\lspan{\F}$, then adding the involutive
standard representations of all the products $x^\mu\star f$ for which
$c_{\mu,f}\neq0$ yields a representation $g=\sum_{f\in\F} P_f\star f$ where
$P_f\in\kk[\mult{X}{L,\F,\prec}{f}]$.  But in general it will not satisfy the
condition $\le{\prec}{(P_f\star f)}\preceq\le{\prec}{g}$ for all $f\in\F$, as
we cannot assume that we started with an ordinary standard representation of
$g$.  The satisfaction of this condition is guaranteed only for involutively
head autoreduced sets, as there it is impossible that the leading terms cancel
(Lemma~\ref{lem:head}).  For such sets the above proof simplifies, as all the
sets $\F_i$ consist of precisely one element and thus no branching is
necessary.

\begin{corollary}\label{cor:localinvpoly}
  For a continuous division\/ $L$ an involutively head autoreduced set\/
  $\F\subset\P$ is involutive, if and only if it is locally involutive.
\end{corollary}

As in the proof of Proposition~\ref{prop:monocompl}, local involution of $\H$
is obviously equivalent to the termination condition $\S=\emptyset$ of the
\texttt{loop} in Algorithm~\ref{alg:polycompl}.  Thus we are now in the
position to prove the following result.

\begin{theorem}\label{thm:polycompl}
  Let\/ $L$ be a constructive Noetherian involutive division and\/
  $(\P,\star,\prec)$ a polynomial algebra of solvable type.  Then
  Algorithm~\ref{alg:polycompl} terminates for any finite input set\/ $\F$
  with an involutive basis of the ideal\/ $\I=\lspan{\F}$.
\end{theorem}

\begin{proof}
  We begin by proving the \emph{correctness} of the algorithm under the
  assumption that it terminates.  The relation $\I=\lspan{\H}$ remains valid
  throughout, although $\H$ changes.  But the only changes are the addition of
  further elements of $\I$ and involutive head autoreductions; both operations
  do not affect the ideal generated by $\H$.  When the algorithm terminates,
  we have $\S=\emptyset$ and thus the output $\H$ is locally involutive and by
  Corollary~\ref{cor:localinvpoly} involutive.
  
  There remains the problem of \emph{termination}.
  Algorithm~\ref{alg:polycompl} produces a sequence $(\H_1,\H_2,\dots)$ with
  $\lspan{\H_i}=\I$.  The set $\H_{i+1}$ is determined from $\H_i$ in Line
  /9/.  We distinguish two cases, namely whether or not during the computation
  of the involutive normal form in Line /8/ the leading exponent changes.  If
  $\le{\prec}{\bar g}=\le{\prec}{g}$, then
  $\lspan{\le{\prec}{\H_i}}=\lspan{\le{\prec}{\H_{i+1}}}$, as
  $\le{\prec}{g}=\le{\prec}{h}+1_j$ for some $h\in\H_i$.  Otherwise we claim
  that $\lspan{\le{\prec}{\H_i}}\subsetneq\lspan{\le{\prec}{\H_{i+1}}}$.
  
  By construction, $g$ is in involutive normal form with respect to the set
  $\H_i$ implying that $\le{\prec}{g}\in\lspan{\le{\prec}{\H_i}}\setminus
  \ispan{\le{\prec}{\H_i}}{L}$.  If we had
  $\lspan{\le{\prec}{\H_i}}=\lspan{\le{\prec}{\H_{i+1}}}$, a polynomial
  $h\in\H_i$ would exist such that $\le{\prec}{g}=\le{\prec}{h}+\mu$ where the
  multi index $\mu$ has a non-vanishing entry $\mu_j$ for at least one
  non-multiplicative index $j\in\nmult{N}{L,\le{\prec}{\H_i}}{h}$.  This
  implies that $\le{\prec}{h}+1_j\preceq\le{\prec}{g}\prec\le{\prec}{\bar g}$.
  But we choose the polynomial $\bar g$ in Line /7/ such that its leading
  exponent is minimal among all non-multiplicative products $x_k\star h$ with
  $h\in\H_i$; hence $\le{\prec}{\bar g}\preceq\le{\prec}{h}+1_j$.  As this is
  a contradiction, we must have
  $\lspan{\le{\prec}{\H_i}}\subsetneq\lspan{\le{\prec}{\H_{i+1}}}$.
  
  So the \texttt{loop} of Algorithm~\ref{alg:polycompl} generates an ascending
  chain of monoid ideals
  $\lspan{\le{\prec}{\H_1}}\subseteq\lspan{\le{\prec}{\H_2}}\subseteq\cdots
  \subseteq\le{\prec}{\I}$.  As $\Nn$ is Noetherian, the chain must become
  stationary at some index $N$.  It follows from the considerations above that
  in all iterations of the \texttt{loop} after the $N$th one $\le{\prec}{\bar
    g}=\le{\prec}{g}$ in Line /8/.  At this stage
  Algorithm~\ref{alg:polycompl} reduces to an involutive completion of the
  monomial set $\le{\prec}{\H_N}$ using Algorithm~\ref{alg:monocompl}---but
  with additional involutive autoreductions after each appearance of a new
  element.  Indeed, in Line /7/ we choose the polynomial $\bar g$ such that
  $\le{\prec}{\bar g}$ is a possible choice for the multi index $\mu$
  Algorithm~\ref{alg:monocompl} adds in Line /8/.  Since we assume that our
  division is Noetherian, it follows now from Proposition~\ref{prop:monocompl}
  together with Remark~\ref{rem:complautored} that
  Algorithm~\ref{alg:polycompl} terminates (and our correctness proof above
  implies that in fact $\lspan{\le{\prec}{\H_N}}=\le{\prec}{\I}$).\qed
\end{proof}

\begin{remark}
  If the division $L$ is not Noetherian, then it may happen that, even when
  the ideal $\I=\lspan{\F}$ does possess a finite involutive basis with
  respect to $L$, Algorithm~\ref{alg:polycompl} does not terminate for the
  input $\F$.  We will see concrete examples for this phenomenon in Part~II
  for the Pommaret division.
  
  The problem is that the existence of an involutive basis for
  $\le{\prec}{\I}$ does not imply that all subideals of it have also an
  involutive basis (as a trivial counter example consider
  $\lspan{xy}\subset\lspan{xy,y^2}$ with the Pommaret division).  In such a
  case it may happen that at some stage of Algorithm~\ref{alg:polycompl} we
  encounter a basis $\H_i$ such that $\lspan{\le{\prec}{\H_i}}$ does not
  possess an involutive basis and then it is possible that the algorithm
  iterates endlessly in an attempt to complete $\le{\prec}{\H_i}$.
  
  This observation entails that variations of Theorem \ref{thm:polycompl} hold
  also for divisions which are not Noetherian.  For example, we could assume
  instead that all subideals of $\le{\prec}{\I}$ possess an involutive basis.
  Alternatively, we could restrict to term orders of type $\omega$.  Then it
  suffices to assume that $\le{\prec}{\I}$ has an involutive basis.  Indeed,
  now it is not possible that Algorithm~\ref{alg:polycompl} iterates endlessly
  within $\le{\prec}{\H_i}$, as sooner or later an element $\bar g$ must be
  selected in Line /7/ with $\le{\prec}{\bar g}\notin\le{\prec}{\H_i}$.\bull
\end{remark}

\begin{corollary}\label{cor:exinvbas}
  For a constructive Noetherian division\/ $L$ every ideal\/ $\I\subseteq\P$
  possesses a finite involutive basis.
\end{corollary}

\begin{example}
  Now we are finally in the position to prove the claims made in Example
  \ref{ex:obstr}.  With respect to the degree reverse lexicographic term order
  the Janet (and the Pommaret) division assigns the polynomial $f_1=z^2-xy$
  the multiplicative variables $\{x,y,z\}$ and the polynomials $f_2=yz-x$ and
  $f_3=y^2-z$ the multiplicative variables $\{x,y\}$.  Thus we must check the
  two non-multiplicative products: $zf_2=yf_1+xf_3$ and $zf_3=yf_2-f_1$.  As
  both possess an involutive standard representation, the set $\S$ in Line /3/
  of Algorithm \ref{alg:polycompl} is empty in the first iteration and thus
  $\F$ is a Janet (and a Pommaret) basis of the ideal it generates.
  
  The situation changes, if we use the degree inverse lexicographic term
  order, as then $\lt{\prec}{f_1}=xy$.  Now $\mult{X}{J,\F,\prec}{f_1}=\{x\}$,
  $\mult{X}{J,\F,\prec}{f_2}=\{x,y,z\}$ and
  $\mult{X}{J,\F,\prec}{f_3}=\{x,y\}$.  In the first iteration we find
  $\S=\{zf_1\}$.  Its involutive normal form is $f_4=z^3-x^2$ and we add this
  polynomial to $\F$ to obtain $\H_1=\{f_1,f_2,f_3,f_4\}$ (the involutive head
  autoreduction does not change the set).  For $f_4$ all variables are
  multiplicative; for the other generators there is one change: $z$ is no
  longer multiplicative for $f_2$.  Thus in the second iteration
  $\S=\{zf_2\}$.  It is easy to check that this polynomial is already in
  involutive normal form with respect to $\H_1$ and we obtain $\H_2$ by adding
  $f_5=yz^2-xz$ to $\H_1$.  In the next iteration $\S$ is empty, so that
  $\H_2$ is indeed the Janet basis of $\lspan{\F}$ for the degree inverse
  lexicographic term order.\bull
\end{example}

Our proof of Theorem~\ref{thm:polycompl} has an interesting consequence which
was first discovered by Apel \cite{apel:alter} for the special case of the
Pommaret division.  Assume that the term order $\prec$ is of type $\omega$,
i.\,e.\ for any two multi indices $\mu$, $\nu$ with $\mu\prec\nu$ only
finitely many multi indices $\rho^{(i)}$ exist with
$\mu\prec\rho^{(1)}\prec\rho^{(2)}\prec\cdots\prec\nu$.  Then even if our
algorithm does \emph{not} terminate, it determines in a finite number of steps
a Gr\"obner basis of the ideal $\I$.

\begin{proposition}\label{prop:infgb}
  Let\/ the term order\/ $\prec$ be of type\/ $\omega$.  Then
  Algorithm~\ref{alg:polycompl} determines for any finite input set\/
  $\F\subset\P$ in a finite number of steps a Gr\"obner basis of the ideal\/
  $\I=\lspan{\F}$.
\end{proposition}

\begin{proof}
  Above we introduced the set $\H_N$ such that
  $\lspan{\le{\prec}{\H_{N+\ell}}}=\lspan{\le{\prec}{\H_N}}$ for all $\ell>0$.
  We claim that $\H_N$ is a Gr\"obner basis of $\I$.
  
  Let $f\in\I$ be an arbitrary element of the ideal.  As $\H_N$ is a basis of
  $\I$, we find for each $h\in\H_N$ a polynomial $g_h\in\P$ such that
  \begin{equation}\label{eq:Hrep}
    f=\sum_{h\in\H_N}g_h\star h\;.
  \end{equation}
  $\H_N$ is a Gr\"obner basis, if and only if we can choose the coefficients
  $g_h$ such that $\le{\prec}{(g_h\star h)}\preceq\le{\prec}{f}$.  Assume that
  for $f$ no such standard representation exists and let
  $\mu=\max_{h\in\H_N}\bigl\{\le{\prec}{g_h}+\le{\prec}{h}\bigr\} \succ
  \le{\prec}{f}$.  If we denote by $\bar\H_N$ the set of all polynomials $\bar
  h\in\H_N$ for which $\le{\prec}{g_{\bar h}}+\le{\prec}{\bar h}=\mu$, then
  the identity $\sum_{\bar h\in\bar\H_N}\lc{\prec}{(g_{\bar h}\star\bar h)}=0$
  must hold and hence $\bar\H_N$ contains at least two elements.  For each
  element $\bar h\in\bar\H_N$ we have $\mu\in\cone{}{\le{\prec}{\bar h}}$.  As
  by construction the set $\H_N$ is involutively head autoreduced, the
  involutive cones of the leading exponents do not intersect and there must be
  at least one generator $\bar h\in\bar\H_N$ such that some non-multiplicative
  variable $x_j\in\nmult{X}{L,\H_N}{\bar h}$ divides $\lt{\prec}{g_{\bar h}}$.
  
  As $\prec$ is of type $\omega$, after a finite number of steps the
  non-multiplicative product $x_j\star\bar h$ is analysed in
  Algorithm~\ref{alg:polycompl}.  Thus for some $n_1\geq0$ the set
  $\H_{N+n_1}$ contains an element $\bar h'$ with $\le{\prec}{\bar
    h'}=\le{\prec}{(x_j\star\bar h)}$.  Let $\mu=\le{\prec}{g_{\bar h}}$,
  $x^{\mu-1_j}\star x_j=cx^\mu+r_1$ and $\bar h'=dx_j\star \bar h+r_2$.  Then
  we may rewrite
  \begin{equation}\label{eq:rewrite}
    g_{\bar h}\star\bar h=
        \frac{\lc{\prec}{g_{\bar h}}}{cd}
            \Bigl[x^{\mu-1_j}\star(\bar h'-r_2)-dr_1\star\bar h\Bigr]
        +\bigl(g_{\bar h}-\lm{\prec}{g_{\bar h}}\bigr)\star\bar h\;.
  \end{equation}
  As $\bar h'$ was determined via an involutive normal form computation
  applied to the product $x_j\star\bar h$ and as we know that at this stage of
  the algorithm the leading exponent does not change during the computation,
  the leading exponent on the right hand side of (\ref{eq:rewrite}) is
  $\le{\prec}{(x^{\mu-1_j}\star\bar h')}$.  If the term $x^{\mu-1_j}$ contains
  a non-multiplicative variable $x_k\in\nmult{X}{L,\H_{N+n_1}}{\bar h'}$, we
  repeat the argument obtaining a polynomial $\bar h''\in\H_{N+n_1+n_2}$ such
  that $\le{\prec}{\bar h''}=\le{\prec}{(x_k\star\bar h')}$.
  
  Obviously, this process terminates after a finite number of steps, even if
  we do it for each $\bar h\in\bar\H_N$.  Thus after $\ell$ further iterations
  we obtain a set $\H_{N+\ell}$ such that, after applying all the found
  relations (\ref{eq:rewrite}), $f$ can be expressed in the form
  $f=\sum_{h\in\H_{N+\ell}}\tilde g_h\star h$ where still
  $\mu=\max_{h\in\H_{N+\ell}}\bigl\{\le{\prec}{\tilde
    g_h}+\le{\prec}{h}\bigr\}$.  Denote again by
  $\bar\H_{N+\ell}\subseteq\H_{N+\ell}$ the set of all polynomials $\bar h$
  achieving this maximum.
  
  By construction, no term $\lt{\prec}{\tilde g_{\bar h}}$ with $\bar
  h\in\bar\H_{N+\ell}$ contains a variable that is non-multiplicative for
  $\bar h$.  Thus we must now have
  $\mu\in\cone{\le{\prec}{(\H_{N+\ell})},L}{\le{\prec}{\bar h}}$ for each
  $\bar h\in\bar\H_{N+\ell}$ implying that $\bar\H_{N+\ell}$ contains at most
  one element.  But then it is not possible that $\mu\succ\le{\prec}{f}$.
  Hence each polynomial $f\in\P$ possesses a standard representation already
  with respect to $\H_N$ and this set is a Gr\"obner basis.\qed
\end{proof}

Note that in the given form this result is only of theoretical interest, as in
general no efficient method exists for checking whether the current basis is
already a Gr\"obner basis.  Using standard criteria would destroy all
potential advantages of the involutive algorithm.  For the special case of
Pommaret bases, Apel \cite{apel:alter} found a simple criterion that allows us
to use a variant of Algorithm \ref{alg:polycompl} for the construction of
Gr\"obner bases independent of the existence of a finite involutive basis.

In contrast to the monomial case, one does not automatically obtain a minimal
involutive basis by making some minor modifications of
Algorithm~\ref{alg:polycompl}.  In particular, it does not suffice to apply it
to a minimal basis in the ordinary sense.  Gerdt and Blinkov \cite{gb:minbas}
presented an algorithm that always returns a minimal involutive basis provided
a finite involutive basis exists.  While it still follows the same basic
strategy of study all products with non-multiplicative variables, it requires
a more complicated organisation of the algorithm.  We omit here the details.

\section{Right and Two-Sided Bases}\label{sec:two}

We now briefly discuss the relation between left and right involutive bases
and the computation of bases for two-sided ideals.  We use in this section the
following notations: the left ideal generated by $\F\subset\P$ is denoted by
$\lspanl{\F}$, the right ideal by $\lspanr{\F}$ and the two-sided ideal by
$\lspant{\F}$ and corresponding notations for the left, right and two-sided
involutive span.

Recall from Remark~\ref{rem:rightnoether} that even with a coefficient field
$\kk$ it is not guaranteed that $\P$ is also right Noetherian and hence
generally the existence of right Gr\"obner bases for right ideals is not
clear.  However, we also noted that the ring $\P$ is always right Noetherian,
if we assume that the maps $\rho_{i}:\kk\rightarrow\kk$ in (\ref{eq:crop}) are
automorphisms.  In the sequel of this section we will always make this
assumption.

From a computational point of view, the theory of right ideals is almost
identical to the corresponding theory for left ideals.  The left-right
asymmetry in our definition of polynomial algebras of solvable type leads
only to one complication.  Suppose that we want to perform a right reduction
of a term $ax^\nu$ with respect to another term $cx^\mu$ with $\mu\mid\nu$.
This requires to find a coefficient $b\in\kk$ such that
$\lc{\prec}{(cx^\mu\star bx^{\nu-\mu})}=c\rho_\mu(b)r_{\mu,\nu-\mu}=a$.
Since, according to the above made assumption, all the maps $\rho_{\mu}$ are
automorphisms, such a $b$ always exists.

It turns out \cite[Sect.~4.11]{hk:solvpoly} that under the made assumption the
results of Kandry-Rodi and Weispfenning \cite[Sects.~4/5]{krw:ncgb} remain
valid for our larger class of non-commutative algebras and can be
straightforwardly extended to involutive bases.  For this reason, we will
only discuss the case of involutive bases and do not treat separately
Gr\"obner bases.

\begin{lemma}\label{lem:lrred}
  Let\/ $(\P,\star,\prec)$ be an arbitrary polynomial algebra of solvable type
  where all the maps\/ $\rho_\mu$ appearing in the commutation relations
  (\ref{eq:crmono1}) are automorphisms.  A polynomial\/ $f\in\P$ is
  (involutively) left reducible modulo a finite set\/ $\F\subset\P$ (with
  respect to an involutive division\/ $L$), if and only if it is
  (involutively) right reducible (with respect to\/ $L$).
\end{lemma}

\begin{proof}
  Because of the made assumptions on the maps $\rho_\mu$, reducibility depends
  solely on the leading exponents.\qed
\end{proof}

\begin{proposition}\label{prop:lrunique}
  Let\/ $\H_l$ be a monic, involutively left autoreduced, minimal left
  involutive set and\/ $\H_r$ a monic, involutively right autoreduced, minimal
  right involutive set for an involutive division\/ $L$.  If\/
  $\lspanl{\H_l}=\lspanr{\H_r}=\I$, then\/ $\H_l=\H_r$.
\end{proposition}

\begin{proof}
  By definition of a minimal basis, the sets $\le{\prec}{\H_l}$ and
  $\le{\prec}{\H_r}$ are both minimal involutive bases of the monoid ideal
  $\le{\prec}{\I}$ and thus are identical.  Assume that
  $(\H_l\setminus\H_r)\cup(\H_r\setminus\H_l)\neq\emptyset$ and let $f$ be an
  element of this set with minimal leading exponent with respect to $\prec$.
  Without loss of generality, we assume that $f\in\H_l\setminus\H_r$.  Because
  of the condition $\lspanl{\H_l}=\lspanr{\H_r}$, we have
  $f\in\ispanr{\H_r}{L,\prec}$.  Thus the (by Proposition
  \ref{prop:invnormalform} unique) right involutive normal form of $f$ with
  respect to $\H_r$ is $0$.  This implies in particular that $f$ is right
  involutively reducible with respect to some $h\in\H_r$ with
  $\le{\prec}{h}\preceq\le{\prec}{f}$.
  
  If $\le{\prec}{h}\prec\le{\prec}{f}$, then $h\in\H_l$, too, as $f$ was
  chosen as a minimal element of the symmetric difference of $\H_l$ and
  $\H_r$.  Hence, by Lemma~\ref{lem:lrred}, $f$ is also left involutively
  reducible with respect to $h$ (because of
  $\le{\prec}{\H_l}=\le{\prec}{\H_r}$ the multiplicative variables of $h$ are
  the same in both cases).  But this contradicts the assumption that $\H_l$ is
  involutively left autoreduced.
  
  If $\le{\prec}{h}=\le{\prec}{f}=\mu$, then we consider the difference
  $g=f-h\in\I$: both the left involutive normal form of $g$ with respect to
  $\H_{l}$ and the right involutive normal form with respect to $\H_{r}$ must
  vanish.  By construction, $\le{\prec}{g}\prec\mu$ and
  $\supp{g}\subseteq(\supp{f}\cup\supp{h})\setminus\{\mu\}$.  Since both
  $\H_{l}$ and $\H_{r}$ are assumed to be involutively autoreduced, no term in
  this set is involutively reducible by $\le{\prec}{\H_l}=\le{\prec}{\H_r}$
  and hence we must have $\supp{g}=\emptyset$, i.\,e.\ $g=0$, a
  contradiction.\qed
\end{proof}

A direct derivation of a theory of two-sided involutive bases along the lines
of Section \ref{sec:invbas} fails, as two-sided linear combinations are rather
unwieldy objects.  A general polynomial $f\in\lspant{\H}$ for some finite set
$\H\subset\P$ is of the form
\begin{equation}\label{eq:twolincomb}
  f=\sum_{h\in\H}\sum_{i=1}^{n_h}\ell_i\star h\star r_i
\end{equation}
with polynomials $\ell_i,r_i\in\P$, i.\,e.\ we must allow several summands
with the same generator $h$.  The definition of a unique involutive standard
representation would require control over the numbers $n_h$ which seems rather
difficult.  Therefore we will take another approach and construct left
involutive bases even for two-sided ideals.  The following proposition is an
involutive version of Theorem 5.4 in \cite{krw:ncgb}.

\begin{proposition}\label{prop:two}
  Let\/ $\H\subset(\P,\star,\prec)$ be a finite set and\/ $L$ an involutive
  division. Then the following five statements are equivalent.
  \begin{description}
  \item[\phantom{ii}{\upshape (i)}] $\H$ is a left involutive basis and
    $\lspanl{\H}=\lspant{\H}$.
  \item[\phantom{i}{\upshape (ii)}] $\H$ is a right involutive basis and
    $\lspanr{\H}=\lspant{\H}$.
  \item[{\upshape (iii)}] $\H$ is a left involutive basis of $\lspanl{\H}$ and
    both $h\star x_i\in\lspanl{\H}$ and $h\star c\in\lspanl{\H}$ for all
    generators $h\in\H$, all variables $x_i$ and all coefficients $c\in\kk$.
  \item[\,{\upshape (iv)}] $\H$ is a right involutive basis of $\lspanr{\H}$
    and both $x_i\star h\in\lspanr{\H}$ and $c\star h\in\lspanr{\H}$ for all
    generators $h\in\H$, all variables $x_i$ and all coefficients $c\in\kk$.
  \item[\phantom{i}\,{\upshape (v)}] A unique generator $h\in\H$ exists for
    every polynomial $f\in\lspant{\H}$ such that
    $\le{\prec}{h}\idiv{L,\le{\prec}{\H}}\le{\prec}{f}$.
  \end{description}
\end{proposition}

\begin{proof}
  We begin with the equivalence of the first two statements.  (i) implies that
  $\ispanl{\H}{L,\prec}=\lspanl{\H}=\lspant{\H}$ and hence trivially
  $\lspanr{\H}\subseteq\lspanl{\H}$.  The same argument as in the proof of
  Proposition \ref{prop:lrunique} shows that we have in fact an equality and
  thus $\ispanr{\H}{L,\prec}=\lspanr{\H}=\lspant{\H}$, i.\,e.\ (ii).  The
  converse goes analogously.
  
  Next we consider the equivalence of (i) and (iii); the equivalence of (ii)
  and (iv) follows by the same argument.  (iii) is a trivial consequence of
  (i).  For the converse, we note that (iii) implies that $f\star
  (ct)\in\lspanl{\H}$ for all $f\in\lspanl{\H}$, all terms $t\in\TT$ and all
  constants $c\in\kk$.  Indeed, as in the proof of Proposition
  \ref{prop:fixmult} we may rewrite the monomial $ct$ as a polynomial in the
  ``terms'' $x^{i_1}\star x^{i_2}\star\cdots\star x^{i_q}$ with $i_1\leq
  i_2\leq\cdots\leq i_q$ and then apply repeatedly our assumptions.
  Obviously, this entails (i).

  The equivalence of (i) or (ii), respectively, with (v) is a trivial
  consequence of the definition of an involutive basis.\qed
\end{proof}

We would like to exploit Statement (iii) for the construction of a left
involutive basis for the two-sided ideal $\lspant{\F}$.  However, if the field
$\kk$ is infinite, then it contains an infinite number of conditions.  In the
sequel we will follow \cite[Sect.~4.11]{hk:solvpoly} and make one further
assumption about the polynomial algebra $\P$.  Let $\kk_0=\{c\in\kk\mid\forall
f\in\P:c\star f=f\star c\}$ be the constant part of the centre of $\P$.

\begin{lemma}\label{lem:cenfield}
  $\kk_0$ is a subfield of\/ $\kk$.
\end{lemma}

\begin{proof}
  It is obvious that $\kk_0$ is a subring.  Thus there only remains to show
  that with $c\in\kk_0^\times$ we have $c^{-1}\in\kk_0$, too.  If $c\in\kk_0$,
  then $x_i\star c=cx_i$, i.\,e.\ $\rho_i(c)=c$ and $h_i(c)=0$, for all $1\leq
  i\leq n$.  Now on one hand $x_i\star(c^{-1}\star c)=x_i$ and on the other
  hand
  \begin{equation}
    (x_i\star c^{-1})\star c=\rho_i(c^{-1})\rho_i(c)x_i+ch_i(c^{-1})
  \end{equation}
  ($h_i(c^{-1})\star c=ch_i(c^{-1})$ since $c\in\kk_0$).  The associativity of
  $\star$ implies now that $\rho_i(c^{-1})=c^{-1}$ and $h_i(c^{-1})=0$.  Hence
  $c^{-1}$ commutes with all variables $x_i$ and it is easy to see that this
  entails $c^{-1}\in\kk_0$.\qed
\end{proof}

We make now the assumption that either $\kk^\times=\{c_1,\dots,c_\ell\}$ is
finite or that the extension $\kk/\kk_0$ is finite, i.\,e.\ that $\kk$ is a
finite-dimensional vector space over $\kk_0$ with basis
$\{c_1,\dots,c_\ell\}$.  In the latter case, it is easy to see that it
suffices in (iii) to require that only all products $h\star c_j$ lie in
$\lspanl{\H}$, as for $c=\sum_{j=1}^\ell\lambda_jc_j$ with $\lambda_j\in\kk_0$
we have $h\star c=\sum_{j=1}^\ell\lambda_j(h\star c_j)$.

These considerations lead to the simple Algorithm \ref{alg:two} below.  It
first constructs in Line /1/ a left involutive basis $\H$ of the left ideal
$\lspanl{\F}$ (using Algorithm \ref{alg:polycompl}).  The \texttt{while} loop
in Lines /2--19/ extends the set $\H$ to a left generating set of the
two-sided ideal $\lspant{\F}$ according to our simplified version of statement
(iii) in Proposition~\ref{prop:two}.  Finally, we complete in Line /20/ this
set to an involutive basis.  Note that in Line /1/ it is not really necessary
to compute a left involutive basis; any left Gr\"obner basis would suffice as
well.  Similarly, an ordinary left normal form could be used in Lines /6/ and
/12/, respectively; the use of $\mathtt{InvLeftNormalForm}_{L,\prec}$
anticipates the final involutive basis computation in Line /20/.

\begin{algorithm}
  \caption{Left Involutive basis for two-sided ideal in $(\P,\star,\prec)$
           \label{alg:two}}
  \begin{algorithmic}[1]
    \REQUIRE finite set $\F\subset\P$, involutive division $L$
    \ENSURE left involutive basis $\H$ of $\lspant{\F}$
    \STATE $\H\leftarrow\mathtt{LeftInvBasis}_{L,\prec}(\F)$;\quad
        $\S\leftarrow\H$
    \WHILE{$\S\neq\emptyset$}
        \STATE $\T\leftarrow\emptyset$
        \FORALL{$f\in\S$}
            \FOR{$i$ \algorithmicfont{from} $1$ \algorithmicfont{to} $n$}
                \STATE $h\leftarrow
                        \mathtt{InvLeftNormalForm}_{L,\prec}(f\star x_i,\H)$
                \IF{$h\neq0$}
                    \STATE $\H\leftarrow\H\cup\{h\}$;\quad
                        $\T\leftarrow\T\cup\{h\}$
                \ENDIF
            \ENDFOR
            \FOR{$j$ \algorithmicfont{from} $1$ \algorithmicfont{to} $\ell$}
                \STATE $h\leftarrow
                        \mathtt{InvLeftNormalForm}_{L,\prec}(f\star c_j,\H)$
                \IF{$h\neq0$}
                    \STATE $\H\leftarrow\H\cup\{h\}$;\quad
                        $\T\leftarrow\T\cup\{h\}$
                \ENDIF
            \ENDFOR
        \ENDFOR
        \STATE $\S\leftarrow\T$
    \ENDWHILE
    \STATE \algorithmicreturn{$\mathtt{LeftInvBasis}_{L,\prec}(\H)$}
  \end{algorithmic}
\end{algorithm}

The termination of the \texttt{while} loop follows from the fact that under
the made assumptions $\P$ is Noetherian and hence a finite generating set of
$\lspant{\F}$ exists.  In principle, we perform here a simple breadth-first
search for it.  The termination of the involutive bases computations in Lines
/1/ and /20/, respectively, depends on the conditions discussed in the last
section.  Thus the termination is guaranteed, if the division $L$ is
constructive and Noetherian.

\section{Involutive Bases for Semigroup Orders}\label{sec:semi}

For a number of applications it is of interest to compute involutive or
Gr\"obner bases with respect to more general orders, namely \emph{semigroup
  orders} (see Appendix~\ref{sec:termord}).  This generalisation does not
affect the basic properties of polynomial algebras of solvable type as
discussed in Sect.~\ref{sec:solvalg}, but if $1$ is no longer the smallest
term, then normal form computations do no longer terminate for all inputs.  So
we can no longer apply Algorithm~\ref{alg:polycompl} directly for the
determination of involutive bases.

\begin{example}
  The \emph{Weyl algebra} $\WW_n$ is the polynomial algebra in the $2n$
  variables $x_1,\dots,x_n$ and $\partial_1,\dots,\partial_n$ with the
  following non-commutative product $\star$: for all $1\leq i\leq n$ we have
  $\partial_i\star x_i=x_i\partial_i+1$ and $\star$ is the normal commutative
  product in all other cases.  It is easy to see that $\WW_n$ is a polynomial
  algebra of solvable type for any monoid order.  For semigroup orders
  compatibility requires that $1\prec x_i\partial_i$ for all $i$.  In
  \cite{sst:grob} such orders are called \emph{multiplicative monomial
    orders}.
  
  An important class of semigroup orders is defined via real weight vectors.
  Let $(\xi,\zeta)\in\RR^n\times\RR^n$ be such that $\xi+\zeta\in\RR^n$ is
  non-negative and let $\prec$ be an arbitrary monoid order.  Then we define
  $x^\mu\partial^\nu\prec_{(\xi,\zeta)}x^\sigma\partial^\tau$, if either
  $\mu\cdot\xi+\nu\cdot\zeta<\sigma\cdot\xi+\tau\cdot\zeta$ or
  $\mu\cdot\xi+\nu\cdot\zeta=\sigma\cdot\xi+\tau\cdot\zeta$ and
  $x^\mu\partial^\nu\prec x^\sigma\partial^\tau$.  This yields a monoid order,
  if and only if both $\xi$ and $\zeta$ are non-negative.  A special case are
  the orders with weight vectors $(\xi,-\xi)$ arising from the action of the
  algebraic torus $(\kk^*)^n$ on the Weyl algebra.  They have numerous
  applications in the theory of $\D$-modules \cite{sst:grob}.\bull
\end{example}

As normal form computations do not necessarily terminate for semigroup orders,
we must slightly modify our definitions of (weak) involutive or Gr\"obner
bases.  The proof of Theorem~\ref{thm:invnormalrep} (and consequently also the
one of Corollary~\ref{cor:charinvbas} showing that a weak involutive basis of
an ideal $\I$ is indeed a basis of $\I$) requires normal form computations and
thus this theorem is no longer valid.  The same problem occurs for Gr\"obner
bases.  Therefore we must explicitly include this condition in our definition.

\begin{definition}\label{def:invbassemi}
  Let\/ $(\P,\star,\prec)$ be a polynomial algebra of solvable type where\/
  $\prec$ is an arbitrary semigroup order.  Let furthermore\/ $\I\subseteq\P$
  be a left ideal.  A \emph{Gr\"obner basis} of\/ $\I$ is a finite set\/ $\G$
  such that\/ $\lspan{\G}=\I$ and\/ $\lspan{\le{\prec}{\G}}=\le{\prec}{\I}$.
  The set\/ $\G$ is a \emph{weak involutive basis} of\/ $\I$ for the
  involutive division\/~$L$, if in addition the set\/ $\le{\prec}{\G}$ is
  weakly involutive for\/ $L$.  It is a \emph{(strong) involutive basis}, if
  it is furthermore involutively head autoreduced.
\end{definition}

In the case of Gr\"obner bases, a classical trick due to Lazard \cite{dl:gb}
consists of homogenising the input and lifting the semigroup order to a monoid
order on the homogenised terms.  One can show that computing first a Gr\"obner
basis for the ideal spanned by the homogenised input and then dehomogenising
yields a Gr\"obner basis with respect to the semigroup order.  Note, however,
that in general we cannot expect that \emph{reduced} Gr\"obner bases exist.

We extend now this approach to involutive bases.  Here we encounter the
additional difficulty that we must lift not only the order but also the used
involutive division.  In particular, we must show that properties like
Noetherity or continuity are preserved by the lift which is non-trivial.  For
the special case of involutive bases in the Weyl algebra, this problem was
first solved in \cite{wms:weyl}.

Let $(\P,\star,\prec)$ be a polynomial algebra of solvable type where $\prec$
is any semigroup order that respects the multiplication $\star$.  We set
$\tilde{\P}=\kk[x_0,x_1,\dots,x_n]$ and extend the multiplication $\star$ to
$\tilde{\P}$ by defining that $x_0$ commutes with all other variables and the
elements of the field $\kk$.  For a polynomial $f=\sum c_\mu x^\mu\in\P$ of
degree $q$, we introduce as usual its \emph{homogenisation} $f^{(h)}=\sum
c_\mu x_0^{q-|\mu|}x^\mu\in\tilde{\P}$.  Conversely, for a polynomial $\tilde
f\in\tilde{\P}$ we denote its projection to $\P$ as $f=\tilde f|_{x_0=1}$.

We denote by $\tilde\TT$ the set of terms in $\tilde{\P}$; obviously, it is as
monoid isomorphic to $\Nno$.  We use in the sequel the following convention.
Multi indices in $\Nno$ always carry a tilde: $\tilde\mu=[\mu_0,\dots,\mu_n]$.
The projection to $\Nn$ defined by dropping the first entry (i.\,e.\ the
exponent of the homogenisation variable $x_0$) is signalled by omitting the
tilde; thus $\mu=[\mu_1,\dots,\mu_n]$.  For subsets $\tilde\N\subset\Nno$ we
also simply write $\N=\{\nu\mid\tilde\nu\in\tilde\N\}\subset\Nn$.

We lift the semigroup order $\prec$ on $\TT$ to a monoid order $\prec_h$ on
$\tilde\TT$ by defining $x^{\tilde\mu}\prec_h x^{\tilde\nu}$, if either
$|\tilde\mu|<|\tilde\nu|$ or both $|\tilde\mu|=|\tilde\nu|$ and $x^\mu\prec
x^\nu$.  It is trivial to check that this yields indeed a monoid order and
that $(\tilde{\P},\star,\prec_h)$ is again a polynomial algebra of solvable
type.  For lifting the involutive division, we proceed somewhat similarly to
the definition of the Janet division: the homogenisation variable $x_0$ is
multiplicative only for terms which have maximal degree in $x_0$.

\begin{proposition}\label{prop:divlift}
  Let\/ $L$ be an arbitrary involutive division on\/ $\Nn$.  For any finite
  set\/ $\tilde\N\subset\Nno$ and every multi index\/ $\tilde\mu\in\tilde\N$,
  we define\/ $\mult{N}{\tilde L,\tilde\N}{\tilde\mu}$ by:
  \begin{itemize}
  \item $0\in\mult{N}{\tilde L,\tilde\N}{\tilde\mu}$, if and only if\/
    $\mu_0=\max_{\tilde\nu\in\tilde\N}\{\nu_0\}$,
  \item $0<i\in\mult{N}{\tilde L,\tilde\N}{\tilde\mu}$, if and only if\/
    $i\in\mult{N}{L,\N}{\mu}$.
  \end{itemize}
  This determines an involutive division\/ $\tilde L$ on\/ $\Nno$.
\end{proposition}

\begin{proof}
  Both conditions for an involutive division are easily verified.  For the
  first one, let $\tilde\rho\in\cone{\tilde L,\tilde\N}{\tilde\mu}\cap
  \cone{\tilde L,\tilde\N}{\tilde\nu}$ with $\tilde\mu,\tilde\nu\in\tilde\N$.
  If $\rho_0=\mu_0=\nu_0$, the first entry can be ignored, and the properties
  of the involutive division $L$ implies the desired result.  If
  $\rho_0=\mu_0>\nu_0$, the index $0$ must be multiplicative for $\tilde\nu$
  contradicting $\mu_0>\nu_0$.  If $\rho_0$ is greater than both $\mu_0$ and
  $\nu_0$, the index $0$ must be multiplicative for both implying
  $\mu_0=\nu_0$.  In this case we may again ignore the first entry and invoke
  the properties of $L$.
  
  For the second condition we note that whether a multiplicative index $i>0$
  becomes non-multiplicative for some element $\tilde\nu\in\tilde\N$ after
  adding a new multi index to $\tilde\N$ is independent of the first entry and
  thus only determined by the involutive division $L$.  If the new multi index
  has a higher first entry than all elements of $\tilde\N$, then $0$ becomes
  non-multiplicative for all elements in $\tilde\N$ but this is permitted.\qed
\end{proof}

Now we check to what extent the properties of $L$ are inherited by the lifted
division $\tilde L$.  Given the similarity of the definition of $\tilde L$ and
the Janet division, it is not surprising that we may reuse many ideas from
proofs for the latter.

\begin{proposition}\label{prop:noelift}
  If\/ $L$ is a Noetherian division, then so is\/ $\tilde L$.
\end{proposition}

\begin{proof}
  Let $\tilde{\N}\subset\Nno$ be an arbitrary finite subset.  In order to
  prove the existence of an $\tilde L$-completion of $\tilde\N$, we first take
  a finite $L$-completion $\hat\N\subset\Nn$ of $\N$ which always exists, as
  by assumption $L$ is Noetherian.  Next, we define a finite subset
  $\tilde\N'\subset\lspan{\tilde\N}$ by setting
   \begin{displaymath}
     \tilde\N' =
         \Bigl\{\,\tilde\mu\in\NN_0^{n+1}\mid
                \mu\in\hat\N\wedge
                \mu_0\leq\max_{\tilde\nu\in\tilde\N}\nu_0\,\Bigr\} 
         \cap\lspan{\tilde\N}\;.
   \end{displaymath}
   We claim that this set $\tilde\N'$ is an $\tilde L$-completion of
   $\tilde\N$.  By construction, we have both
   $\tilde\N'\subset\lspan{\tilde\N}$ and $\tilde\N\subseteq\tilde\N'$, so
   that we must only show that $\tilde\N'$ is involutive.
   
   Let $\tilde\mu\in\lspan{\tilde\N'}$ be arbitrary.  By construction of
   $\tilde\N'$, we can find $\tilde\nu\in\tilde\N'$ with
   $\nu\idiv{L,{\hat\N}}\mu$.  Moreover, the definition of $\tilde\N'$
   guarantees that we can choose $\tilde\nu$ in such a way that either
   $\nu_0=\mu_0$ or $\nu_0=\max_{\tilde\rho\in\tilde\N'}\rho_0<\mu_0$ holds.
   In the former case, we trivially have $\tilde\nu\idiv{\tilde
     L,\tilde\N'}\tilde\mu$; in the latter case we have $0\in\mult{N}{\tilde
     L,\tilde\N}{\tilde\nu}$ (see the proof of
   Proposition~\ref{prop:divlift}).  Thus in either case
   $\tilde\mu\in\lspan{\tilde\N'}_{\tilde L}$.\qed
\end{proof}

\begin{proposition}\label{prop:contlift}
  If\/ $L$ is a continuous division, then so is\/ $\tilde L$.
\end{proposition}

\begin{proof}
  Let $(\tilde\nu^{(1)},\dots,\tilde\nu^{(r)})$ with
  $\tilde\nu^{(i)}\in\tilde\N$ be a finite sequence as described in the
  definition of continuity.  We first note that the integer sequence
  $(\nu^{(1)}_0,\dots,\nu^{(r)}_0)$ is monotonically increasing.  If
  $\nu^{(i)}_0$ is not maximal among the entries $\mu_0$ for
  $\tilde\mu\in\tilde\N$, no multiplicative divisor of $\tilde\nu^{(i)}+1_j$
  in $\tilde\N$ can have a smaller first entry: if $\nu^{(i)}_0$ is maximal,
  the index $0$ is multiplicative for $\tilde\nu^{(i)}$ and any involutive
  divisor in $\tilde\N$ must also be maximal in the zero entry.  Thus it
  suffices to look at those parts of the sequence where equality in the zero
  entries holds.  But there the inequality of the multi indices
  $\tilde\nu^{(i)}$ follows from the continuity of the underlying division
  $L$.\qed
\end{proof}

Unfortunately, it is much harder to show that constructivity is preserved.  We
will do this only for globally defined divisions and the Janet division.

\begin{proposition}
  If the continuous division\/ $L$ is either globally defined or the Janet
  division, then the lifted division\/ $\tilde L$ is constructive.
\end{proposition}

\begin{proof}
  We give a proof only for the case of a globally defined division.  For the
  Janet division $J$ one must only make a few modifications of the proof that
  $J$ itself is constructive.  We omit the details; they can be found in
  \cite{wms:weyl}.
  
  We select a finite set $\tilde\N\subset\Nno$, a multi index
  $\tilde\mu\in\tilde\N$ and a non-multiplicative index $i$ of $\tilde\mu$
  such that the conditions in the definition of constructivity are fulfilled.
  Assume that there exists a $\tilde\rho\in\tilde\N$ such that
  $\tilde\mu+1_i=\tilde\rho+\tilde\sigma+\tilde\tau$ with
  $\tilde\rho+\tilde\sigma\in\cone{\tilde L,\tilde\N}{\tilde\rho}$ and
  $\tilde\rho+\tilde\sigma+\tilde\tau\in\cone{\tilde L,\tilde
    \N\cup\{\tilde\rho+\tilde\sigma\}}{\tilde\rho+\tilde\sigma}$.  Let $L$ be
  a globally defined division.  If $i=0$, then
  $\mu_0+1=\rho_0+\sigma_0+\tau_0$ implies that $\sigma_0=\tau_0=0$: for
  $\sigma_0>0$, we would have ($0$ is multiplicative for $\tilde\rho$)
  $\rho_0>\mu_0\geq\rho_0+\sigma_0>\rho_0$.  For $\sigma_0=0$ and $\tau_0>0$ a
  similar contradiction appears.  If $i>0$, the argumentation is simple.  A
  global division is always constructive, as adding further elements to $\N$
  does not change the multiplicative indices.  But the same holds for the
  indices $k>0$ in the lifted division $\tilde L$.  Thus under the above
  conditions $\tilde\mu+1_i\in\ispan{\tilde\N}{\tilde L}$ contradicting the
  made assumptions.\qed
\end{proof}

Based on these results, Algorithm \ref{alg:polycompl} can be extended to
semigroup orders.  Given a finite set $\F\in\P$, we first determine its
homogenisation $\F^{(h)}\in\tilde{\P}$ and then compute an involutive basis of
$\lspan{F^{(h)}}$ with respect to $\tilde L$ and $\prec_h$.  What remains to
be done is first to show that the existence of a finite involutive basis is
preserved under the lifting to $\tilde{\P}$ and then to study the properties
of the dehomogenisation of this basis.

\begin{proposition}\label{prop:finite}
  If the left ideal\/ $\I=\lspan{\F}\subseteq\P$ possesses an involutive basis
  with respect to the Noetherian division\/~$L$ and the semigroup
  order~$\prec$, then the left ideal\/
  $\tilde\I=\lspan{\F^{(h)}}\subseteq\tilde{\P}$ generated by the
  homogenisations of the elements in the finite set\/ $\F$ possesses an
  involutive basis with respect to the lifted division\/~$\tilde L$ and the
  monoid order\/~$\prec_h$.
\end{proposition}

\begin{proof}
  By Theorem \ref{thm:noetherfield}, the ideal $\tilde\I\subseteq\tilde{\P}$
  possesses a Gr\"obner basis $\tilde\G$ with respect to the monoid order
  $\prec_h$.  By Proposition \ref{prop:noelift}, a finite $\tilde
  L$-completion $\tilde\N$ of the set $\le{\prec_h}{\tilde\G}$ exists.
  Moreover, as $\tilde\G$ is a Gr\"obner basis of $\tilde\I$, the monoid
  ideals $\lspan{\le{\prec_h}{\tilde\G}}$ and $\le{\prec_h}{\tilde\I}$
  coincide.  Thus $\tilde\N$ is an involutive basis of
  $\le{\prec_h}{\tilde\I}$ with respect to the lifted division $\tilde L$ and
  an involutive basis $\tilde\H$ of $\tilde\I$ with respect to the division
  $\tilde L$ is given by
  \begin{equation}
    \tilde\H = \bigl\{\,x^{\tilde\mu}\star\tilde g\mid
                        {\tilde g}\in\tilde\G\ \wedge\ 
     \le{\prec_h}{(x^{\tilde\mu}\star\tilde g)}\in\tilde\N\,\bigr\}\;.
  \end{equation}
  This set is obviously finite.\qed
\end{proof}

Hence our lifting leads to a situation where we can apply
Theorem~\ref{thm:polycompl}.  Unfortunately, the dehomogenisation of the
strong involutive basis computed in $\tilde{\P}$ does not necessarily lead to
a \emph{strong} involutive basis in $\P$, but we obtain always at least a weak
involutive basis and thus in particular a Gr\"obner basis.

\begin{theorem}\label{thm:hominv}
  Let\/ $\tilde\H$ be a strong involutive basis of the left ideal\/
  $\tilde\I\subseteq\tilde{\P}$ with respect to\/ $\tilde L$ and $\prec_h$.
  Then the dehomogenisation\/ $\H$ is a weak involutive basis of the left
  ideal\/ $\I\subseteq\P$ with respect to\/ $L$ and $\prec$.
\end{theorem}

\begin{proof}
  For any $f\in\I$ an integer $k\geq0$ exists such that $\tilde
  f=x_0^kf^{(h)}\in\tilde\I$.  The polynomial $\tilde f$ possesses a unique
  involutive standard representation
  \begin{equation}\label{eq:isrf}
    \tilde f=\sum_{\tilde h\in\tilde\H}\tilde P_{\tilde h}\tilde h
  \end{equation}
  with $\tilde P_{\tilde h}\in\kk[\mult{X}{\tilde
    L,\le{\prec_h}{\tilde\H}}{\tilde h}]$ and $\le{\prec_h}{(\tilde P_{\tilde
      h}\tilde h)}\preceq_h\le{\prec_h}{\tilde f}$.  Setting $x_0=1$ in
  (\ref{eq:isrf}) yields a representation of $f$ with respect to the
  dehomogenised basis\footnote{Note that the dehomogenised basis $\H$ is in
    general smaller than $\tilde\H$, as some elements of $\tilde\H$ may differ
    only in powers of $x_0$.} $\H$ of the form $f=\sum_{h\in\H}P_hh$ where
  $P_h\in\kk[\mult{X}{L,\le{\prec}{\H}}{h}]$ by the definition of the lifted
  division $\tilde L$.  This obviously implies that $\lspan{\H}=\I$.  By the
  definition of the lifted order $\prec_h$ and the homogeneity of the lifted
  polynomials, we have furthermore that
  $\le{\prec}{(P_hh)}\preceq\le{\prec}{f}$ and hence that $\le{\prec}{\H}$ is
  a weak involutive basis of $\le{\prec}{\I}$.  Since all conditions of
  Definition~\ref{def:invbassemi} are satisfied, the set $\H$ is therefore
  indeed a weak involutive basis of the ideal $\I$.\qed
\end{proof}

\begin{remark}\label{rem:pomlaz}
  For the Pommaret division $P$ the situation is considerably simpler.  There
  is no need to define a lifted division $\tilde P$ according to Proposition
  \ref{prop:divlift}.  Instead we renumber $x_0$ to $x_{n+1}$ and then use the
  standard Pommaret division on $\Nno$.  This approach implies that for all
  multi indices $\tilde\mu\in\Nno$ with $\mu\neq0$ the equality
  $\mult{N}{P}{\tilde\mu}=\mult{N}{P}{\mu}$ holds, as obviously $n+1$ is
  multiplicative only for multi indices of the form $\tilde\mu=\ell_{n+1}$,
  i.\,e.\ for which $\mu=0$.  One easily sees that the above proof of Theorem
  \ref{thm:hominv} is not affected by this change of the division used in
  $\Nno$ and hence remains true.\bull
\end{remark}

It is not a shortcoming of our proof that in general we do not get a strong
involutive basis, but actually some ideals do not possess strong involutive
bases.  In particular, there is no point in invoking
Proposition~\ref{prop:weakbasis} for obtaining a strong basis.  While we may
surely obtain by elimination a subset $\H'\subseteq\H$ such that
$\le{\prec}{\H'}$ is a strong involutive basis of $\lspan{\le{\prec}{\H}}$, in
general $\lspan{\H'}\subsetneq\I$.

\begin{example}\label{ex:nostrongbas}
  Consider in the Weyl algebra ${\mathbbm W}_2=\kk[x,y,\partial_x,\partial_y]$
  the left ideal generated by the set
  $\F=\{\underline{1}+x+y,\underline{\partial_y}-\partial_x\}$.  We take the
  semigroup order induced by the weight vector $(-1,-1,1,1)$ and refined by a
  term order for which $\partial_y\succ\partial_x\succ y\succ x$.  Then the
  underlined terms are the leading ones.  One easily checks that $\F$ is a
  Gr\"obner basis for this order.  Furthermore, all variables are
  multiplicative for each generator with respect to the Pommaret division and
  thus $\F$ is a weak Pommaret basis, too.
  
  Obviously, the set $\F$ is neither a reduced Gr\"obner basis nor a strong
  Pommaret basis, as $1$ is a (multiplicative) divisor of $\partial_y$.
  However, it is easy to see that the left ideal $\I=\lspan{\F}$ does not
  possess a reduced Gr\"obner basis or a strong Pommaret basis.  Indeed, we
  have $\le{\prec}{\I}=\NN_0^4$ and thus such a basis had to consist of only a
  single generator; but $\I$ is not a principal ideal.\bull
\end{example}

A special situation arises for the Janet division.  Recall from
Remark~\ref{rem:autored} that any finite set $\N\subset\Nn$ is automatically
involutively autoreduced with respect to the Janet division.  Thus any weak
Janet basis is a strong basis, if all generators have different leading
exponents.  If we follow the above outlined strategy of applying
Algorithm~\ref{alg:polycompl} to a homogenised basis and then dehomogenising
the result, we cannot generally expect this condition to be satisfied.
However, with a minor modification of the algorithm we can achieve this goal.

\begin{theorem}\label{thm:strjanbas}
  Let\/ $(\P,\star,\prec)$ be a polynomial algebra of solvable type where\/
  $\prec$ is an arbitrary semigroup order.  Then every left ideal\/
  $\I\subseteq\P$ possesses a strong Janet basis for\/ $\prec$.
\end{theorem}

\begin{proof}
  Assume that at some intermediate stage of Algorithm~\ref{alg:polycompl} the
  basis $\tilde\H$ contains two polynomials $\tilde f$ and $\tilde g$ such
  that $\le{\prec_h}{(\tilde g)}=\le{\prec_h}{(\tilde f)}+1_0$, i.\,e.\ the
  leading exponents differ only in the first entry.  If $\tilde g=x_{0}\tilde
  f$, we will find $f=g$ after the dehomogenisation and no obstruction to a
  strong basis appears.  Otherwise we note that, by definition of the lifted
  Janet division $J_h$, the homogenisation variable $x_0$ is
  non-multiplicative for $\tilde f$.  Thus at some later stage the algorithm
  must consider the non-multiplicative product $x_0\tilde f$ (if it was
  already treated, $\tilde\H$ would not be involutively head autoreduced).
  
  In the usual algorithm, we then determine the involutive normal form of the
  polynomial $x_0\tilde f$; the first step of this computation is to replace
  $x_0\tilde f$ by $x_0\tilde f-\tilde g$.  Alternatively, we may proceed
  instead as follows.  The polynomial $\tilde g$ is removed from the basis
  $\tilde\H$ and replaced by $x_0\tilde f$.  Then we continue by analysing
  the involutive normal form of $\tilde g$ with respect to the new basis.
  Note that this modification concerns only the situation that a
  multiplication by $x_0$ has been performed and that the basis $\tilde\H$
  contains already an element with the same leading exponent as the obtained
  polynomial.
  
  If the final output $\tilde\H$ of the thus modified completion algorithm
  contains two polynomials $\tilde f$ and $\tilde g$ such that
  $\le{\prec_h}{(\tilde g)}$ and $\le{\prec_h}{(\tilde f)}$ differ only in the
  first entry, then either $\tilde g=x_0^k\tilde f$ or $\tilde f=x_0^k\tilde
  g$ for some $k\in\NN$.  Thus the dehomogenisation yields a basis $\H$ where
  all elements possess different leading exponents and $\H$ is a strong Janet
  basis.  Looking at the proof of Theorem~\ref{thm:polycompl}, it is easy to
  see that this modification does not affect the correctness and the
  termination of the algorithm.  As the Janet division is Noetherian, these
  considerations prove together with Proposition~\ref{prop:noelift} the
  assertion.\qed
\end{proof}

Note that our modification only achieves its goal, if we really restrict in
Algorithm~\ref{alg:polycompl} to head reductions.  Otherwise some other terms
than the leading term in $x_0\tilde f$ might be reducible but not the
corresponding terms in $\tilde f$.  Then we could still find after the
dehomogenisation two generators with the same leading exponent.

\begin{example}\label{ex:strjanbas}
  Let us consider in the Weyl algebra $\WW_3$ with the three variables $x$,
  $y$, $z$ the left ideal generated by the set
  $\F=\{\partial_z-y\partial_x,\,\partial_y\}$.  If we apply the usual
  involutive completion Algorithm~\ref{alg:polycompl} (to the homogenisation
  $\F^{(h)}$), we obtain for the weight vector $(-1,0,0,1,0,0)$ refined by the
  degree reverse lexicographic order and the Janet division the following weak
  basis with nine generators:
  \begin{equation}
    \H_1=\bigl\{\ \partial_x,\ \partial_y,\ \partial_z,\ \partial_x\partial_z,\
    \partial_y\partial_z,\  y\partial_x,\  y\partial_x+\partial_z,\
    y\partial_x\partial_z,\ y\partial_x\partial_z+\partial_z^2\
    \bigr\}\;.
  \end{equation}
  As one easily sees from the last four generators, it is not a strong basis.
  
  Applying the modified algorithm for the Janet division yields the following
  basis with only seven generators:
  \begin{equation}
    \H_2=\bigl\{\ \partial_x+\partial_y\partial_z,\ \partial_y,\
      \partial_z,\ \partial_x\partial_z,\ \partial_y\partial_z,\
    y\partial_x+\partial_z,\ y\partial_x\partial_z+\partial_z^2\ \bigr\}\;.
  \end{equation}
  Obviously, we now have a strong basis, as all leading exponents are
  different.
  
  This example also demonstrates the profound effect of the homogenisation.  A
  strong Janet or Pommaret basis of $\lspan{\F}$ is simply given by
  $\H=\{\partial_x,\,\partial_y,\,\partial_z\}$ which is simultaneously a
  reduced Gr\"obner basis.  In $\lspan{\F^{(h)}}$ many reductions are not
  possible because the terms contain different powers of $x_0$.  However, this
  is a general problem of all approaches to Gr\"obner bases for semigroup
  orders using homogenisation and not specific for the involutive approach.
  
  In this particular case, one could have applied the involutive completion
  algorithm directly to the original set $\F$ and it would have terminated
  with the minimal basis $\H$, although we are using a order which is not a
  monoid order.  Unfortunately, it is not clear how to predict when infinite
  reduction chains appear in normal form computations with respect to such
  orders, so that one does not know in advance whether one may dispense with
  the homogenisation.\bull
\end{example}

\section{Involutive Bases for Semigroup Orders II: Mora's Normal Form}
\label{sec:semi2}

One computational disadvantage of the approach outlined in the previous
section is that the basis $\tilde\H$ in the homogenised algebra $\tilde{\P}$
is often much larger than the final basis $\H$ in the original algebra $\P$,
as upon dehomogenisation generators may become identical.  Furthermore, we
have seen that it is difficult to prove the constructivity of the lifted
division $L_h$ which limits the applicability of this technique.  Finally, for
most divisions we are not able to determine strong bases.

An alternative approach for Gr\"obner bases computations in the ordinary
polynomial ring was proposed first by Greuel and Pfister \cite{gp:standard}
and later independently by Gr\"abe \cite{hgg:tangent,hgg:local}; extensive
textbook discussions are contained in \cite[Chapt.~4]{clo:uag} and
\cite[Sect.~1.6]{gp:singular}.  It allows us to dispense completely with
computing in the homogenised algebra $\tilde{\P}$.  Two ideas are the core of
this approach: we modify the normal form algorithm using ideas developed by
Mora \cite{tm:tangent} for the computation of tangent cones and we work over a
ring of fractions of $\P$.  We will now show that a generalisation to
arbitrary polynomial algebras of solvable type and to involutive normal forms
is possible and removes all the mentioned problems.

The central problem in working with semigroup orders is that they are no
longer well-orders and hence normal form computations in the classical form do
not necessarily terminate.  Mora \cite{tm:tangent} introduced the notion of
the \emph{\'ecart} of a polynomial $f$ as the difference between the lowest
and the highest degree of a term in $f$ and based a new normal form algorithm
on it which always terminates.  The main differences between it and the usual
algorithm lie in the possibility to reduce also with respect to intermediate
results (see Line /9/ in Algorithm~\ref{alg:mora} below) and that it computes
only a ``weak'' normal form (cf.\ Proposition~\ref{prop:mora} below).

Mora's approach is valid only for tangent cone orders where the leading term
is always of minimal degree.  Greuel and Pfister \cite{gp:standard} noticed
that a slight modification of the definition of the \'ecart allows us to use
it for arbitrary semigroup orders.  So we set for any polynomial
$f\in\P\setminus\{0\}$ and any semigroup order $\prec$
\begin{equation}\label{eq:ecart}
  \ec{f}=\deg{f}-\deg{\lt{\prec}{f}}\;.
\end{equation}

The extension of the Mora normal form to an involutive normal form faces one
problem.  As already mentioned, one allows here also reductions with respect
to some intermediate results and thus one must decide on the assignment of
multiplicative variables to these.  However, it immediately follows from the
proof of the correctness of the Mora algorithm how this assignment must be
done in order to obtain in the end an involutive standard representation with
respect to the set $\G$ (one should stress that this assignment is \emph{not}
performed according to some involutive division in the sense of
Definition~\ref{def:invdiv}).

In Algorithm~\ref{alg:mora} below we use the following approach.  To each
member $g$ of the set $\hat{\G}$ with respect to which we reduce we assign a
set $N[g]$ of multiplicative indices.  We write
$\le{\prec}{g}\idiv{N}\le{\prec}{h}$, if the multi index $\le{\prec}{h}$ lies
in the restricted cone of $\le{\prec}{g}$ defined by $N[g]$.  The set $\S$
collects all generators $g\in\G$ which have already been used for reductions
and the set $\N$ is the intersection of the corresponding sets of
multiplicative indices.  If a new polynomial $h$ is added to $\hat{\G}$, it is
assigned as multiplicative indices the current value of $\N$.

\begin{algorithm}
  \caption{Involutive Mora normal form for a semigroup 
           order $\prec$ on $\P$\label{alg:mora}}
  \begin{algorithmic}[1]
    \REQUIRE polynomial $f\in\P$, finite set $\G\subset\P$, 
             involutive division $L$
    \ENSURE involutive Mora normal form $h$ of $f$ with respect to $\G$
    \STATE $h\leftarrow f$;\quad $\hat{\G}\leftarrow\G$
    \FORALL{$g\in\G$}
        \STATE $N[g]\leftarrow\mult{N}{L,\le{\prec}{\G}}{\le{\prec}{g}}$
    \ENDFOR
    \STATE $\N\leftarrow\{1,\dots,n\};\quad\S\leftarrow\emptyset$
    \WHILE{$(h\neq0)\wedge
            (\exists\,g\in\hat{\G}:\le{\prec}{g}\idiv{N}\le{\prec}{h})$}
        \STATE choose $g$ with $\ec{g}$ minimal among all $g\in\hat{\G}$ such
             that $\le{\prec}{g}\idiv{N}\le{\prec}{h}$
        \IF{$(g\in\G)\wedge(g\notin\S)$}
            \STATE $\S\leftarrow\S\cup\{g\};\quad\N\leftarrow\N\cap N[g]$
        \ENDIF
        \IF{$\ec{g}>\ec{h}$}
            \STATE $\hat{\G}\leftarrow\hat{\G}\cup\{h\};\quad N[h]\leftarrow\N$
        \ENDIF
        \STATE $\mu\leftarrow\le{\prec}{h}-\le{\prec}{g}$;
        \quad $h\leftarrow h-
            \frac{\lc{\prec}{h}}{\lc{\prec}{(x^\mu\star g)}}x^\mu\star g$
    \ENDWHILE
    \STATE \algorithmicreturn{$h$}
  \end{algorithmic}
\end{algorithm}

\begin{proposition}\label{prop:mora}
  Algorithm~\ref{alg:mora} always terminates.  Let\/ $(\P=\kk[X],\star,\prec)$
  be a polynomial algebra of solvable type (for an arbitrary semigroup order\/
  $\prec$) such that\/ $\kk[X']$ is a subring of\/ $\P$ for any subset\/
  $X'\subset X$.  Then the output\/ $h$ is a\/ \emph{weak involutive normal
    form} of the input\/ $f$ with respect to the set\/ $\G$ in the sense that
  there exists a polynomial\/ $u\in\P$ with\/ $\le{\prec}{u}=0$ such that the
  difference\/ $u\star f-h$ possesses an involutive standard representation
  \begin{equation}\label{eq:morarep}
     u\star f-h=\sum_{g\in\G}P_g\star g
  \end{equation}
  and none of the leading exponents\/ $\le{\prec}{g}$ involutively divides\/
  $\le{\prec}{h}$.  If\/ $\prec$ is a monoid order, then\/ $u=1$ and\/ $h$ is
  an involutive normal form in the usual sense.
\end{proposition}

\begin{proof}
  As the proof is almost identical with the one for the non-involutive version
  of the Mora normal form given by Greuel and Pfister
  \cite{gp:standard,gp:singular}, we only sketch the required modifications;
  full details are given in \cite{wms:invol}.  For the termination proof no
  modifications are needed.  For the existence of the involutive standard
  representation one uses the same induction as in the non-involutive case and
  keeps track of the multiplicative variables.  The key point is that if a
  reduction with respect to a polynomial $\hat{g}\in\hat{\G}\setminus\G$ is
  performed, then this polynomial is multiplied only with terms which are
  multiplicative for all $g\in\G$ appearing in $\hat{g}$.  This ensures that
  in the end indeed each non-zero coefficient $P_g$ is contained in
  $\kk[\mult{X}{L,\G,\prec}{g}]$.\qed
\end{proof}

\begin{remark}
  The assumption about $\P$ in Proposition \ref{prop:mora} is necessary,
  because the coefficients $P_g$ in (\ref{eq:morarep}) are the result of
  multiplications.  While the above considerations ensure that each factor
  lies in $\kk[\mult{X}{L,\G,\prec}{g}]$, it is unclear in a general
  polynomial algebra whether this remains true for their product.  Simple
  examples for polynomial algebras of solvable type satisfying the made
  assumption are rings of linear difference or differential operators.  In the
  case of the Pommaret division, the assumption can be weaken a bit and every
  iterated polynomial algebra of solvable type in the sense of Definition
  \ref{def:iteralg} is permitted, too.\bull
\end{remark}

We move now to a larger ring of fractions where all polynomials with leading
exponent $0$ are units.  In such a ring it really makes sense to call $h$ a
(weak) normal form of $f$, as we multiply $f$ only by a unit.

\begin{proposition}\label{prop:multset}
  Let\/ $(\P,\star,\prec)$ be a polynomial algebra of solvable type where\/
  $\prec$ is a semigroup order.  Then the subset
  \begin{equation}\label{eq:multset}
    \Sp=\{f\in\P\mid\le{\prec}{f}=0\}\;.
  \end{equation}
  is multiplicatively closed and the left localisation\/
  $\locp=\Sp^{-1}\star\P$ is a well defined ring of left fractions.
\end{proposition}

\begin{proof}
  Obviously, $1\in\Sp$.  If $1+f$ and $1+g$ are two elements in $\Sp$, then
  the compatibility of the order $\prec$ with the multiplication $\star$
  ensures that their product is of the form $(1+f)\star(1+g)=1+h$ with
  $\le{\prec}{h}\prec0$.  Hence the set $\Sp$ is multiplicatively closed.
  
  As polynomial algebras of solvable type do not possess zero divisors, a
  sufficient condition for the existence of the ring of left fractions
  $\Sp^{-1}\star\P$ is that for all $f\in\Sp$ and $g\in\P$ the intersection
  $(\P\star f)\cap(\Sp\star g)$ is not empty \cite[Sect.~12.1]{pmc:alg2}.  But
  this can be shown using minor modifications of our proof of
  Proposition~\ref{prop:ore}.
  
  We first choose coefficients $r_0,s_0\in\R$ such that in $\bar h_1=r_0g\star
  f-s_0f\star g$ the leading terms cancel, i.\,e.\ we have $\le{\prec}{\bar
    h_1}\prec\le{\prec}{f}+\le{\prec}{g}=\le{\prec}{g}$.  Then we compute with
  (the non-involutive form of) Algorithm \ref{alg:mora} a weak normal form
  $h_1$ of $\bar h_1$ with respect to the set $\F_0=\{f,g\}$.  By Proposition
  \ref{prop:mora} this yields a standard representation $u_1\star\bar
  h_1-h_1=\phi_0\star f+\psi_0\star g$ where $\le{\prec}{u_1}=0$.  Assume that
  $\le{\prec}{\psi_0}\succeq0$.  Then we arrive at the contradiction
  $\le{\prec}{(\psi_0\star g)}\succeq\le{\prec}{g}\succ\le{\prec}{\bar h_1}=
  \le{\prec}{(u_1\star\bar h_1)}$.  Thus $\le{\prec}{\psi_0}\prec0$.  If
  $h_1=0$, then $(u_1\star r_0g-\phi_0)\star f=(u_1\star s_0f+\psi_0)\star g$
  and by the considerations above on the leading exponents $u_1\star
  s_0f+\psi_0\in\Sp$ so that indeed $(\P\star f)\cap(\Sp\star
  g)\neq\emptyset$.
  
  If $h_1\neq0$, we proceed as in the proof of Proposition~\ref{prop:ore}.  We
  introduce $\F_1=\F_0\cup\{h_1\}$ and choose $r_1,s_1\in\R$ such that in
  $\bar h_2=r_1h_1\star f-s_1f\star h_1$ the leading terms cancel.  If we
  compute a weak Mora normal form $h_2$ of $\bar h_2$, then we obtain a
  standard representation $u_2\star\bar h_2-h_2=\phi_1\star f+\psi_1\star
  g+\rho_1\star h_1$ where again $\le{\prec}{u_2}=0$.  The properties of a
  standard representation imply now that
  $\le{\prec}{\psi_1}+\le{\prec}{g}\preceq\le{\prec}{\bar h_2}$ and
  $\le{\prec}{\rho_1}+\le{\prec}{h_1}\preceq\le{\prec}{\bar h_2}$.  Together
  with the inequalities $\le{\prec}{\bar h_2}\prec
  \le{\prec}{f}+\le{\prec}{h_1}= \le{\prec}{h_1}\prec\le{\prec}{g}$ this
  entails that both $\le{\prec}{\psi_1}\prec0$ and $\le{\prec}{\rho_1}\prec0$.
  Thus for $h_2=0$ we have found $\phi\in\P$ and $\psi\in\Sp$ such that
  $\phi\star f=\psi\star g$.  If $h_2\neq0$, similar inequalities in the
  subsequent iterations ensure that we always have $\psi\in\Sp$.\qed
\end{proof}

As any localisation of a Noetherian ring is again Noetherian, $\locp$ is
Noetherian, if $\P$ is so.  One sees immediately that the units in $\locp$ are
all those fractions where not only the denominator but also the numerator is
contained in $\Sp$.  Given an ideal $\I\subseteq\locp$, we may always assume
without loss of generality that its generated by a set $\F\subset\P$ of
polynomials, as multiplication of a generator by a unit does not change the
span.  Hence in all computations we will exclusively work with polynomials and
not with fractions.

As all elements of $\Sp$ are units in $\locp$, we may extend the notions of
leading term, monomial or exponent: if $f\in\locp$, then we can choose a unit
$u\in\Sp$ with $\lc{\prec}{u}=1$ such that $u\star f\in\P$ is a polynomial;
now we define $\le{\prec}{f}=\le{\prec}{(u\star f)}$ etc.  One easily verifies
that this definition is independent of the choice of $u$.

Following Greuel and Pfister \cite{gp:singular}, one can now construct a
complete theory of involutive bases over $\locp$.
Definition~\ref{def:invbassemi} of Gr\"obner and involutive bases can be
extended without changes from the ring $\P$ to $\locp$.
Theorem~\ref{thm:noetherfield} on the existence of Gr\"obner bases generalises
to $\locp$, as its proof is only based on the leading exponents and a simple
normal form argument remaining valid due to our considerations above.

Note that even if the set $\G$ is involutively head autoreduced, we cannot
conclude in analogy to Proposition~\ref{prop:invnormalform} that the
involutive Mora normal form is unique, as we only consider the leading term in
Algorithm~\ref{alg:mora} and hence the lower terms in $h$ may still be
involutively divisible by the leading term of some generator $g\in\G$.
However, Theorem~\ref{thm:invnormalrep} remains valid.

\begin{theorem}\label{thm:mora}
  Let\/ $(\P=\kk[X],\star,\prec)$ be a polynomial algebra of solvable type
  (for an arbitrary semigroup order\/ $\prec$) such that\/ $\kk[X']$ is a
  subring of\/ $\P$ for any subset\/ $X'\subset X$.  Furthermore, let\/ $L$ be
  a constructive Noetherian division.  For a finite set\/ $\F\subset\P$ of
  polynomials let\/ $\I=\lspan{\F}$ be the left ideal generated by it in the
  localisation\/ $\locp$.  If we apply Algorithm \ref{alg:polycompl} with the
  involutive Mora normal form instead of the usual one to the set\/ $\F$, then
  it terminates with an involutive basis of the ideal\/ $\I$.
\end{theorem}

\begin{proof}
  The termination of Algorithm~\ref{alg:polycompl} under the made assumptions
  was shown in Proposition~\ref{prop:localinvpoly} and
  Theorem~\ref{thm:polycompl}.  One easily verifies that their proofs are not
  affected by the substitution of the normal form algorithm, as they rely
  mainly on Theorem~\ref{thm:invnormalrep} and on the fact that the leading
  term of the normal form is not involutively divisible by the leading term of
  any generator.  Both properties remain valid for the Mora normal form.\qed
\end{proof}

\begin{remark}
  Note that Theorem \ref{thm:mora} guarantees the existence of \emph{strong}
  involutive bases.  Due to the extension to $\locp$, Example
  \ref{ex:nostrongbas} is no longer a valid counter example.  As the first
  generator in $\F$ is now a unit, we find that $\lspan{\F}=\locp$ and
  $\{1\}$ is a trivial strong Pommaret basis.\bull
\end{remark}

\begin{example}
  We continue Example~\ref{ex:strjanbas}.  Following the approach given by
  Theorem~\ref{thm:mora}, we immediately compute as Janet basis of
  $\lspan{\F}$ (over $\locp$) the minimal basis
  $\H_3=\{\partial_x,\partial_y,\partial_z\}$.  Obviously, it is considerably
  smaller than the bases obtained with Lazard's approach (over $\P$).  This
  effect becomes even more profound, if we look at the sizes of the bases in
  the homogenised Weyl algebra: both $\tilde\H_1$ and $\tilde\H_2$ consist of
  $21$ generators.\bull
\end{example}

\section{Involutive Bases over Rings}\label{sec:ring}

Finally, we consider the general case that $\P=\R[x_1,\dots,x_n]$ is a
polynomial algebra of solvable type over a (left) Noetherian ring $\R$.  In
the commutative case, Gr\"obner bases for such algebras have been studied in
\cite{gtz:pridec,tr:groe} (see \cite[Chapt.~4]{al:gb} for a more extensive
textbook discussion); for {\small PBW} extensions (recall
Example~\ref{ex:pbw}) a theory of Gr\"obner bases was recently developed in
\cite{gry:pbw}.  We will follow the basic ideas developed in these references
and assume in the sequel that linear equations are solvable in the coefficient
ring $\R$, meaning that the following two operations can be effectively
performed:
\begin{description}
\item[\phantom{i}{\upshape (i)}] given elements $s,r_1,\dots,r_k\in\R$, we can
  decide whether $s\in\lspan{r_1,\dots,r_k}_{\R}$ (the left ideal in $\R$
  generated by $r_1,\dots,r_k$);
\item[{\upshape (ii)}] given elements $r_1,\dots,r_k\in\R$, we can construct a
  finite basis of the module $\syz{}{r_1,\dots,r_k}$ of left syzygies
  $s_1r_1+\cdots+s_kr_k=0$.
\end{description}

The first operation is obviously necessary for the algorithmic reduction of
polynomials with respect to a set $\F\subset\P$.  The necessity of the second
operation will become evident later.  Compared with the commutative case,
reduction is a more complicated process, in particular due to the possibility
that in the commutation relations (\ref{eq:crmono}) for the multiplication in
$\P$ the maps $\rho_\mu$ may be different from the identity on $\R$ and the
coefficients $r_{\mu\nu}$ unequal one.

Let $\G\subset\P$ be a finite set.  We introduce for any polynomial $f\in\P$
the sets $\G_f=\{g\in\G\mid\ \le{\prec}{g}\mid\le{\prec}{f}\}$ and
\begin{equation}
  \bar\G_f=\bigl\{x^\mu\star g\mid g\in\G_f\wedge
                  \mu=\le{\prec}{f}-\le{\prec}{g}\wedge
                  \le{\prec}{(x^\mu\star g)}=\le{\prec}{f}\bigr\}
\end{equation}
Note that the last condition in the definition of $\bar\G_f$ is redundant
only, if the coefficient ring $\R$ is an integral domain.  Otherwise it may
happen that $|\bar\G_f|<|\G_f|$, namely if $\rho_\mu(r)r_{\mu\nu}=0$ where
$\lm{\prec}{g}=rx^\nu$.  The polynomial $f$ is \emph{head reducible} with
respect to $\G$, if $\lc{\prec}{g}\in\lspan{\lc{\prec}{\bar\G_f}}_{\R}$ (note
that we use $\bar\G_f$ here so that the reduction comes only from the leading
terms and is not due to some zero divisors as leading coefficients).
\emph{Involutive head reducibility} is defined analogously via sets $\G_{f,L}$
and $\bar\G_{f,L}$ where only involutive divisors with respect to the division
$L$ on $\Nn$ are taken into account, i.\,e.\ we set
\begin{equation}
  \G_{f,L}=\{g\in\G\mid
             \le{\prec}{f}\in\cone{L,\le{\prec}{\G}}{\le{\prec}{g}}\}\;.  
\end{equation}
Thus the set $\G$ is \emph{involutively head autoreduced}, if
$\lc{\prec}{g}\notin\lspan{\lc{\prec}{(\bar\G_{g,L}\setminus\{g\})}}_{\R}$ for
all polynomials $g\in\G$.  This is now a much weaker notion than before; in
particular, Lemma~\ref{lem:head} is no longer valid.

\begin{definition}\label{def:gbring}
  Let\/ $\I\subseteq\P$ be a left ideal in the polynomial algebra\/
  $(\P,\star,\prec)$ of solvable type over a ring\/ $\R$ in which linear
  equations can be solved.  A finite set\/ $\G\subset\P$ is a\/
  \emph{Gr\"obner basis} of\/ $\I$, if for every polynomial $f\in\I$ the
  condition\/ $\lc{\prec}{f}\in\lspan{\lc{\prec}{\bar\G_f}}_{\R}$ is
  satisfied.  The set\/ $\G$ is a\/ \emph{weak involutive basis} for the
  involutive division\/ $L$, if for every polynomial $f\in\I$ the condition\/
  $\lc{\prec}{f}\in\lspan{\lc{\prec}{\bar\G_{f,L}}}_{\R}$ is satisfied.  A
  weak involutive basis is a\/ \emph{strong involutive basis}, if every set\/
  $\bar\G_{f,L}$ contains precisely one element.
\end{definition}

It is easy to see that the characterisation of (weak) involutive bases via the
existence of involutive standard representations
(Theorem~\ref{thm:invnormalrep}) remains valid.  Indeed, only the first part
of the proof requires a minor change: the polynomial $f_1$ is now of the form
$f_1=f-\sum_{h\in\H_{f,L}}r_hh$ where the coefficients $r_h\in\R$ are chosen
such that $\le{\prec}{f_1}\prec\le{\prec}{f}$.  

Clearly, a necessary condition for the existence of Gr\"obner and thus of
(weak) involutive bases for arbitrary left ideals $\I\subset\P$ is that the
algebra $\P$ is a (left) Noetherian ring.  As we have seen in
Section~\ref{sec:hilbert}, this assumption becomes non-trivial, if the
coefficient ring $\R$ is not a field.  In this section, we will assume
throughout that $\P$ is a polynomial algebra of solvable type over a left
Noetherian ring $\R$ with centred commutation relations (cf.\
Definition~\ref{def:ccr}) so that Theorem \ref{thm:subfieldnoether} asserts
that $\P$ is left Noetherian, too.\footnote{The case of an iterated polynomial
  algebra of solvable type (cf.\ Definition~\ref{def:iteralg}) will be
  considered in Part II, after we have developed a syzygy theory for
  involutive bases.} A very useful side effect of this assumption is that the
scalars appearing in the commutation relations (\ref{eq:crmono}) are units
and thus not zero divisors which is important for some arguments.

\begin{example}\label{ex:xyz}
  As in the previous two sections, we cannot generally expect strong
  involutive bases to exist.  As a simple concrete example, also demonstrating
  the need of the second assumption on $\R$, we consider in $\kk[x,y][z]$
  (with the ordinary multiplication) the ideal $\I$ generated by the set
  $\F=\{x^2z-1,y^2z+1\}$.  Obviously, both generators have the same leading
  exponent $[1]$; nevertheless none is reducible by the other one due to the
  relative primeness of the coefficients.  Furthermore, the syzygy
  $\Sv=x^2\ev_2-y^2\ev_1\in\kk[x,y]^2$ connecting the leading coefficients
  leads to the polynomial $x^2+y^2\in\I$.  It is easy to see that a Gr\"obner
  and weak Janet basis of $\I$ is obtained by adding it to $\F$.  A strong
  Janet basis does not exist, as none of these generators may be removed from
  the basis.  \bull
\end{example}

This example shows that simply applying our completion Algorithm
\ref{alg:polycompl} will generally not suffice.  Obviously, with respect to
the Janet division $z$ is multiplicative for both elements of $\F$ so that no
non-multiplicative variables exist and thus it is not possible to generate the
missing generator by multiplication with a non-multiplicative variable.  We
must substitute in Algorithm \ref{alg:polycompl} the involutive head
autoreduction by a more comprehensive operation.\footnote{In the classical
  case of commutative variables over a coefficient field, it is not difficult
  to show that for any finite set $\F$ the syzygy module
  $\syz{}{\lm{\prec}{\F}}$ of the leading \emph{monomials} can be spanned by
  binomial generators corresponding to the $S$-polynomials in the Buchberger
  algorithm.  In Part II we will show that in any such syzygy at least one
  component contains a non-multiplicative variable, so that implicitly the
  involutive completion algorithm also runs over a generating set of this
  syzygy module.  When we move on to coefficient rings, it is well-known that
  additional, more complicated syzygies coming from the coefficients must be
  considered.  For these we can no longer assume that one component contains a
  non-multiplicative variable.  Hence \emph{partially} we must follow the same
  approach as in the generalisation of the Buchberger algorithm and this leads
  to the notion of $\R$-saturation where some syzygies not reachable via
  non-multiplicative variables are explicitly considered.}

\begin{definition}\label{def:rsat}
  Let\/ $\F\subset\P$ be a finite set and\/ $L$ an involutive division.  We
  consider for each\/ $f\in\F$ the syzygies\/ $\sum_{\bar
    f\in\bar\F_{f,L}}s_{\bar f}\,\lc{\prec}{\bar f}=0$ connecting the leading
  coefficients of the elements of the set\/ $\bar\F_{f,L}$.  The set\/ $\F$
  is\/ \emph{involutively $\R$-saturated} for the division\/ $L$, if for any
  such syzygy\/ $\Sv$ the polynomial $\sum_{\bar f\in\bar\F_{f,L}}s_{\bar
    f}\bar f$ possesses an involutive standard representation with respect
  to\/ $\F$.
\end{definition}

For checking involutive $\R$-saturation, it obviously suffices to consider a
finite basis of each of the finitely many syzygy modules
$\syz{}{\lc{\prec}{\bar\F_{f,L}}}$ so that such a check can easily be
performed effectively.  An element $f\in\F$ is involutively head reducible by
the other elements of $\F$, if and only if $\syz{}{\lc{\prec}{\bar\F_{f,L}}}$
contains a syzygy with $s_{f}=1$.  For this reason it is easy to combine an
involutive $\R$-saturation with an involutive head autoreduction leading to
Algorithm~\ref{alg:rsat}.

The \texttt{for} loop in Lines /5-13/ takes care of the involutive head
autoreduction (the call $\mathtt{HeadReduce}_{L,\prec}(f,\H)$ involutively
head reduces $f$ with respect to the set $\H\setminus\{f\}$ but with
multiplicative variables determined with respect to the full set
$\H$---cf.~Remark \ref{rem:auto}).  The \texttt{for} loop in Lines /17-22/
checks the involutive $\R$-saturation.  Each iteration of the outer
\texttt{while} loop analyses from the remaining polynomials (collected in
$\S$) those with the highest leading exponent.  The set $\S$ is reset to the
full basis, whenever a new element has been put into $\H$; this ensures that
all new reduction possibilities are taken into account.  In Line /15/ it does
not matter which element $f\in\S_\nu$ we choose, as the set $\H'_{f,L}$
depends only on $\le{\prec}{f}$ and all elements of $\S_\nu$ possess by
construction the same leading exponent $\nu$.

\begin{algorithm}
  \caption{Involutive $\R$-saturation (and head autoreduction)\label{alg:rsat}}
  \begin{algorithmic}[1]
    \REQUIRE finite set $\F\subset\P$, involutive division $L$ on $\Nn$
    \ENSURE involutively $\R$-saturated and head autoreduced set $\H$ with
            $\lspan{\H}=\lspan{\F}$
    \STATE $\H\leftarrow\F$;\quad $\S\leftarrow\F$
    \WHILE{$S\neq\emptyset$}
        \STATE $\nu\leftarrow\max_\prec\le{\prec}{\S}$;\quad 
               $\S_\nu\leftarrow\{f\in\H\mid\le{\prec}{f}=\nu\}$
        \STATE $\S\leftarrow\S\setminus\S_\nu$;\quad $\H'\leftarrow\H$
        \FORALL{$f\in\S_\nu$}
            \STATE $h\leftarrow\mathtt{HeadReduce}_{L,\prec}(f,\H)$
            \IF{$f\neq h$}
                \STATE $\S_\nu\leftarrow\S_\nu\setminus\{f\}$;\quad
                       $\H'\leftarrow\H'\setminus\{f\}$
                \IF{$h\neq0$}
                    \STATE $\H'\leftarrow\H'\cup\{h\}$
                \ENDIF
            \ENDIF
        \ENDFOR
        \IF{$\S_\nu\neq\emptyset$}
            \STATE choose $f\in\S_\nu$ and determine the set $\bar\H'_{f,L}$
            \STATE compute basis $\B$ of $\syz{}{\lc{\prec}{\bar\H'_{f,L}}}$
            \FORALL{$\Sv=\sum_{\bar f\in\bar\H'_{f,L}}s_{\bar f}
                              \ev_{\bar f}\in\B$}
                \STATE $h\leftarrow\mathtt{NormalForm}_{L,\prec}
                    (\sum_{\bar f\in\bar\H'_{f,L}}s_{\bar f}\bar f,\H')$
                \IF{$h\neq0$}
                    \STATE $\H'\leftarrow\H'\cup\{h\}$
                \ENDIF
            \ENDFOR
        \ENDIF
        \IF{$\H'\neq\H$}
            \STATE $\H\leftarrow\H'$;\quad $\S\leftarrow\H$
        \ENDIF
    \ENDWHILE
    \STATE \algorithmicreturn{$\H$}
  \end{algorithmic}
\end{algorithm}

\begin{proposition}\label{prop:rsat}
  Under the made assumptions about the polynomial algebra\/ $\P$,
  Algorithm~\ref{alg:rsat} terminates for any finite input set\/ $\F\subset\P$
  with an involutively\/ $\R$-saturated and head autoreduced set\/ $\H$ such
  that\/ $\lspan{\H}=\lspan{\F}$.
\end{proposition}

\begin{proof}
  The correctness of the algorithm is trivial.  The termination follows from
  the fact that both $\R$ and $\Nn$ are Noetherian.  Whenever we add a new
  polynomial $h$ to the set $\H'$, we have either that
  $\le{\prec}{h}\notin\lspan{\le{\prec}{\H'}}_{\Nn}$ or
  $\lc{\prec}{h}\notin\lspan{\lc{\prec}{\H'_{h,L}}}_{\R}$.  As neither in
  $\Nn$ nor in $\R$ infinite ascending chains of ideals are possible, the
  algorithm must terminate after a finite number of steps.\qed
\end{proof}

An obvious idea is now to substitute in the completion
Algorithm~\ref{alg:polycompl} the involutive head autoreduction by an
involutive $\R$-saturation.  Recall that Proposition~\ref{prop:localinvpoly}
(and Corollary~\ref{cor:localinvpoly}) was the crucial step for proving the
correctness of Algorithm~\ref{alg:polycompl}.  Our next goal is thus to show
that under the made assumptions for involutively $\R$-saturated sets local
involution implies weak involution.

\begin{proposition}\label{prop:localinvring}
  Under the made assumptions about the polynomial algebra\/ $\P$, a finite,
  involutively\/ $\R$-saturated set\/ $\F\subset\P$ is weakly involutive, if
  and only if it is locally involutive.
\end{proposition}

\begin{proof}
  We first note that Proposition~\ref{prop:localinvpoly} remains true under
  the made assumptions.  Its proof only requires a few trivial modifications,
  as all appearing coefficients (for example, when we rewrite
  $x^\mu\rightarrow x^{\mu-1_j}\star x_j$) are units in the case of centred
  commutation relations and thus we may proceed as for a field.  Hence if $\F$
  is locally involutive, then $\I=\lspan{\F}=\ispan{\F}{L,\prec}$ implying
  that any polynomial $g\in\I$ may be written in the form
  $g=\sum_{f\in\F}P_f\star f$ with $P_f\in\R[\mult{X}{L,\F,\prec}{f}]$.
  Furthermore, it follows from this proof that for centred commutation
  relations we may assume that the polynomials $P_f$ satisfy
  $\le{\prec}{(P_f\star f)}=\le{\prec}{P_f}+\le{\prec}{f}$.  We are done, if
  we can show that they can be chosen such that additionally
  $\le{\prec}{(P_f\star f)}\preceq\le{\prec}{g}$, i.\,e.\ such that we obtain
  an involutive standard representation of $g$.
  
  If the representation coming out of the proof of
  Proposition~\ref{prop:localinvpoly} already satisfies this condition on the
  leading exponents, nothing has to be done.  Otherwise we set
  $\nu=\max_\prec\bigl\{\le{\prec}{(P_f\star f)}\mid f\in\F\bigr\}$ and
  $\F_\nu=\{f\in\F\mid\le{\prec}{(P_f\star f)}=\nu\}$.  As by construction
  $\nu\in\bigcap_{f\in\F_\nu}\cone{L,\le{\prec}{\F}}{\le{\prec}{f}}$, the
  properties of an involutive division imply that we can write
  $\F_\nu=\{f_1,\dots,f_k\}$ with
  $\le{\prec}{f_1}\mid\le{\prec}{f_2}\mid\cdots\mid\le{\prec}{f_k}$ and hence
  $\F_\nu\subseteq\F_{f_k,L}$.  Since we have assumed that
  $\le{\prec}{(P_f\star f)}=\le{\prec}{P_f}+\le{\prec}{f}$, we even find
  $\F_\nu\subseteq\bar\F_{f_k,L}$.
  
  By construction, the equality $\sum_{f\in\F_\nu}\lc{\prec}{(P_f\star f)}=0$
  holds.  If we now set $\lm{\prec}{f}=r_fx^{\nu_f}$ and
  $\lm{\prec}{P_f}=s_fx^{\mu_f}$, then we obtain under the made assumptions:
  $\lc{\prec}{(P_f\star f)}= s_{f}\rho_{\mu_{f}}(r_{f})r_{\mu_{f}\nu_{f}}=
  \bigl[s_{f}\bar\rho_{\mu_{f}}(r_{f})r_{\mu_{f}\nu_{f}}\bigr]r_{f}$ and hence
  the above equality corresponds to a syzygy of the set
  $\lc{\prec}{\F_{f_k,L}}$.  As the set $\F$ is involutively $\R$-saturated,
  there exists an involutive standard representation
  \begin{equation}
    \sum_{i=1}^k
        \bigl[s_{f_{i}} \bar\rho_{\mu_{f_{i}}}(r_{f_{i}})
              r_{\mu_{f_{i}}\nu_{f_{i}}}\bigr] \bar f_i
    =\sum_{f\in\F}Q_f\star f
  \end{equation}
  with $Q_f\in\kk[\mult{X}{L,\F,\prec}{f}]$ and $\le{\prec}{(Q_f\star
    f)}=\le{\prec}{Q_f}+\le{\prec}{f}\prec\nu_{f_k}$.
  
  Introducing now the polynomials $Q'_f=Q_f-
  \bigl[s_{f}\bar\rho_{\mu_{f}}(r_{f}) r_{\mu_{f}\nu_{f}}\bigr]
  x^{\nu_{f_k}-\nu_f}$ for $f\in\F_\nu$ and $Q'_f=Q_f$ otherwise, we get the
  syzygy $\sum_{f\in\F}Q'_f\star f=0$.  If we set
  $P'_f=P_f-c_f^{-1}x^{\nu-\nu_{f_k}}\star Q'_f$ with
  $c_f=\bar\rho_{\nu-\nu_{f_k}}\bigl(s_{f}\bar\rho_{\mu_{f}}(r_{f})
  r_{\mu_{f}\nu_{f}}\bigr) \bar\rho_{\mu_{f}}(r_{f})r_{\mu_{f}\nu_{f}}$, then,
  by construction, $g=\sum_{f\in\F}P'_f\star f$ is another involutive
  representation of the polynomial $g$ with
  $\nu'=\max_\prec\bigl\{\le{\prec}{(P'_f\star f)}\mid f\in\F\bigr\}\prec\nu$.

  Repeating this procedure for a finite number of times obviously yields an
  involutive standard representation of the polynomial $g$.  As $g$ was an
  arbitrary element of the ideal $\I=\lspan{\F}$, this implies that $\F$ is
  indeed weakly involutive.\qed
\end{proof}

\begin{theorem}\label{thm:rsat}
  Let\/ $\P$ be a polynomial algebra of solvable type satisfying the made
  assumptions.  If the subalgorithm\/ $\mathtt{InvHeadAutoReduce}_{L,\prec}$
  is substituted in Algorithm~\ref{alg:polycompl} by Algorithm~\ref{alg:rsat},
  then the completion will terminate with a weak involutive basis of\/
  $\I=\lspan{\F}$ for any finite input set\/ $\F\subset\P$ such that the
  monoid ideal\/ $\le{\prec}{\I}$ possesses a weak involutive basis.
\end{theorem}

\begin{proof}
  The correctness of the modified algorithm follows immediately from
  Proposition~\ref{prop:localinvring}.  For the termination we may use the
  same argument as in the proof of Theorem~\ref{thm:polycompl}, as it depends
  only on the leading exponents.\qed
\end{proof}

\section{Conclusions}\label{sec:conc}

We studied involutive bases for a rather general class of non-commutative
polynomial algebras.  Our approach was very closely modelled on that of
Kandry-Rody and Weispfenning \cite{krw:ncgb} and subsequently Kredel
\cite{hk:solvpoly}.  However, we believe that the third condition in
Definition \ref{def:solvalg} (the compatibility between term order $\prec$ and
non-commutative product $\star$) is more natural than the stricter axioms in
\cite{krw:ncgb}.  In particular, we could not see where Kandry-Rody and
Weispfenning needed their stricter conditions, as all their main results hold
in our more general situation, as later shown by Kredel.

Comparing with \cite{apel:diss,bgv:algnoncomm,hk:solvpoly,vl:diss}, one must
say that the there used approach is more constructive than ours.  More
precisely, all these authors specify the non-commutative product via
commutation relations and thus have automatically a concrete algorithm for
evaluating any product.  As we have seen in the proof of Proposition
\ref{prop:fixmult}, the same data suffices to fix our axiomatically described
product, but it does not provide us with an algorithm.  However, we showed
that we can always map to their approach via a basis transformation.

We showed that the polynomial algebras of solvable type form a natural
framework for involutive bases.  This is not surprising, if one takes into
account that the main part of the involutive theory happens in the monoid
$\Nn$ and the decisive third condition in Definition~\ref{def:solvalg} of a
polynomial algebra of solvable type ensures that its product~$\star$ does not
interfere with the leading exponents.

We extended the theory of involutive bases to semigroup orders and to
polynomials over coefficient rings.  It turned out that the novel concept of a
\emph{weak} involutive basis is crucial for such generalisations, as in both
cases strong bases rarely exist.  These weak bases are still Gr\"obner bases
and involutive standard representations still exist (though they are no longer
unique).  It seems that in such computations the Janet division has a
distinguished position, as by Theorem \ref{thm:strjanbas} strong Janet bases
always exist.  If one is only interested in using
Algorithm~\ref{alg:polycompl} as an alternative to Buchberger's algorithm,
weak bases are sufficient.  However, most of the more advanced applications of
involutive bases studied in Part II will require strong involutive bases.

Concerning involutive bases over rings, we will study in Part II the special
case that the coefficient ring is again a polynomial algebra of solvable type.
Using the syzygy theory that will be developed there, we will be able to
obtain stronger results and a ``purely involutive'' completion algorithm.  The
current approach contains hidden in the concept of $\R$-saturation parts of
the Buchberger algorithm for the construction of Gr\"obner bases over rings.

Definition \ref{def:invdiv} represents the currently mainly used definition of
an involutive division.  While it appears quite natural, one problem is that
in some sense too many involutive divisions exist, in particular rather weird
ones with unpleasant properties.  This effect has lead us to the introduction
of such technical concepts like continuity and constructivity.  One could
imagine that there should exist a stricter definition of involutive divisions
that automatically ensures that Algorithm \ref{alg:monocompl} terminates
without having to resort to these technicalities.

Most of these weird divisions are globally defined and multiplicative indices
are assigned only to finitely many multi indices.  Such divisions are
obviously of no interest, as more or less no monoid ideal possesses an
involutive basis for them.  One way to eliminate these divisions would be to
require that for every $q\in\NN_0$ the monoid ideal $(\Nn)_{\geq
  q}=\{\nu\in\Nn\mid q\leq|\nu|\}$ has an involutive basis.  All the
involutive divisions used in practice satisfy this condition, but it is still
a long way from this simple condition to the termination of Algorithm
\ref{alg:monocompl}.

We did not discuss the efficiency of the here presented algorithms.  Much of
the literature on involutive bases is concerned with their use as an
alternative approach to the construction of Gr\"obner bases.  In particular,
recent experiments by Gerdt et al.\ \cite{gby:janbas2} comparing a specialised
\texttt{C/C++} program for the construction of Janet bases with the Gr\"obner
bases package of \textsc{Singular} \cite{gps:singular} indicate that the
involutive approach is highly competitive.  This is quite remarkable, if one
takes into account that \textsc{Singular} is based on the results of many
years of intensive research on Gr\"obner bases by many groups, whereas
involutive bases are still very young and only a few researchers have actively
worked on them.  The results in Part II will offer some heuristic explanations
for this observation.

Finally, we mention that most of the algorithms discussed in this article have
been implemented (for general polynomial algebras of solvable type) by
M.~Hausdorf \cite{wms:invbas1,wms:invbas2} in the computer algebra system
\MuPAD.\footnote{For more information see \texttt{www.mupad.de}.} The
implementation does not use the simple completion
Algorithm~\ref{alg:polycompl} but a more optimised version yielding minimal
bases developed by Gerdt and Blinkov \cite{gb:minbas}.  It also includes the
modified algorithm for determining strong Janet bases in local rings.

\appendix
\section{Term Orders}\label{sec:termord}

We use in this article non-standard definitions of some basic term orders.
More precisely, we revert the order of the variables: our definitions become
the standard ones, if one transforms
$(x_1,\dots,x_n)\rightarrow(x_n,\dots,x_1)$.  The reason for this reversal is
that this way the definitions fit better to the conventions in the theory of
involutive systems of differential equations.  Furthermore, they appear more
natural in some applications like the determination of the depth in Part II.

A \emph{term order} $\prec$ is for us a total order on the set $\TT$ of all
terms $x^\mu$ satisfying the following two conditions: (i) $1\preceq t$ for
all terms $t\in\TT$ and (ii) $s\preceq t$ implies $r\cdot s\preceq r\cdot t$
for all terms $r,s,t\in\TT$.  If a term order fulfils in addition the
condition that $s\prec t$ whenever $\deg s<\deg t$, it is called \emph{degree
  compatible}.  As $\TT$ and $\Nn$ are isomorphic as monoids, we may also
speak of term orders on $\Nn$.  In fact, most term orders are defined via
multi indices.

A more appropriate name for term orders might be \emph{monoid orders}, as the
two conditions above say nothing but that these orders respect the monoid
structure of $\TT$.  A more general class of (total) orders are
\emph{semigroup orders} where we skip the first condition, i.\,e.\ we only
take the semigroup structure of $\TT$ into account.  It is a well-known
property of such orders that they are no longer well-orders.  This implies in
particular the existence of infinite descending sequences so that normal form
algorithms do not necessarily terminate.

The \emph{lexicographic} order is defined by $x^\mu\prec_{\mbox{\scriptsize
    lex}}x^\nu$, if the last non-vanishing entry of $\mu-\nu$ is negative.
Thus $x_2^2x_3 \prec_{\mbox{\scriptsize lex}} x_1x_3^2$.  With respect to the
\emph{reverse lexicographic} order $x^\mu\prec_{\mbox{\scriptsize
    revlex}}x^\nu$, if the first non-vanishing entry of $\mu-\nu$ is positive.
Now we have $x_1x_3^2 \prec_{\mbox{\scriptsize revlex}} x_2^2x_3$.  However,
$\prec_{\mbox{\scriptsize revlex}}$ is only a semigroup order, as it violates
the first condition: $x_1\prec_{\mbox{\scriptsize revlex}}1$.  Degree
compatible versions of these orders exist, too.
$x^\mu\prec_{\mbox{\scriptsize deglex}}x^\nu$, if $|\mu|<|\nu|$ or if
$|\mu|=|\nu|$ and $x^\mu\prec_{\mbox{\scriptsize lex}}x^\nu$.  Similarly,
$x^\mu\prec_{\mbox{\scriptsize degrevlex}}x^\nu$, if $|\mu|<|\nu|$ or if
$|\mu|=|\nu|$ and $x^\mu\prec_{\mbox{\scriptsize revlex}}x^\nu$.  Obviously
$\prec_{\mbox{\scriptsize degrevlex}}$ is a term order.  It possesses the
following useful characterisation which is easy to prove.
  
\begin{lemma}\label{lem:revlex}
  Let\/ $\prec$ be a degree compatible term order such that the condition\/
  $\lt{\prec}{f}\in\lspan{x_1,\dots,x_k}$ is equivalent to\/
  $f\in\lspan{x_1,\dots,x_k}$ for every homogeneous polynomial\/ $f\in\P$.
  Then $\prec$ is the degree reverse lexicographic order
  $\prec_{\mbox{\scriptsize\upshape degrevlex}}$.
\end{lemma}

We say that a term order \emph{respects classes}, if for multi indices $\mu$,
$\nu$ of the same length $\cls{\mu}<\cls{\nu}$ implies $x^\mu\prec x^\nu$.  It
is now easy to see that by Lemma~\ref{lem:revlex} on terms of the same degree
any class respecting term order on $\TT$ coincides with the degree reverse
lexicographic order.  If we consider free polynomial modules, class respecting
orders have the same relation to {\small TOP} lifts \cite{al:gb} of
$\prec_{\mbox{\scriptsize degrevlex}}$.

\begin{acknowledgement}
  The author would like to thank V.P.~Gerdt for a number of interesting
  discussions on involutive bases.  M.~Hausdorf and R.~Steinwandt participated
  in an informal seminar at Karlsruhe University where most ideas of this
  article were presented and gave many valuable comments.  The constructive
  remarks of the anonymous referees were also very helpful.  This work
  received partial financial support by Deutsche Forschungsgemeinschaft, INTAS
  grant 99-1222 and NEST-Adventure contract 5006 (\emph{GIFT}).
\end{acknowledgement}

\bibliography{../../BIB/DiffEq,../../BIB/Algebra,../../BIB/Seiler,%
../../BIB/Groebner,../../BIB/DiffAlg,../../BIB/Misc,../../BIB/Lie,%
../../BIB/CA,../../BIB/ZZProc}
\bibliographystyle{plain}

\end{document}